\documentclass[12pt]{amsart}

\usepackage{amsmath}
\usepackage{amsthm}

\usepackage{amsfonts}
\usepackage{amssymb}
\usepackage{graphicx}
\usepackage[pagewise]{lineno}
\newcommand\cyr{%
\renewcommand\rmdefault{wncyr}%
\renewcommand\sfdefault{wncyss}%
\renewcommand\encodingdefault{OT2}%
\normalfont
\selectfont}
\DeclareTextFontCommand{\textcyr}{\cyr}

\usepackage{tikz-cd}
\usepackage{ amssymb }
\usepackage{tikz}
\usepackage{stix}
\usepackage{enumitem}
\usetikzlibrary{matrix,arrows}\usepackage[OT2,OT1]{fontenc}

\newcommand{\G}{\mathcal{G}}

\newcommand{\Z}{\mathbb{Z}}

\newcommand{\noi}{\noindent}

\newcommand{\bs}{\bigskip}

\newcommand{\ZZ}{\mathbb{Z}/3\mathbb{Z}}
\newcommand{\M}{\mathcal{M}}
\newcommand{\HH}{\mathcal H}

\input xy
\xyoption{all}

\usepackage[top =.6in, bottom=.6in, left=.5in, right=.7in]{geometry}

\input xy
\xyoption{all}
 
 \LARGE

\title{Cohomlogical Kernels of Elementary Abelian Degree $p^2$ Extensions}

\author{ Bill Jacob and Nathan Schley} 
\curraddr{Department of Mathematics, University of California, Santa Barbara, 
 Santa Barbara, California, USA 93106 }
\email{jacob@math.ucsb.edu, schley@math.ucsb.edu}

\begin{document}

\begin{abstract}
Let $p$ be an odd prime,
$F$  a field
with a primitive $p^2$th root of unity, and
$E=F(\sqrt[p]{b_1},\sqrt[p]{b_2})$  an elementary abelian extension of degree $p^2$. 
This paper studies the cohomological
kernel $H^n(E/F,\Z/p\Z):={\rm ker}\left(H^n(F,\Z/p\Z)\rightarrow H^n(E,\Z/p\Z)\right)$
for all $n$.   When $p=3$, using tools of Positselski, a six-term
exact sequence is given that is analogous to the $p=2$ case. As an application the  quotient $H^n(E/F,\ZZ)/{\rm Dec}^n(E/F,\ZZ)$ where 
${\rm Dec}^n(E/F,\ZZ)$ is the ``expected kernel'' is described.  This quotient group is of interest because computations of Tignol 
that show when $n=2$
nontrivial elements give rise to indecomposible division algebras of
exponent $3$ and index $9$.
\end{abstract}

\maketitle

\bs

\noi {\bf Introduction.} Suppose that
$p>2$ is prime,
$F$ is a field
containing a primitive $p^2$th root of unity, and suppose
$E=F(\sqrt[p]{b_1},\sqrt[p]{b_2})$ is an elementary abelian extension of degree $p^2$. The relative Brauer group 
$H^2(E/F,\Z/p\Z):=$ ${\rm ker}\left(H^2(F,\Z/p\Z)\rightarrow H^2(E,\Z/p\Z)
\right)$
is of interest because computations of Tignol [T]  showed
it can be larger than the ``expected'' kernel 
${\rm Dec}^2(E/F):=(b_1)\smile H^1(F,\Z/p\Z)+(b_2)\smile H^1(F,\Z/p\Z)$
where these new elements give rise to indecomposible division algebras of
exponent $p$ and index $p^2$ (see  [T] or [TW] as well as [J], [AJ1],[AJO1],[AJO2],[AJO3] for related work.)
 In this paper we give a four-term sequence of Galois modules
 $$0\rightarrow\M_1\rightarrow\M_2\rightarrow\M_3\rightarrow\M_4\rightarrow0,$$
 each of which are free as $\Z$-modules and for which homotopies $h_i$
 giving multiplication by $p$ are constructed.  Passing to 
 the $\otimes_\Z\Z/p\Z$ and $\otimes_\Z\Z/p^2\Z$ versions
  a theory developed by Positselski [P] can be applied to give
a characterization of the 
kernel $H^n(E/F,\Z/p\Z)$.
Similar applications of  [P] can be found in
[AJ2], [AJ3], [AJO1],[AJO2] and [S].

When $p=3$ further results are obtained.  Using the same four-term sequence
one obtains a six term exact sequence of Galois cohomology groups that 
more thoroughly describe the cohomological kernels. These in turn
are used to characterize the Brauer Kernels as well as the 
quotient group $H^2(E/F,\ZZ)/{\rm Dec}^2(E/F,\ZZ)$ in this case. Elements of
$H^2(E/F,\ZZ)$ are explicitly described as $2$-cocyle classes where
the deviation of a class from ${\rm Dec}^2(E/F,\ZZ)$ is measured by
an invariant in $H^1(F,\ZZ)$.  We expect these results to be true for 
odd $p>3$ but more work is needed.  The complexity is the computation of
certain Bochstein maps that is possible when $p=3$ because of
the simpler structure of $\M_4$ in that case. 

The paper is organized as follows.
The first section gives a four-term exact sequence of $\G:=$ Gal($E/F)$-modules
with multiplication by $p$ homotopies for which the theory developed
by Positselski [P] can be applied.  In the second section, the structure of the
fourth module $\M_4$ is developed in detail which enables a characterization
of the cohomological kernels of interest. The results in the first two sections are given for arbitrary odd
primes $p$ with applications to be explored in a subsequent paper. In the third,
fourth and fifth  sections
the results are specialized to $p=3$ where the description of $H^n(E/F,\ZZ)$
is given.  This is rather technical and requires a computational analysis
of $H^{n-1}(F,\M_4\otimes_\Z\ZZ)$.
 
 In order to apply the Positelski theory, it is necessary that the
 Bochstein maps associated with the modules $\M_i$ are zero.
 Although this is clear for $\M_1$, $\M_2$ and $\M_3$ it is not true
 for $\M_4$.  It turns out that the vanishing of the Bochstein for the first
 three modules is sufficient for the applications in the first three sections.
 The fourth section shows that when $p=3$, although the fourth Bochsterin ${\mathcal{B}}_4$ is nonzero, the composite $h_3\circ{\mathcal{B}}_4$ is zero and this
 is enough for the Positselski machinery to work (see [S, Rem. 2.5].)  In this way the
 six-term sequences of Theorems 4.7 and 4.8 are obtained.  The final section
 gives further analysis of the quotient $H^n(E/F,\ZZ)/{\rm Dec}^n(E/F,\ZZ)$ making
 use of the results from Section 4.

 The paper uses essentially no field theory and consists of computations within
 the category of pro-$p$ groups that satisfy certain conditions known to be true
 for the pro-$p$ Galois $\G$ of a field $F$  containing a primitive $p^2$th root of unity. The conditions assumed are non trivial; in particular coinduced modules corresponding to closed subgroups are assumed to
 satisfy a Bockstein triviality condition, a result true for pro-$p$ Galois groups of fields
 and is a consequence the proof of the Bloch-Kato conjecture,
 which as $\mu_{p^2}\subset F$
gives an isomorphism $K^M_nF/p^rK^M_nF\rightarrow H^n(\G,\Z/p^r\Z)$ 
for $r=1,2$ (see [HW, Chap. 1] for more discussion.)
 Specifically we assume for $\G$ as well as its closed subgroups that
  the reduction (mod $p$) maps
 $H^n(\G,\Z/p^2\Z)\rightarrow H^n(\G,\Z/p\Z)$  are surjective, which for fields is evident 
 in view of the Bloch-Kato Conjecture and the surjectivity
 of $K^M_nF/p^2K^M_nF\rightarrow K^M_nF/pK^M_nF$.   
 We also will assume, corresponding to the fact that $\mu_{p^2}\subset F$ enables 
 cyclic $p$-extensions of $F$ extend to cyclic $p^2$-extenisons, that every
 normal subgroup $\HH$ of $\G$ with $\G/\HH\cong\Z/p\Z$ contains a 
 normal subgroup $\widetilde{\mathcal H}$ of $\G$ with 
 $\G/\widetilde{\mathcal H}\cong\Z/p^2\Z$. If a pro-$p$ group  $\G$ satisfies these
 conditions we shall call it a {\em Voevodsky} group.
 With these assumptions, most of the discussion in the paper deals with group cohomology, although
 regular digressions will connect the results to the field case.
 
 \bs
 
 \bs
 
 \noi \S1. {\bf The Four-term Sequence with Homotopies.}

\bs

We assume $p$ is an odd prime and $\G$ is a Voevodsky pro-$p$ 
group as defined above. In this section we define the $\G$-modules used in the
four-term sequence. Initially we work with free $\Z$-modules and from the results
there we obtain the necessary $\Z/p\Z$- and $\Z/p^2\Z$-modules needed to apply
the Positeselski machinery. We lay out some basic notation next.

\bs

\noi {\bf Notation 1.1.}
We set $G=\langle\tau_1,\tau_2\rangle\cong\Z/p\Z\oplus\Z/p\Z$
and set  $t_i:=\tau_i-1\in \Z[G]$ for $i=1,2$.  We define the ``trace'' elements, $T_1=1+\tau_1+\tau_1^2+\cdots+\tau_1^{p-1}$, $T_2=1+\tau_2+\tau_2^2+\cdots+\tau_2^{p-1}$, and for $3\leq i\leq p+1$ we set  $\tau_k:=\tau_1\tau_2^{k-2}$, $t_k=\tau_k-1$, and
$T_k=\sum_{j=0}^{p-1}(\tau_1\tau_2^{k-2})^j$. We denote
$T_G=\sum_{g\in G}g$.
  We fix a surjection
  $\G\rightarrow G$ through which we can view all $\Z[G]$-modules as
  $\G$-modules.

\bs

Our ``mindset'' is to think of $G$ as the Galois group
of  $E/F$ where $E:=F(\sqrt[p]{b_1},\sqrt[p]{b_2})$ and $G$ acts with
$\tau_i(\sqrt[p]{b_i})=\zeta\sqrt[p]{b_i}$. We will be interested
in subfields $L_i:=F(\sqrt[p]{b_i})$ where
 for $j\neq i$ we have that
$T_j:E\rightarrow L_i$ is the trace or norm depending on if it is interpreted
additively or multiplicatively.  For the   calculations,
however, the field is not involved, and we use the 
Positselski theory assuming $\G$ is a Voevodsky group. 
We note that $T_i^2=pT_i$
and if $i\neq j$ then $T_iT_j=T_G$.  Also, $t_iT_i=0$ and
$t_i^{p-1}\equiv T_i$ (mod $p\Z[G]$). 

\bs

\noi {\bf The Four-Term Sequence and Homotopies.} 
We will use the following four-term 
sequence,
$$(*) \hspace{.5in}   0\longrightarrow \Z[G]\cdot T_G\stackrel{d_1}{\longrightarrow}
\Z[G]\oplus\Z[G]\cdot T_G\stackrel{d_2}{\longrightarrow}\hspace{.8in}
$$ $$\hspace{.7in}
\left(
\Z[G]\oplus\Z[G]\oplus
\Z[G]\cdot T_3\oplus
\Z[G]\cdot T_{p+1}\right)\stackrel{d_{3}}{\longrightarrow}
{\rm cok}(d_2)\longrightarrow 0.$$
Here, 
\begin{eqnarray*}
d_1(n\cdot T_G)&=&(n\cdot T_G,\,-pn\cdot T_G),\\
d_{2}(\eta,n\cdot T_G)&=&(-\eta\cdot  t_1,\,-\eta\cdot  t_2,\,\eta\cdot  T_3+n\cdot T_G,
\,\eta\cdot  T_{p+1}+n\cdot T_G),
\end{eqnarray*} 
and $d_{3}$ is the projection onto the cokernel. 
Each of the $d_i$ clearly are $G$-module maps. 
To make notation easier we will abbreviate the modules in this sequence
\begin{eqnarray*}
\M_1&:=&\Z[G]\cdot T_G,\\ \M_2&:=&\Z[G]\oplus\Z[G]\cdot T_G,\\
\M_3&:=&\Z[G]\oplus\Z[G]\oplus
\Z[G]\cdot T_3\oplus
\Z[G]\cdot T_{p+1},\ \  {\rm  and}\\ 
\M_4&:=&{\rm cok}(d_2).
\end{eqnarray*} 
We remark that (as checked below) the sequence is exact 
and each of the $\M_i$ are free
as $\Z$-modules of ranks $1$, $p^2+1$, $2p^2+p$ and $p^2+p$ respectively.   

For the homotopies we use
\begin{eqnarray*}
h_1(\xi,nT_G)&=&\xi\cdot T_G+(p-1)n\cdot T_G\\
h_2(\xi_1 ,\,\xi_2 ,\,\xi_3T_3,\,\xi_4T_{p+1})&=&(
\xi_1\cdot\kappa_1+\xi_2
\cdot\kappa_2
- \frac{p-1}{2} \xi_3T_3-\frac{p-1}{2}\xi_4T_{p+1},\,\xi_3 T_3\cdot T_G)\\
{h}_3([m])&=&p\cdot m-d_2\circ h_2(m).
\end{eqnarray*}
where $\kappa_1, \kappa_2\in\Z[G]$ solve the equation
$$t_1\kappa_1+t_2\kappa_2+\frac{p-1}{2}  T_3+\frac{p-1}{2}  T_{p+1}=T_G-p,$$
and the well-definition of
$h_3$ is checked below.  We first must show that the $\kappa_i$ exist, where we
remark that they are not unique but a pair can be selected and is then fixed 
throughout the rest of this work.

\bs

\noi {\bf Lemma 1.2.} {\em There exist $\kappa_1,\kappa_2\in\Z[G]$ such that
$$t_1\kappa_1+t_2\kappa_2+\frac{p-1}{2}  T_3+\frac{p-1}{2}  T_{p+1}=T_G-p.$$
Moreover, $t_1t_3T_{p+1}\kappa_1+t_2t_3T_{p+1}\kappa_2=-\frac{p(p+1)}{2}t_3T_{p+1}$,
$t_1T_2\kappa_1=T_G-pT_2$, and $T_1t_2\kappa_1=T_G-pT_1$. It follows that
\begin{eqnarray*}
h_2(t_1T_2,0,0,0)&\equiv&(T_G,0) \ \ ({\rm mod} \ p\cdot\M_2)\\
h_2(0,T_1t_2,0,0)&\equiv&(T_G,0) \ \ ({\rm mod} \ p\cdot\M_2) \\
h_2(t_1^iT_2,0,0,0)&\equiv&(0,0) \ \ ({\rm mod} \ p\cdot\M_2) \ \ {\rm for} \ i>1\\
h_2(0,T_1t_2^j,0,0)&\equiv&(0,0)  \ \ ({\rm mod} \ p\cdot\M_2) \ \ {\rm for} \ j>1\\
h_2(t_1t_3T_{p+1},t_2t_3T_{p+1},0,0)&\equiv&(0,0)  \ \ ({\rm mod} \ p\cdot\M_2).\\
\end{eqnarray*}
When $p=3$ the values of $\kappa_1$ and $\kappa_2$ can be taken as
$\kappa_1=-(\tau_1+2\tau_1^2)=t_1+t_1^2-3(t_1+1)^2$ and $\kappa_2=-(\tau_2+2\tau_2^2)
=t_2+t_2^2-3(t_2+1)^2$.
}

\bs

\noi {\bf Proof.}   By a direct calculation
$$f:=2p-(T_1+T_2+T_4+T_5+\cdots+T_p)+\frac{p-3}{2}T_3+\frac{p-3}{2}T_{p+1}\in {\mathcal{I}}$$
where ${\mathcal{I}}=\Z[G]t_1+\Z[G]t_2={\rm ker}(T_G\cdot)$.
For, adding the coefficients of $f$ we find that $2p-(p-1)p+(p-3)p=0.$  So there exist
$\kappa_1,\kappa_2$ with 
$f=(1-\tau_1)\kappa_1+(1-\tau_2)\kappa_2$.
Since
$\sum_{i=1}^{p+1}T_i=p+T_G$, rearranging gives the result. 

We now have 
$$t_1\kappa_1+t_2\kappa_2=T_G-\frac{p-1}{2}  T_3-\frac{p-1}{2}  T_{p+1}-p.$$
Multiplication by $t_3T_{p+1}$ gives 
$t_1t_3T_{p+1}\kappa_1+t_2t_3T_{p+1}\kappa_2=-\frac{p(p+1)}{2}t_3T_{p+1}$; multiplication by $T_2$ gives $t_1T_2\kappa_1=T_G-pT_2$; and multiplication by $T_1$ gives
 $T_1t_2\kappa_1=T_G-pT_1$.
The final statements follow by direct computation from these and 
the definition of $h_2$.  $\Box$


\bs

We now check the well-definition, exactness and  homotopy conditions for the
four-term sequence.

\bs

\noi {\bf Theorem 1.3.} {\em The four term sequence $(*)$ given
above with maps $d_i$ is exact and the maps $h_i$ satisfy the
homotopy condition with multiplication by $3$. The modules
$\M_1,\M_2,\M_3,\M_4$ are free as $\Z$-modules of ranks
$1,p^2+1,2p^2+p,p^2+p$ respectively.  The $d_i$ and $h_i$ are
$G$-module maps and $h_3$ is well-defined.}

\bs

\noi {\bf Proof.} That all the maps considered are $G$-module maps
follows as they are given by sums of multiplication by 
elements of $\Z[G]$. The rank and freeness assertions are clear except for
$\M_4$ which is proved in the next section. To check exactness, by definition $d_1$ is injective, $d_3$ is surjective and
${\rm ker(d_3})={\rm im}(d_2)$.  Now we suppose that
$(\eta,nT_G)\in {\rm ker}(d_2)$.  By the definition of $d_{2}$ we see that
$\eta(\tau_1-1)=\eta(\tau_2-1)=0$.  From this we know $\eta=mT_G$ for some $m$
(since the coefficients of $\tau_1^i\tau_2^j$ must be the same as $\tau_k\cdot\tau_1^i\tau_2^j$
for all $i,j\in\{0,1,\ldots,p-1\}$ and $k=1,2$.)
By the definition of $d_{2}$, since $T_kT_G=pT_G$ we also find for $3\leq k\leq p+1$ that $\eta T_k+n\cdot T_G=
(pm+n)T_G=0$.  It follows that $n=-pm$ and this gives that ${\rm ker}(d_2)=\Z\cdot T_G={\rm im}(d_1)$ so exactness follows. 

To check the homotopy conditions, 
we have $d_1(n\cdot T_G)=(n\cdot T_G,-pn\cdot T_G)$.
Therefore as the first homotopy
$h_1:\Z[G]\oplus\Z[G]\cdot T_G {\longrightarrow}\,\Z[G]\cdot T_G$  
is given by $h_1(\eta,nT_G)=\eta\cdot T_G+(p-1)nT_G$ we will have
$h_1\circ d_1(T_G)=h_1(T_G,-pT_G)=p^2T_G-p(p-1)T_G=p\cdot T_G$
as desired.  

As $$h_2(\eta_1 ,\eta_2 ,\eta_3T_3,\eta_4T_4)=(
\eta_1\cdot\kappa_1+\eta_2
\cdot\kappa_2
-\frac{p-1}{2}  \eta_3T_3-\frac{p-1}{2}  \eta_4T_4,\eta_3 T_3\cdot T_G),$$
and
$$-t_1\kappa_1-t_2\kappa_2-\frac{p-1}{2}  T_3
-\frac{p-1}{2}  T_{p+1}=p-T_G,$$
we have
$$h_2(d_2(1,0))=h_2(-t_1,-t_2,T_3,T_{p+1})=(p-T_G,pT_G).
$$
As
$d_1\circ h_1(1,0)=d_1(T_G)=(T_G,-pT_G)$ we find that
$$(d_1\circ h_1+ {h}_2\circ d_2)(1,0)=(T_G,-pT_G)+(p-T_G,pT_G)=
(p,0).$$
We also compute 
$$h_2(d_2(0,T_G))=h_2(0,0,T_G,T_G)=(-(p-1)T_G,p^2T_G).$$
As $d_1\circ h_1(0,T_G)=d_1((p-1)T_G)=((p-1)T_G,-p(p-1)T_G)$ we find
$$(d_1\circ h_1+ {h}_2\circ d_2)(0,T_G)=
(0,pT_G).$$
This shows that  $(d_1\circ h_1+ {h}_2\circ d_2)=p\cdot$ since
the result is true for the $G$-module generators of ${\mathcal{M}}_2$.

Finally to define $h_3$ we define
$\tilde{h}_3:\M_3\rightarrow \M_3$  for $m\in \M_3$ by
$\tilde{h}_3(m):=pm-d_2\circ h_2(m)$.  We know that
${\rm im}(d_2)\subset \M_3$ is generated as a $G$-module by $\Delta:=(-t_1,-t_2,T_3,T_{p+1})=d_2(1,0)\in \M_3$
and $\Delta_G=(0,0,T_G,T_G)=d_2(0,T_G)\in \M_3$.
By the preceding paragraph, as $T_G\cdot T_3=T_G\cdot T_{p+1}=pT_G$  we have
\begin{eqnarray*}\tilde{h}_3(\Delta)&=&p\Delta-d_2\circ h_2(\Delta)=
 p\Delta-d_2(p-T_G,pT_G)\\&=&p\Delta-(-pt_1,-pt_2,(p-T_G)T_3+pT_G,(p-T_G)T_{p+1}+pT_G)\\&=&
 p\Delta-p\Delta=0,
 \end{eqnarray*}
 and
\begin{eqnarray*}
\tilde{h}_3(\Delta_G)&=&p\Delta_G-d_2\circ h_2(\Delta_G)=
 p\Delta_G-d_2(-(p-1)T_G,p^2T_G)\\&=&
 p\Delta_G-(0,0,-p(p-1)T_G+p^2T_G,-p(p-1)T_G+p^2T_G)\\
 &=& p\Delta_G-p\Delta_G=0.
 \end{eqnarray*}
 This shows that $\tilde{h}_3$ vanishes on the $G$-module ${\rm im}(d_2)$ and therefore
 induces a $G$-map $h_3:{\rm cok}(d_2)=\M_4\rightarrow \M_3$. 
 Denoting by $[ \ ]$ the equivalence
 classes in ${\rm cok}(d_2)$
we have established the well-definition of $${h}_3([m])=pm-d_2\circ h_2(m)$$
claimed above.  From this, using $d_3\circ d_2=0$,
 $$(d_3\circ{h}_3)([m])=d_3(pm- d_2\circ h_2(m))
 =[pd_3(m)-d_3\circ d_2\circ h_2(m)]=
 p[m]$$ showing that $d_3\circ{h}_3=p\cdot$. Moreover, by
 definition
$$
(d_2\circ h_2+ {h}_3\circ d_3)(m)=d_2\circ h_2(m)+(pm
-d_2\circ h_2(m))=pm$$
so we find that $d_2\circ h_2+ {h}_3\circ d_3=p\cdot$ as well. We have   verified all the homotopy prism conditions proving the theorem. $\Box$

\bs

In the next section we verify that  $\M_4$ is free
as a $\Z$-module.  Moreover, we will show 
each $\M_i$ has a $\Z$-basis that subdivides into
two sets, one of which lies in the kernel of the $d_i$ and the other maps to
a subset of a basis for $\M_{i+1}$ and therefore taking $\otimes_\Z\Z/p\Z$ and  $\otimes_\Z\Z/p^2\Z$
give exact sequences of $\Z/p\Z$ and $\Z/p^2\Z$-modules
as required to apply the theory given in [P].  
Also,  $\G$   is a  Voevodsky group
and  each of the $\M_i$ where $i=1,2,3$
are direct sums of coinduced modules for $\G$. This means, as in [S], 
that the Bochstein maps associated with the short exact sequence
$$0\rightarrow\M_i\otimes_\Z\Z/p\Z\rightarrow\M_i\otimes_\Z\Z/p^2\Z
\rightarrow\M_i\otimes_\Z\Z/p\Z\rightarrow0,
$$ 
${\mathcal{B}}_i:H^n(\G,\M_i\otimes_\Z\Z/p\Z)\rightarrow
H^{n+1}(\G,\M_i\otimes_\Z\Z/p\Z)$, are zero for $i=1,2,3$.
  It is {\em not} the case
that this is true for $\G$ and $\M_4$ and  Section 4
will develop the tools we need to circumvent this.

\bs

  We now give the main application of the
four-term sequence.  The Positselski connecting map
$\eta:H^{n-1}(F,\M_4\otimes_\Z\Z/p\Z){\longrightarrow}
H^{n}(F,\Z/p\Z)$ is defined in [P] (see also [AJ 1,2], [AJO 1,2] and [S].)

\bs

\noi {\bf Theorem 1.4.} {\em Suppose that $F$ is a field containing a primitive $p^2$th 
root of unity and that $E=F(\sqrt[p]{b_1}\sqrt[p]{b_2})$
is a $\Z/p\Z\oplus\Z/p\Z$ extension of $F$.  Then the Positselski connecting
map $\eta$ fits into an exact sequence
$$H^{n-1}(F,\M_4\otimes_\Z\Z/p\Z)\stackrel{\eta}{\longrightarrow}
H^{n}(F,\Z/p\Z)\stackrel{i_{E/F}}{\longrightarrow}H^{n}(F,\Z/p\Z)$$
that determines the cohomological kernel of this $p$-$p$ extension.
}

\bs  

When $p=3$ we will show in Section 4 more is true, namely that the composite
map $h_3\circ{\mathcal{B}}:H^n(\G,\M_4\otimes_\Z\ZZ)
\rightarrow H^{n+1}(\G,\M_3\otimes_\Z\ZZ)$ is zero, which is sufficient
to apply Positselski's machinery to obtain a six-term exact
sequence (see Theorem 4.8.)
However, although we conjecture it is so, we do not know yet if this
composite vanishes for $p>3$.
The fact that  the composite suffices for the six-term sequence  can be found by careful inspection
of the proof given in [P, Th.6].

\bs
 
\bs

\noi {\bf \S 2. \ The structure of the module $\M_4$.}

\bs

This first lemma collects results about the modules 
$\M_1$, $\M_2$ and $M_3$ described above.  
The proofs are elementary  so only a brief discussion is provided to enable the reader
to follow the identifications used later. The final statement in part (vi)
is necessary to apply the Positselski machinery, and as noted above,
relies in the field case on the proof of the Bloch-Kato conjecture (that
in turn relies on Voevodsky's work), which corresponds to our nontrivial assumptions
about the group $\G$.

\bs

\noi {\bf Lemma 2.1.} {\em
(i) The module $\Z[G]$ is $\Z$-free of rank $p^2$.  As a $\G$-module 
it is the coinduced module ${\rm CoInd}^G_{\{1\}}(\Z)=Hom_{\Z}(G,\Z)$.\\
(ii) The modules $\Z[G]\cdot T_i$ are $\Z$-free and have
$\Z$-rank $p$. As a $\G$-module 
it is the coinduced module ${\rm CoInd}^G_{\langle\tau_i\rangle}(\Z)
=Hom_{\Z[\langle\tau_i\rangle]}(G,\Z)$.\\
(iii) Consequently the modules
$\M_1,\M_2,\M_3$ are free as $\Z$-modules of ranks
$1,p^2+1,2p^2+p,$ respectively.\\
(iv) The module $\Z[G]\cdot\{T_1,T_2\}$ is $\Z$-free of rank $2p-1$ and the
module $\Z[G]\cdot t_1t_2$ is $\Z$-free of rank $(p-1)^2$ fitting into a short exact
sequence
$$0\rightarrow\Z[G]\cdot\{T_1,T_2\}\longrightarrow\ZZ[G]
\stackrel{\cdot t_1t_2}{\longrightarrow}\ZZ[G]\cdot t_1t_2\rightarrow0.$$
(v) The results in (i), (i), (ii), and (iv)  apply to the modules obtained
tensoring  with $\otimes_\Z\Z/p\Z$ and $\otimes_\Z\Z/p^2\Z$.\\
(vi) Assuming $\G$ is a Voevodsky group, or is the pro-$p$ Galois group of a field $F$ containing
a primitive $p^2$th root of unity, for $i=1,2,3$  the Bochstein maps
${\mathcal B}_i:H^n(\G,\M_i\otimes_\Z\Z/p\Z)\rightarrow
H^{n+1}(\G,\M_i\otimes_\Z\Z/p\Z)$ are zero.
}

\bs

\noi {\bf Proof.} For (i) and (ii) the rank assertions are easy. The coinduced module structures are well-known, where $f=\sum_{g\in G}f_g\cdot g\in \Z[G]$ maps to the function
$\tilde{f}\in Hom_{\Z}(G,\Z)$ given by $\tilde{f}(g)=f_g$, and we note that
when $f\in \Z[G]T_i$ then as $f_g=f_{\tau_ig}$ we must also have
$\tilde{f}\in Hom_{\Z[\langle\tau_i\rangle]}(G,\Z)$.
The assertions in (iii) are clear by the definitions. 

For (iv),
$B:=\{t_1^it_2^j|0\leq i,j\leq p-1
\}$ is a $\Z$-basis for $\Z/p\Z[G]$ which we partition into 
$B_1:=\{t_1^it_2^{p-1}|0\leq i\leq p-1\}\cup$ $\{t_1^{p-1}t_2^j|0\leq j\leq p-2\}$ and
$B_2:=\{t_1^it_2^j|0\leq i,j\leq p-2
\}$.
We have $B_1\subset\ZZ[G]\cdot t_1t_2$ and $B_2\cdot t_1t_2=
\{t_1^it_2^j|1\leq i,j\leq p-1
\}$ are $\Z$-linearly independent in $\Z[G]$. From this the free rank claims and the
exactness of the sequence follow.

The discussions in 
(i), (ii), (iii) and (iv) clearly apply after
tensoring  with $\otimes_\Z\Z/p\Z$ and $\otimes_\Z\Z/p^2\Z$ giving (v). 
For the final
assertion in (vi) we note by Shapiro's Lemma the groups
$H^n(\G,\M_i\otimes_\Z\Z/p\Z)$ or $H^n(F,\M_i\otimes_\Z\Z/p\Z)$
are direct sums of $H^n(\HH,\Z/p\Z)$ where the $\HH$ are closed subgroups
of $\G$, or the Galois groups of  field extensions of $F$.
So the vanishing of the Bochstein maps follow  giving the result.
$\Box$

\bs

We now turn our attention to the module $\M_4:={\rm cok}(
d_2:\M_2\rightarrow\M_3)$.  It is somewhat technical, and
we summarize the results in the case of $p=3$ at the end of the
section to increase the readability of the applications when
$p=3$ in the subsequent sections.

\bs

\noi {\bf Definition 2.2.} We define $\alpha_1:=[(1,0,0,0)], \
\alpha_2:=[(0,1,0,0)]$, $ \beta_3:=[(0,0,T_3,0)]$,
 $\beta_4:=[(0,0,0,T_{p+1})] \in\M_4$. Here the brackets $[ \ ] $ 
 denote the class in $\M_4$.
 
 \bs
 
Clearly $\alpha_1,\alpha_2,\beta_3,\beta_4$ 
generate $\M_4$ as a $G$-module.
But by definition $0=[d_2(1,0)]=-t_1\alpha_1-t_2
\alpha_2+\beta_3+\beta_4$
so we can reduce our generating set to 
$\alpha_1$, $\alpha_2$ and $\beta_3$ (or $\alpha_1$, $\alpha_2$ and $\beta_4$).
We also note that $d_2(T_2,-T_G)=(-T_2t_1,0,0,0)$
and therefore $t_1T_2\alpha_1=0\in\M_4$ and  
by symmetry  
$T_1t_2\alpha_2=0\in\M_4$.
With  this
information we can describe $\M_4$ more thoroughly.

\bs

\noi {\bf Lemma 2.3.} {\em
The module
$\M_4$  is $\Z$-free, has $\Z$-rank $p^2+2p$ and is generated as a $\Z$-module by $\beta_4$ over
the submodule $\Z[G]\cdot\alpha_1+\Z[G]\cdot\alpha_2$ which is  $\Z$-free and has $\Z$-rank $p^2+2p-1$.  Also the $\Z[G]\cdot\alpha_i\cong\Z[G]/\Z[G]t_iT_{3-i}$ have rank $p^2-p+1$ and the $\Z[G]\cdot\beta_i\cong \Z[G]\cdot T_i$ are
$\Z$-free with $\Z$-rank $p$.}

\bs

\noi {\bf Proof.}
To see this
we first note that as $\M_1$, $\M_2$ and $\M_3$ are $\Z$-free
of ranks $1$, $p^2+1$ and $2p^2+2p$ respectively, it suffices to show that
$\M_4$ is generated by $p^2+2p$ elements, where we will see that their pullbacks to $\M_3$
together with the image of a basis of $d_2({\mathcal B})$ in $\M_3$ of
a basis ${\mathcal B}$ of $\M_2$ form a basis for $\M_3$. 

We begin by considering the subgroup $(0,\Z[G],0,0)\subset\M_3$
 and the corresponding subgroup
 $\Z[G]\cdot\alpha_2=\overline{(0,\Z[G],0,0)}\subset\M_4$. 
 By the the  previous observation that $T_1t_2\alpha_2=0\in\M_4$, we 
have a zero sequence 
$$
 0\rightarrow \Z[G]\cdot T_{1}t_2\rightarrow\Z[G]\stackrel{\cdot\alpha_2}{\longrightarrow}
\Z[G]\cdot\alpha_2\rightarrow0.$$
(This sequence will be seen to be exact at the end of the proof.)
One readily checks that  
$\Z[G]\cdot  T_{1}t_2$ has the $\Z$-basis, $(\tau_2-1)T_{1}$, 
$(\tau_2^2-1)T_{1},\ldots,$ $(\tau_2^{p-1}-1)T_{1}$ of $p-1$ elements, and therefore
$\Z[G]\cdot\alpha_2$ is generated as a $\Z$-module by the $p^2-p+1$ elements 
$\tau_1^i\tau_2^j\cdot \alpha_2$ with $(i,j)\in {\mathcal J}_2$ where
$${\mathcal J}_2:=\{(i,j)|0\leq i\leq p-2,0\leq j\leq p-1\}\cup\{(p-1,0)\}).$$
Pulling these elements back to $\M_3$ gives a set of $p^2-p+1$ elements
$${\mathcal A}_1:=\{(0,\tau_1^i\tau_2^j,0,0)\in\M_3|  (i,j)\in {\mathcal J}_2\}\subset\M_3,$$ which combined with the $p-1$ elements 
\begin{eqnarray*}
{\mathcal D}_1&:=&\{d_2(-(\tau_2^{k-1}+\cdots+\tau_2+1)T_1,kT_G)|
1\leq k\leq p-1\}\\
&=&\{(0,(\tau_2^k-1)T_1,0,0)|1\leq k\leq p-1\}
\end{eqnarray*}
  give a $\Z$-basis ${\mathcal A}_1\cup{\mathcal D}_1$ for $(0,\Z[G],0,0)\subset\M_3$. (To see this note that the $(\tau_2^k-1)T_1$ for $1\leq k\leq p-1$ have lead coefficients $\tau_1^{p-1}\tau_2^k$, precisely what is needed to fill out the 
  $\tau_1^i\tau_2^j$ where $(i,j)\in {\mathcal J}_2$ and give a basis.) We also note for use below that ${\mathcal D}_1\subset d_2(\M_2)$.

We next turn to the quotient $(\Z[G],\Z[G],0,0)/(0,\Z[G],0,0)$
and the corresponding quotient
$(\Z[G]\alpha_1+\Z[G]\alpha_2)/\Z[G]\alpha_2\subset
\M_4/ \Z[G]\alpha_2$. We start with $3p-2$ 
multiples of $\alpha_1$, $\{\tau_1^k\tau_2^\ell\alpha_1|(k,\ell)\in{\mathcal J}_1\}$ where
$${\mathcal J}_1:=\{(0,\ell) |0\leq\ell\leq p-1\}\cup\{(k,\ell)|1\leq k\leq2,0\leq \ell\leq p-2\}.$$
We note that
\begin{eqnarray*}
d_2(t_3t_{p+1},0)&=&(-t_1t_3t_{p+1},-t_2t_3t_{p+1},0,0)\\
&=&(\tau_1^3-\tau_1^2(1+\tau_2+\tau_2^{-1})+\tau_1(1+\tau_2+\tau_2^{-1})-1,-t_2t_3T_{p+1},0,0).
\end{eqnarray*}
and consider the $(p-3)p$ elements $\{d_2(\tau_1^i\tau_2^jt_3t_{p+1},0)|(i,j)\in{\mathcal J}_3\}\subset\M_3$  where  $${\mathcal J}_3:=\{(i,j)|0\leq i\leq p-4 \ \ {\rm and} \ \  0\leq j\leq p-1\}.$$  These have corresponding top terms 
$\tau_1^{i+3}\tau_2^{j\pm1}$, where we note that 
$i+3\in
\{3,4,\ldots,p-1\}$ and we view $j\pm1\in
\{0,1,2,\ldots,p-1\}$. We next note that 
\begin{eqnarray*}
d_2(-T_2,T_G)&=&((\tau_1-1)T_2,0,0,0) \ \ {\rm and}\\d_2(-(\tau_1+1)T_2,2T_G)&=&((\tau_1^2-1)T_2,0,0,0)
\end{eqnarray*}
 and these have corresponding lead terms $\tau_1\tau_2^{p-1}$ and 
$\tau_1^2\tau_2^{p-1}$.  

Altogether this information shows that if
\begin{eqnarray*}
{\mathcal A}_2&:=&\{(\tau_1^k\tau_2^\ell,0,0,0)| (k,\ell)\in{\mathcal J}_1\}\\
{\mathcal D}_2&:=&\{d_2(\tau_1^i\tau_2^jt_3t_{p+1},0)|(i,j)\in{\mathcal J}_3\}\\&=&\{(-\tau_1^i\tau_2^jt_1t_3t_{p+1},
-\tau_1^i\tau_2^jt_2t_3t_{p+1},0,0) 
| (i,j)\in{\mathcal J}_3\}\\
{\mathcal D}_3&:=&\{d_2(-T_2,T_G),d_2(-(\tau_1+1)T_2,2T_G)\}\\ &=&\{((\tau_1-1)T_2,0,0,0),((\tau^2_1-1)T_2,0,0,0)\},
\end{eqnarray*}
and if we denote by $\overline{\mathcal A}_2,\overline{\mathcal D}_2,
\overline{\mathcal D}_3$ the images of their elements, respectivly, in the
quotient module
$(\Z[G],\Z[G],0,0)/(0,\Z[G],0,0)$,
then the union $\overline{\mathcal A}_2\cup\overline
{\mathcal D}_2\cup\overline{\mathcal D}_3$ is
 a spanning set of $p^2$ elements for the quotient $(\Z[G],\Z[G],0,0)/(0,\Z[G],0,0)$,
for which the $p^2-3p+2$ elements in 
 ${\mathcal D}_2\cup{\mathcal D}_3$
lie in $d_2(\M_2)\subset\M_3$. 

Combing these two bases we have a spanning set for $(\Z[G],\Z[G],0,0)$ consisting of the of $2p^2$  elements 
${\mathcal A}_1\cup{\mathcal A}_2\cup{\mathcal D}_1\cup{\mathcal D}_2\cup{\mathcal D}_3$
 of which the  $p^2-3p+2+(p-1)=p^2-2p+1$ elements in
 ${\mathcal D}_1\cup{\mathcal D}_2\cup{\mathcal D}_3$
  lie in $d_2(\M_2)$.  This in turn means that we have $2p^2-(p^2-2p+1)=p^2+2p-1$ elements 
  $\{\tau_1^k\tau_2^\ell\cdot\alpha_1|(k,\ell)\in{\mathcal J}_1\}\cup
  \{\tau_1^i\tau_2^j\cdot\alpha_2|(i,j)\in{\mathcal J}_2\}$
that span $\Z[G]\alpha_1+\Z[G]\alpha_2\subset\M_4$.
(Note, adding the numbers of ${\mathcal J}_1$
and ${\mathcal J}_2$  elements  gives $(p^2-p+1)+(3p-2)=
p^2+2p-1$, the same number.)

	Next we note that  $\beta_3+\beta_4=t_1\alpha_1+t_2\alpha_2\in \Z[G]\alpha_1+\Z[G]\alpha_2$ and therefore $\M_4$ is spanned by  $Z[G]\alpha_1+\Z[G]\alpha_2$ and $\beta_4$ over $\Z$. Further, 
$$
d_2((\tau_1\tau_2-1),0)=(-t_1(\tau_1\tau_2-1),-t_2(\tau_1\tau_2-1),0,(\tau_1\tau_2-1)T_{p+1})$$
and therefore
$(\tau_1\tau_2-1)\beta_4\in\Z[G]\alpha_1+\Z[G]\alpha_2$, that is,
$\tau_1\tau_2$ acts trivially on the class $\overline{\beta}_4$ in $\M_4/(\Z[G]\alpha_1+\Z[G]\alpha_2)$.  But also, $\tau_1\tau_2^{-1}\cdot\beta_4=\beta_4$ so the quotient module 
$\overline{\Z[G]\beta_4}\subset\M_4/(\Z[G]\alpha_1+\Z[G]\alpha_2)$
is isomorphic to $\Z$.
This means that $\M_4$ spanned as a $\Z$-module by $\Z[G]\alpha_1+\Z[G]\alpha_2$ and $\beta_4$ and we have created a $\Z$-spanning set for $\M_4$ of $p^2+2p$ elements.  In fact these form a $\Z$-free basis for $\M_4$, for the exactness of the main sequence shows that torsion free part of $\M_4$ must have $\Z$-rank $p^2+2p$.
	
	We continue to fill out the basis for $\M_3/(\Z[G],\Z[G],0,0)$. We have the following $2p$ generators 
\begin{eqnarray*}
{\mathcal A}_3&:=&\{(0,0,0,T_{p+1})\}\\
{\mathcal D}_4&:=&\{d_2((\tau_1\tau_2)^\ell-1,0)|1\leq\ell\leq p-1\}\\&=&
\{(-((\tau_1\tau_2)^\ell-1)t_1,-((\tau_1\tau_2)^\ell-1)t_2,0,-((\tau_1\tau_2)^\ell-1)T_{p+1})|1\leq\ell\leq p-1\}\\
{\mathcal D}_5&:=&\{d_2((\tau_1\tau_2^{-1})^\ell-1,0)|1\leq\ell\leq p-1
\}\\&=&\{(-((\tau_1\tau_2^{-1})^\ell-1)t_1,-((\tau_1\tau_2^{-1})^\ell-1)t_2,-((\tau_1\tau_2^{-1})^\ell-1)T_3,0)|1\leq\ell\leq p-1\}\\
{\mathcal D}_6&:=&\{d_2(1,0)\}=\{(-t_1,-t_2,T_3,T_{p+1})\}
\end{eqnarray*}
Clearly these $2p$ generators generate $(0,0,\Z[G]T_3,\Z[G]T_{p+1})$ where $2p-1$ of them in ${\mathcal D}_4\cup{\mathcal D}_5\cup{\mathcal D}_6$, lie in $d_2(\M_2)$. 

Altogether we have shown that we have $2p^2+2p$ $\Z$-generators for $\M_3$, and hence a $\Z$-basis for $\M_3$ with  $(p^2-2p+1)+(2p-1)=p^2$  in $\cup_{j=1}^6{\mathcal D}_j\subset d_2(\M_2)$ and the remaining $p^2+2p$ in $\cup_{i=1}^3{\mathcal A}_i$ map to a basis for $\M_4$.  

It remains to check that the $\Z$-rank of $\Z[G]\cdot\beta_4$ is $p$. If $\mu\beta_4=0\in\M_4$, then there would exist $(\eta,nT_G) \in \M_2$ with
$d_2(\eta,nT_G)=(0,0,0,\mu T_{p+1})$. But then we would
have $\eta t_1=\eta t_2=0$ which by the exactness arguments in the 
proof of Theorem 1.3 gives $\eta=m T_G$ for some $m\in \Z$. We would have $d_2(\eta,nT_G)=(0,0,(pm+n)T_G,(pm+n)T_G)$ so that $\mu T_{p+1}=(pm+n)T_G=0$. Therefore $\Z[G]\cdot\beta_4\cong \Z[G]\cdot T_{p+1}$.
As $\Z[G]\cdot T_{p+1}$ is free of rank $p$ this proves the result. $\Box$

\bs

In the following we denote by ${\mathcal D}:=\cup_{j=1}^6{\mathcal D}_j$,
$ {\mathcal A}:=\cup_{i=1}^3{\mathcal A}_j$ and by 
$\overline{\mathcal A}$ the classes of elements in ${\mathcal A}$ in
$\M_4$.  
	
\bs

\noi {\bf Corollary 2.4.} {\em (i) We have that
$ {\mathcal A}\cup {\mathcal D}$ is a $\Z$-basis for $\M_3$ where
$ \overline{\mathcal A}$ is a $\Z$-basis for $\M_4$
and ${\mathcal D}$ is a $\Z$-basis for $d_2(\M_2)$.\\
(ii) We have a short exact sequence
of $G$-modules that are $\Z$-free as $\Z$-modules
$$0\longrightarrow (\Z[G]\cdot\alpha_1+\Z[G]\alpha_2)\longrightarrow
{\mathcal M}_4\stackrel{\pi_4}{\longrightarrow}
\Z[G]\cdot T_G\longrightarrow0$$
where $\pi_4$ is given by
$\pi_4([(\eta_1 t_1,\eta_2 t_2,\eta_3 T_3,\eta_4 T_{p+1})]=
(\eta_3-\eta_4) T_G$.}

\bs

\noi {\bf Proof.}  For (i) the basis is extracted from the proof of Lemma 2.3.

For (ii) the exactness of the sequence  follows from (i) once it is shown
that the  $\pi_4$ is well-defined. To check the well-definition of
$\pi_4$ we note that if $(\eta_1 t_1,\eta_2 t_2,\eta_3 T_3,\eta_4 T_{p+1})\in {\mathcal D}:=\cup_{j=1}^6{\mathcal D}_j$ then
$(\eta_4-\eta_3) T_G=0$.
$\Box$

\bs

\noi {\bf Corollary 2.5.} {\em 
If $(\theta_1,\theta_2,0,\theta_4T_{p+1})\in d_2(\M_2\otimes_\Z\Z/p\Z)$ then 
it is a $\Z/p\Z$ linear combination of ${\mathcal D}_1,{\mathcal D}_2,{\mathcal D}_3$ and ${\mathcal D}_4$ basis elements  (mod $p$), and therefore 
we can
express 
$\theta_4T_{p+1}=\sum_{k=1}^{p-1}a_k((\tau_1\tau_2)^\ell-1)T_{p+1}$ for $a_k\in\Z/p\Z$.}

\bs

\noi {\bf Proof.}   First we note that (mod $p$) 
only the basis elements in ${\mathcal D}_4,{\mathcal D}_5,{\mathcal D}_6$ can
contribute to an element in $d_2(\M_2/p)$ where at least one of latter
two coordinates are nonzero.  However, ${\mathcal D}_4$ elements have third coordinate zero, the third coordinates of 
${\mathcal D}_5$ are $-((\tau_1\tau_2^{-1})^\ell-1)T_3$ where $1\leq\ell\leq p-1$
and the third coordinates of the ${\mathcal D}_6$ element
 is $T_3$.  These third coordinates are 
$\Z$-linearly independent (and so of course are
$\Z/p\Z$-linearly independent)
and therefore no nontrivial linear combination of these can have $\theta_3T_3$ term zero. It follows that any element of
$d_2(\M_2\otimes\Z/p\Z)$   with zero third coordinate must be a 
$\Z/p\Z$ sum of elements in ${\mathcal D}_1,{\mathcal D}_2,{\mathcal D}_3$ and ${\mathcal D}_4$.  The required form for $\theta_4$
follows because it must be a sum of 
fourth entries from ${\mathcal D}_4$ elements. This proves the corollary.  $\Box$

\bs

The next Corollary  describes the
$\Z$-module generated by ${\mathcal D}_1,{\mathcal D}_2,{\mathcal D}_3$ and ${\mathcal D}_4$
with a set of $\Z[G]$-module generators. This result, of course, also applies when
viewing the corresponding $\Z/p\Z[G]$ modules.

 \bs
 
 \noi {\bf Corollary 2.6.} {\em We set
 $${\mathcal P}=\Z[G]\cdot(t_1T_2,0,0,0)+
\Z[G]\cdot(0,T_1t_2,0,0)+\Z[G]\cdot
(-t_1t_3,-t_2t_3,0,t_3T_{p+1})\subset \M_3.$$
Then ${\mathcal P}$ is the set of $\Z$-linear combinations of
 ${\mathcal D}_1,{\mathcal D}_2,{\mathcal D}_3$ and ${\mathcal D}_4$ basis elements
 and so we have the exact sequence
$$0\rightarrow{\mathcal P}\rightarrow
 \Z[G]\oplus\Z[G]\oplus0\oplus\Z[G]\beta_4\rightarrow
\M_4\rightarrow0.
$$}

\noi {\bf Proof.} We have $(0,T_1t_2,0,0)\in{\mathcal D}_1$,
$(t_1T_2,0,0,0)\in {\mathcal D}_3$ and 
$(-t_1t_3,-t_2t_3,0,t_3T_{p+1})\in {\mathcal D}_4$. 
From this it follows that ${\mathcal P}$ is a subset of the $\Z[G]$-module which is the 
$\Z$-linear combinations of
 ${\mathcal D}_1,{\mathcal D}_2,{\mathcal D}_3$ and ${\mathcal D}_4$ basis elements.
 For the reverse inclusion, one readily checks that  ${\mathcal D}_1\subseteq \Z[G]\cdot(0,T_1t_2,0,0)$,
 ${\mathcal D}_2\subseteq \Z[G]\cdot
(-t_1t_3,-t_2t_3,0,t_3T_{p+1})$, ${\mathcal D}_3\subseteq \Z[G]\cdot(t_1T_2,0,0,0)$
and ${\mathcal D}_4\subseteq \Z[G]\cdot
(-t_1t_3,-t_2t_3,0,t_3T_{p+1})$.  So these $\Z[G]$-modules agree. 
  By 
 Corollary 2.5
 ${\mathcal P}=d_2(\M_2)\cap(\Z[G],\Z[G],0,\Z[G]T_4)\subset\M_3$
 and the result follows. $\Box$

\bs 

This final Corollary focuses on $(\Z[G],\Z[G],0,0)$ and the image of a variation of the previous bases inside this submodule.

\bs

\noi {\bf Corollary 2.7.} {\em A basis for $d_2(\M_2)\cap (\Z[G],\Z[G],0,0)$ in $\M_3$ consists of the  $p^2-2p+1$ elements in
${\mathcal D}_0:={\mathcal D}_1\cup{\mathcal D}_2\cup{\mathcal D}_3$. Alternatively
we can take as a basis the following elements
\begin{eqnarray*}
d_2(-T_2,T_G)&=&(t_1T_2,0,0,0)\\
d_2(-T_1,T_G)&=&(0,T_1t_2,0,0)\\
d_2(-t_1^iT_2,0)&=&(t_1^{i+1}T_2,0,0,0) \ \ {\rm for} \ \ 1\leq i\leq p-2\\
d_2(-T_1t_2^j,0)&=&(0,T_1t_2^{j+1},0,0) \ \ {\rm for} \ \ 1\leq j\leq p-2\\
d_2(-\tau_1^i\tau_2^jt_3t_{p+1},0)&=&(\tau_1^i\tau_2^jt_1t_3t_{p+1},
\tau_1^i\tau_2^jt_2t_3t_{p+1},0,0)
\  {\rm for} \ (i,j)\in{\mathcal J}'_3
 \ \ {\rm  where} \\&& \ \  {\mathcal J}'_3:=\{(i,j)|0\leq i\leq p-4
\ \ {\rm and} \ \ 0\leq j\leq p-2\}.
\end{eqnarray*}
Here, except for the first two basis elements where $h_1$ has value $T_G$, the value of $h_1$ is $0$. }

\bs

\noi {\bf Proof.} The first statement follows by first checking that 
${\mathcal D}_0\subset d_2(\M_2)\cap (\Z[G],\Z[G],0,0)$ and then noting
by inspection that if a $\Z$-linear combination of elements in 
${\mathcal D}_4\cup{\mathcal D}_5\cup{\mathcal D}_6$ lies in $(\Z[G],\Z[G],0,0)$ 
then it is necessarily $0$.

For the second statement, since the list has $p^2-2p+1$ elements, it suffices to show that its $\Z$-span contains ${\mathcal D}_0$.   
For this
first we note that the second and fourth elements in the list are
extracted from and  a minor variation of ${\mathcal D}_1$ which is 
therefore included in their span.  Next we note that
${\mathcal J}_3={\mathcal J}_3'\,\dot\cup\,\{(i,p-1)|0\leq i\leq p-4\}$ which means
that the linear combinations of elements in ${\mathcal D}_3$ with 
$(i,j)\not\in {\mathcal J}_3'$ each have as first entry some nonzero $\tau_1^i\tau_2^{p-1}$ term when viewed as a polynomial in $\tau_1,\tau_2$.  For such a 
combination we can subtract off a combination of the $p-1$ first and third entries of the list so that the highest degree $\tau_2$-term of 
the first entry is at most $p-2$.  These are necessarily combinations
of the fifth set of elements.
The final statement is a direct calculation using the definition of $h_1$.
$\Box$

	\bs

We close this section with proofs of versions of Lemma 2.3 and Corollary 2.4 in the
case where $p=3$.  These results will be used in the
next section where the cohomological kernels of a elementary abelian
degree $9$ extensions are considered, and in fact is all that is necessary
for the remainder of this paper.  Hopefully the other results in this section
will be useful as generalizations to the case of $p>3$
are investigated.

\bs

\noi {\bf Lemma 2.3.} {\em ($p=3$ case.)
The module
$\M_4$  is $\Z$-free and has $\Z$-rank $15$. Moreover,
the submodules $\Z[G]\cdot\alpha_i$ are  $\Z$-free and have $\Z$-rank $7$ while $\Z[G]\cdot\beta_i\cong \Z[G]\cdot T_i$ are
$\Z$-free and have $\Z$-rank $3$.}

\bs

\noi {\bf Proof.}
To see this
we first note that as $\M_1$, $\M_2$ and $\M_3$ are $\Z$-free
of ranks $1$, $10$ and $24$ respectively, it suffices to show that
$\M_4$ is generated by $15$ elements, where we will see that their pullbacks to $\M_3$
together with the image of a basis of $d_2({\mathcal B})$ in $\M_3$ of
a basis ${\mathcal B}$ of $\M_2$ form a basis for $\M_3$. 

As  previously observed, we 
have a zero sequence
$$(*_{\alpha_1}) \ \ \ \ \ 0\rightarrow \Z[G]\cdot t_1 T_2\rightarrow\Z[G]\stackrel{\cdot\alpha_1}{\longrightarrow}
\Z[G]\cdot\alpha_1\rightarrow0$$
and moreover using Lemma 2.1 (ii) one can check that  
$\Z[G]\cdot t_1 T_2$ has the $\Z$-free basis, $(\tau_1-1)T_2$, 
$(\tau_1^2-1)T_2$ and therefore 
$\Z[G]\cdot\alpha_1$ is generated as a $\Z$-module by the seven elements $\alpha_1,\tau_2\alpha_1,\tau_2^2\alpha_1,\tau_1\alpha_1,\tau_1\tau_2\alpha_1,
 \tau_1^2\alpha_1,\tau_1^2\tau_2\alpha_1$.
Combining this with the same result for $\alpha_2$ shows that
$\Z[G]\cdot\alpha_1+\Z[G]\cdot\alpha_2$ is generated 
as a $\Z$-module by the 
fourteen elements $\alpha_1,\tau_2\alpha_1,\tau_2^2\alpha_1,\tau_1\alpha_1,\tau_1\tau_2\alpha_1,$
 $ \tau_1^2\alpha_1,\tau_1^2\tau_2\alpha_1,
 \alpha_2,\tau_1\alpha_2,\tau_1^2\alpha_2,\tau_2\alpha_2,\tau_1\tau_2\alpha_2,
 \tau_2^2\alpha_2,\tau_1\tau_2^2\alpha_2$.

	Next we note that  $\beta_3$ generates
	$\M_4$ over $\Z[G]\cdot\alpha_1+\Z[G]\cdot\alpha_2$
as a $\Z$-module. For this, as $\tau_1\tau_2T_3=T_3$ we have
$\tau_1\tau_2\beta_3=\beta_3$. Also as $(1-\tau_1\tau_2^2)T_4=0$ we have
$$d_2(1-\tau_1\tau_2^2,0)=(-(1-\tau_1\tau_2^2)t_1,
-(1-\tau_1\tau_2^2)t_2,(1-\tau_1\tau_2^2)T_3,0)$$ and we find that
$(1-\tau_1\tau_2^2)\beta_3=(1-\tau_1\tau_2^2)(\alpha_1+
\alpha_2)\in \Z[G]\cdot\alpha_1+\Z[G]\cdot\alpha_2$.
This shows that $\tau_2\beta_3
=\tau_2(\tau_1\tau_2)\beta_3=\tau_1\tau_2^2\beta_3$
and $\tau_1\beta_3=\tau_1(\tau_1\tau_2)^2\beta_3=\tau_2^2\beta_3$ and  both lie in the
$\Z$-module generated by $\alpha_1,\alpha_2$ and
$\beta_3$. Therefore the fifteen elements $\alpha_1,\tau_2\alpha_1,\tau_2^2\alpha_1,\tau_1\alpha_1,\tau_1\tau_2\alpha_1,$
 $ \tau_1^2\alpha_1,\tau_1^2\tau_2\alpha_1,
 \alpha_2,\tau_1\alpha_2,\tau_1^2\alpha_2,\tau_2\alpha_2,\tau_1\tau_2\alpha_2,$
 $ \tau_2^2\alpha_2,\tau_1\tau_2^2\alpha_2,\beta_3$ 
generate $\M_4$ as required.

These fifteen generators for $\M_4$ pull back to the fifteen elements
\begin{eqnarray*}
\{1,\tau_2,\tau_2^2,\tau_1,\tau_1\tau_2,\tau_1^2,\tau_1^2\tau_2\}\cdot(1,0,0,0),&&\\
\{1,\tau_1,\tau_1^2,\tau_2,\tau_1\tau_2,\tau_2^2,\tau_1\tau_2^2\}\cdot(0,1,0,0),&&\\(0,0,T_3,0)&\in &\M_3
\end{eqnarray*} 
which we will denote by ${\mathcal A}$.
Next we consider ${\mathcal A}_0:=\{1,\tau_1\tau_2,\tau_1^2\tau_2^2,\tau_1^2\tau_2,
\tau_1\tau_1^2,$ $(1+\tau_1+\tau_1^2),$  $(\tau_2+\tau_2\tau_1+\tau_2\tau_1^2),$
$ (1+\tau_2+\tau_2^2),(\tau_1+\tau_1\tau_2+\tau_1\tau_2^2)\}$. Inspecting
the top degree terms of ${\mathcal A}_0$ we see that it is a
$\Z$-basis for $\Z[G]$.
This means we have $9$ 
 $\Z$-independent elements $${\mathcal D}_0:=
 \{(\{1,\tau_1\}T_2,-T_G),(\{1,\tau_2\}T_1,-T_G),
 (\{1,\tau_3,\tau_3^2\},0),(\{\tau_4,\tau_4^2\},0)\}\subset\M_2$$
 and using 
 $(T_1,-T_G)+(T_2,-T_G)+(\tau_1\tau_2+\tau_1^2\tau_2^2+\tau_1^2\tau_2+
\tau_1\tau_1^2-1,0)= (T_G,-2T_G)$
we see that $(T_G,-2T_G)\in$ span $({\mathcal D}_0)\subset\M_2$.
But then  if we subtract $d_1(T_G)=(T_G,-3T_G)\in {\rm im}(d_1)
\subset\M_2$ we obtain
$(0,T_G)\in$ span (${\mathcal D}_0\cup\{(T_G,-3T_G)\}$), and
 from this we see that ${\mathcal D}_0\cup\{(T_G,-3T_G)\}$ is a $\Z$-basis
 for $\M_2$.
 We now consider the nine images under $d_2$,
\begin{eqnarray*}
d_2((\{1,\tau_1\}T_2,-T_G))&=&(-\{1,\tau_1\}T_2t_1,0,0)\\
d_2((\{1,\tau_2\}T_1,-T_G))&=&(0,-\{1,\tau_2\}T_1t_2,0,0)\\
d_2((\{1,\tau_3,\tau_3^2\},0))&=&(-\{1,\tau_3,\tau_3^2\}t_1,
-\{1,\tau_3,\tau_3^2\}t_2,T_3,\{1,\tau_3,\tau_3^2\}T_4)\\
d_2((\{\tau_4,\tau_4^2\},0))&=&(-\{\tau_4,\tau_4^2\}t_1,
-\{\tau_4,\tau_4^2\}t_2,\{\tau_4,\tau_4^2\}T_3,T_4),
\end{eqnarray*}
which we denote by ${\mathcal D}$.
Altogether combining the $15$ elements in ${\mathcal A}$ (the pull backs of the generators
for $\M_4$ in $\M_3$) with these latter $9$ elements in 
${\mathcal D}\subset d_2(\M_2)\subset\M_3$ gives ${\mathcal A}\cup{\mathcal D}$ which
can be seen to be
a $\Z$-basis of $\M_3$ as follows. First note that
bases for $(\Z[G],0,0,0)$ and $(0,\Z[G],0,0)$, ($9$ elements
in total)
are included, and then linear combinations of these can be subtracted off
the remaining six elements to yield a basis for $(0,0,\Z[G]T_3,\Z[G]T_4)$,
first by subtracting $(0,0,T_3,0)$ from  the three elements
in $d_2((\{1,\tau_3,\tau_3^2\},0))$ to obtain a basis for $(0,0,0,\Z[G]T_4)$; and
then using three remaining for a total of six elements as a basis for
$(0,0,\Z[G]T_3,\Z[G]T_4)$.

As a consequence of these observations we find 
$\alpha_1,\tau_2\alpha_1,\tau_2^2\alpha_1,\tau_1\alpha_1,$
$\tau_1\tau_2\alpha_1,
 \tau_1^2\alpha_1,\tau_1^2\tau_2\alpha_1$
are linearly independent in $\M_4$ and therefore
 the sequence $(*_{\alpha_1})$ is exact and 
 $\Z[G]\cdot\alpha_1\subset\M_4$ has rank $7$.  It also follows that 
$\Z[G]\cdot\alpha_1\cap \Z[G]\cdot\alpha_2=\{0\}$.  Finally  by Lemma 2.1 (i),
 $\Z[G]\cdot\beta_3$ has at most $\Z$-rank $3$. But using $T_1\alpha_i=T_2\alpha_i=0$ our calculations show 
  $\tau_2\beta_3=\tau_1\tau_2^2\beta_3=\beta_3-
(1-\tau_1\tau_2^2)(\alpha_1+\alpha_2)=\beta_3-
(1+\tau_1+\tau_1\tau_2)(\alpha_1+\alpha_2)$ and 
 $\tau_2^2\beta_3=\beta_3-(\tau_1+\tau_2+\tau_1\tau_2)(\alpha_1+\alpha_2)$
 showing that $\Z[G]\cdot\beta_3$ has $\Z$-rank at least $3$. In particular,
 $\Z[G]\cdot\beta_3\cong \Z[G]\cdot T_3$. $\Box$

\bs

\noi {\bf Corollary 2.4.} {\em ($p=3$ case.) We have a short exact sequence
of $G$-modules that are $\Z$-free as $\Z$-modules
$$0\longrightarrow (\Z[G]\cdot\alpha_1\oplus\Z[G]\alpha_2)\longrightarrow
{\mathcal M}_4\stackrel{\pi_4}{\longrightarrow}
\Z[G]\cdot T_G\longrightarrow0$$
where $\pi_4$ is given by
$\pi_4([\eta_1 t_1,\eta_2 t_2,\eta_3 T_3,\eta_4 T_4)]=
(\eta_3-\eta_4) T_G$.
The analogous result is true (mod $3$), e.g tensoring
$\otimes_\Z\ZZ$. }

\bs

\noi {\bf Proof.}  The directness of the sum $\Z[G]\cdot\alpha_1\oplus\Z[G]\alpha_2$ 
follows from the $14$ elements in the first fifteen generators for $\M_3$ in the
proof that $\M_4$ has rank $15$. The cokernel of the inclusion
$\Z[G]\cdot\alpha_1\oplus\Z[G]\alpha_2\hookrightarrow \M_4$ then has rank
$1$ it remains to show that the $\pi_4$ is well-defined. This is a 
consequence of $\beta_3+\beta_4=t_1\alpha_1+t_2\alpha_2\in \M_4$.
$\Box$

\bs

\bs

\noi {\bf \S 3.  \ Analysis of 
$H^n(\G,\M_4\otimes_\Z\ZZ)$ and the Cohomological Kernel.}

\bs

We now assume $p=3$ and study the cohomology 
group $H^n(\G,\M_4\otimes_\Z\ZZ)$. 
The main result in this section is Theorem 3.14 where the cohomological kernel of an elementary abelian $3$-$3$ extension is characterized.  As discussed in Section 1, throughout this
paper we assume that $\G$ is what we call a Voevodsky group, which gives  that for $i=1,2,3$
 the Bochsterin maps
${\mathcal{B}}_i:H^n(\G,\M_i\otimes_\Z\Z/3\Z)\rightarrow
H^{n+1}(\G,\M_i\otimes_\Z\Z/3\Z)$ are zero
and that every cyclic quotient of index three extends to a (non unique) cyclic
quotient of index nine.  

To proceed we set up some notation involving cyclic quotients of $\G$
of degree $9$. 
Suppose  $0\neq\theta\in H^1(\G,\ZZ)$
and $\tilde{\theta}\in H^1(\G,\Z/9\Z)$ is a lift (so that
$\theta=\pi\circ\tilde{\theta}$ where $\pi:\Z/9\Z\rightarrow\ZZ$ is the 
natural projection.) 
We set  ${\mathcal H}_{\theta=0}:={\rm ker}(\theta)\subset\G$ and
  ${{\mathcal H}}_{\tilde{\theta}=0}:={\rm ker}(\tilde{\theta})\subset
{\mathcal H}_{\theta}=0\subset\G$, both  normal in $\G$. 
We then have: \\ 
(i)
$\G/{{\mathcal H}}_{\tilde{\theta}=0}\cong\Z/9\Z$
with $H^1(\G/{{\mathcal H}}_{\tilde{\theta}=0},\Z/9\Z)=\langle\tilde{\theta}\rangle$;
\\ (ii)
$\G/{{\mathcal H}}_{\theta=0}\cong\Z/3\Z$
with $H^1(\G/{{\mathcal H}}_{{\theta}=0},\Z/3\Z)=\langle{\theta}\rangle$; 
\\ (iii)
 ${{\mathcal H}}_{\theta=0}/{{\mathcal H}}_{\tilde{\theta}=0}
\cong\Z/3\Z\cong\langle\overline{\sigma}^3\rangle$
with $H^1({{\mathcal H}}_{\theta=0}/{{\mathcal H}}_{\tilde{\theta}=0}
,3\Z/9\Z)=\langle\tilde{\theta}|_{{{\mathcal H}}_{\theta=0}}\rangle$. \\
We remark at this point
that the map $\ZZ\rightarrow \Z/9\Z$ given by $1\mapsto 3$ corresponds
in the field case
to the embeddings in Galois cohomology where we identify or
primitive cube root of unity with $\zeta^3$ where $\zeta$ is a fixed
primitive $9$th root of unity, and that is why the notation with
$H^1$ is spelled out as it is.  We also note that
in this setting, if we identify $\G/{{\mathcal H}}_{\theta=0}\cong\langle\tau\rangle$
where $\tau^3=1$ the $\G$-module $\Z[\langle\tau\rangle]$ induced by
this quotient is in fact the co-induced module for the subgroup
${{\mathcal H}}_{\theta=0}$.

 Since we will be computing with cocyles and not cohomology classes we
 need the following (somewhat unconventional) notation
 to increase the
 readability of the computations.  
 
 \bs

\noi {\bf Notation 3.1} (i) For functions $x\in C^1(\G,\ZZ)$ 
and $y\in C^{n-1}(\G,\ZZ)$ we denote
by $x\cdot y'\in  C^n(\G,\ZZ)$ the function 
$(x\cdot y')(\gamma,\gamma'):=x(\gamma)\cdot y(\gamma')$
for $\gamma\in \G$ and $\gamma'\in \G^{n-1}$.
For notational convenience, the ``dot'', $\cdot$, will usually be deleted.\\
(ii) For $c\in  Z^n(\G,\ZZ)$ we denote by
$[c]\in  H^n(\G,\ZZ)$ the cohomology class of $c$.\\
(iii) For $\psi\in Z^1(\G,\ZZ)$ and $\chi\in Z^{n-1}(\G,\ZZ)$ we use the conventional
notation for the cup product $[\psi]\smile[\chi]
\in H^n(\G,\ZZ)$. In terms of the previous notation this
means that $[\psi]\smile[\chi]=[\psi\cdot\chi']$. To increase the readability of the computations
we will regularly drop the brackets and when no confusion can arise
use the same notation for a cohomology class as a cocycle it represents. For example the cohomology class of $\chi\in Z^{n-1}(\G,\ZZ)$ may also be denoted by 
$\chi\in H^{n-1}(\G,\ZZ)$ if no confusion can arise. \\
(iv) We denote by $\lambda^2_i:\ZZ\rightarrow\ZZ$
the function $\lambda^2_i:=\binom{i}{2}$ (mod $3$), that is
$\lambda^2_i=0$ if $i=0,1$ and $\lambda^2_2=1$.

\bs

\noi {\bf Remark 3.2.} The purpose of the above notation, especially the prime notation
in (i) where the prime is used to separate out variables, will become apparent as the
computations proceed and we need to look at the actual functions
not just cohomology classes. We will often compose $\lambda^2_i$ with other
functions in $C^1(\G,\ZZ)$ to define more complicated
functions; for example we will use the notation  $\lambda^2_\theta$ to denote the
composite of the character $\theta$ with $\lambda^2_i$.
In subsequent work where $p>3$,
for $1\leq k<p$ we set $\lambda^k_i:\Z/p\Z\rightarrow\Z/p\Z$ to be
the function $\lambda^k_i:=\binom{i}{k}$ \, (mod $p$), where the notation is chosen
so that the entries of the binomial coefficients appear upside down to remind us that
we are working (mod $p$).

\bs

We will
need to lift functions in $C^n(\G,\M_4\otimes_\Z\ZZ)$ to functions
in $C^n(\G,\M_4\otimes_\Z\Z/9\Z)$  in two different ways and so two notations are needed. 
These are described next.

\bs

\noi {\bf Definitions 3.3.} (i)
For $f\in C^1(\G,\ZZ)$ we
 set $\widehat{{f}}\in C^1(F,\Z/9\Z)$ to
be the extension given by the requirement that $\widehat{{f}}(x)\in
\{0,1,2\}$ (mod $9$) and ${\widehat{{f}}}(x)\equiv
{{f}}(x)$ (mod $3$).\\
(ii) For  ${{\theta}}\in  Z^1(\G,\ZZ)$  we will denote by
    $\tilde{{{\theta}}}\in Z^1(\G,\Z/9\Z)$ a lift, which are not unique but for 
 which in the applications any fixed choice will suffice. \\
(iii) For $f\in C^n(\G,\ZZ)$ it will be convenient to denote by $3f\in
C^n(\G,\Z/9\Z)$ the function whose values are $3f(x)\in \Z/9\Z$, which are determined
by the embeding $\ZZ\rightarrow\Z/9\Z$ given by $1$
(mod $3$) $\mapsto 3$ (mod $9)$.\\
 (iv)
With these definitions
we can extend $\lambda^2_\theta$ to $\lambda^2_{\widehat{\theta}}:=
\binom{\widehat{\theta}}{2}$ (mod $9$) and $\lambda^2_{\tilde{\theta}}:=
\binom{\tilde{\theta}}{2}$ (mod $9$).\\
(v)
Suppose that $\tilde{\theta}\in H^1(\G,\Z/9\Z)$ is a
fixed  lift
of $\theta\in H^1(\G,\ZZ)$, ${{\mathcal H}}_{\tilde{\theta}=0}={\rm ker}(\tilde{\theta})
\subset \G$ and that 
$\G/ {{\mathcal H}}_{\tilde{\theta}=0}$ is generated by $\sigma {{\mathcal H}}_{\tilde{\theta}=0}$
for a fixed $\sigma\in \G$.  Define
$s_\theta\in  C^1(\G,\ZZ)$  by
$s_\theta(\sigma^kh)=[k/3]\in \ZZ$ where $h\in 
{{\mathcal H}}_{\tilde{\theta}=0}$, $k\in \{0,1,\ldots,8\}$  and the bracket denotes the greatest
integer function but taking the value $($mod $3)$. 

\bs

For
${{u}}\in  C^1(\G,\M\otimes_Z\ZZ)$ we will consider
 lifts to $\tilde{{{u}}}\in C^1(\G,\M\otimes_\Z\Z/9\Z)$, often  built out of lifts as defined above.
 In 
 cases where the Bochstein is known to be trivial (e.g.
  $H^n(\G,\M\otimes_\Z\Z/9\Z)\rightarrow
  H^n(\G,\M\otimes_\Z\ZZ)$  is a surjection), 
 for  $u\in  Z^n(\G,\M\otimes_\Z\ZZ)$ we will 
 need to consider lifts to
$\tilde{u}\in Z^{n}(\G,\M\otimes_\Z\Z/9\Z)$. 
  
  \bs

\noi {\bf Conventions 3.4.}
Functions $f$ without any hats or tildes over them
will be understood as $f\in C^n(\G,\ZZ)$,
while the functions $\widehat{f}$, $\tilde{f}$ or $3f$ will be
understood to lie in $C^n(\G,\Z/9\Z)$.
We do this so there
is no ambiguity about how we use $+$, that is, whether one is adding in $\ZZ$ or $\Z/9\Z$.
 We will use this notation regularly
to remove hats and tildes one we reach multiples of three; for example
for $f,g\in C^1(\G,\ZZ)$, instead of $\widehat{f}+3\widehat{g}$ we will write $\widehat{f}+3g$.  In particular, because of this slight, but convenient, abuse of notation, $3f\neq0$ in general, and the intended meaning will be clear
from context.

\bs

As usual, $\delta$ denotes the coboundary
 $\delta:C^m(\G,\Z/p^\ell\Z)\rightarrow C^{m+1}(\G,\Z/p^\ell\Z)$.
We will need to compute $\delta$ values repeatedly and some
basic values  are given next. 

\bs

\noi {\bf Lemma 3.5.} If ${{\theta}}\in  Z^1(\G,\ZZ)$ then\\
(i) $\delta(s_{{\theta}})=-
({{\theta}}\lambda^2_{{{\theta}}'}+\lambda^2_{{\theta}}{{\theta}}')
\in Z^2(\G,\ZZ).$\\
(ii) 
$  \delta(\widehat{{\theta}})
=3({{\theta}}\lambda^2_{{{\theta}}'}+\lambda^2_{{\theta}}{{\theta}}')
\in Z^2(\G,\Z/9\Z).$\\
(iii) $\delta(\lambda^2_{{\theta}})=-{\theta}{\theta}'
\in Z^2(\G,\ZZ)$
and $\delta(\lambda^2_{\tilde{\theta}})=-\tilde{\theta}\tilde{\theta}'
\in Z^2(\G,\Z/9\Z)$.\\
(iv) $\delta(\theta s_\theta)
=-(\lambda^2_\theta\lambda^2_{\theta'}+s_\theta\theta'+\theta s_{\theta'})
\in Z^2(\G,\ZZ)$.\\
(v) $\delta(\lambda^2_{\widehat{\theta}})=-\widehat{\theta}\,\widehat{\theta}'+3\left(
\lambda^2_{{\theta}}\lambda^2_{{\theta}'}
+\theta\lambda^2_{\theta'}+\lambda^2_\theta\theta'\right)
\in Z^2(\G,\Z/9\Z).
$\\
(vi) $\delta(\widehat{s_\theta})=-(\widehat{\theta}\lambda^2_{\widehat{\theta}'}+\lambda^2_
{\widehat{\theta}}\widehat{\theta}')+3\lambda^2_\theta\lambda^2_{\theta'}
\in Z^2(\G,\Z/9\Z).$\\
(vii) For $\gamma\in\langle\tau\rangle$ and $t:=\tau-1$ we have  $(\gamma-1)\cdot 1=\widehat{\theta}(\gamma)t+
\lambda^2_{\widehat{\theta}(\gamma)}t^2\in\Z/9\Z[\langle\tau\rangle]$.

\bs

\noi {\bf Proof.} 
For (i) suppose $\gamma=\sigma^kh$, 
 $\gamma'=\sigma^{k'}h'$ so $\gamma\gamma'=\sigma^{k+k'}h''$.
 Set $\ell,\ell'\in\{0,1,2\}$ so that $k\equiv \ell$
 (mod $3$) and $k'\equiv \ell'$
 (mod $3$). We note that $(\delta s)(\gamma,\gamma')=
 [k'/3]-[(k+k')/3]+[k/3]$ is unchanged if either $k$ or $k'$ 
 is changed by
 a multiple of $3$.  So there are only nine values to examine.
 If $\ell=0$ or $\ell'=0$, or if $\ell=\ell'=1$, then $(\delta s)(\gamma,\gamma')=0$ because of the additivity of $[ \ ]$ for these values.
 But if one or both of $\ell$ or $\ell'$ is $2$ then
 $[(k+k')/3]=[k'/3]+[k/3]+1$ so $(\delta s)(\gamma,\gamma')=-1$.
As $\theta(\sigma^kh)=k\in \ZZ$ this coincides with the values of
$-(\theta \lambda^2_{\theta'}+\lambda^2_{\theta}\theta')$, as required.

  For (ii), 
as $\widehat{\theta}={\rm res}_3(\tilde{\theta})\in\{0,1,2\}\subset\Z/9\Z$,
we have $\tilde{\theta}=\widehat{{\theta}}+3s_{{\theta}}$,
and so by (i) as $\delta(\tilde{\theta})=0$
 we have
$$  \delta(\widehat{{\theta}})
=-3\delta(s_{{\theta}})=3({{\theta}}\lambda^2_{{{\theta}}'}+\lambda^2_{{\theta}}{{\theta}}').$$
For (iii), by definition
$$\delta(\lambda^2_{\tilde{\theta}})=\frac{\tilde{\theta}(\tilde{\theta}-1)}{2}+
\frac{\tilde{\theta}'(\tilde{\theta}'-1)}{2}-
\frac{(\tilde{\theta}+\tilde{\theta})'(\tilde{\theta}+\tilde{\theta}'-1)}{2}=
-\tilde{\theta}\tilde{\theta}'.$$

For (iv), one readily checks $\lambda^2_{\widehat{\theta}}={\rm res}_3(\lambda^2_{\tilde{\theta}})
\in\{0,1,2\}\subset\Z/9\Z$, and from this 
\begin{eqnarray*}\lambda^2_{\tilde{\theta}+3}&=&
\frac{(\tilde{\theta}+3)(\tilde{\theta}+2)}{2} \ = \ 
\frac{\tilde{\theta}^2+5\tilde{\theta}+6}{2}\\&=&
\frac{\tilde{\theta}(\tilde{\theta}-1)}{2}+3(\theta+1)
 \ = \ \lambda^2_{\tilde{\theta}}+3(\theta+1)
\end{eqnarray*}
and from this it follows (also by directly comparing values) that
$$\lambda^2_{\tilde{\theta}}=\lambda^2_{\widehat{\theta}}+3(\widehat{\theta}+1)s_\theta.$$
As
$\delta(\lambda^2_{\tilde{\theta}})=
-\tilde{\theta}\tilde{\theta}'$
and using $\theta^2=\theta-\lambda^2_\theta$ and $\theta\lambda^2_\theta=2\lambda^2_\theta$ (mod $3$) we find,
\begin{eqnarray*}
\delta(\theta s_\theta)&=&\theta' s_{\theta'}-(\theta+\theta') s_{\theta+\theta'}+\theta s_\theta 
 \ = \  
\theta'(s_{\theta'}- s_{\theta+\theta'})+\theta(s_{\theta}- s_{\theta+\theta'})
\\&=&\theta'(\delta(s_{\theta})- s_{\theta})+\theta(\delta(s_{\theta})- s_{\theta'})
\\&=&\theta'(-\theta\lambda_{\theta'}- \lambda_{\theta}\theta'- s_{\theta})+\theta(-\theta\lambda_{\theta'}- \lambda_{\theta}\theta'- s_{\theta'})
\\&=&-((\theta'-\lambda^2_{\theta'})\lambda^2_\theta +2\theta'\lambda^2_\theta)-
((\theta-\lambda^2_\theta)\lambda^2_{\theta'} +2\theta\lambda^2_{\theta'})
-s_\theta\theta'-\theta s_{\theta'}\\
&=&2\lambda^2_\theta\lambda^2_{\theta'}-s_\theta\theta'-\theta s_{\theta'}
\end{eqnarray*}
giving (iv). 

Using $\lambda^2_{\tilde{\theta}}=\lambda^2_{\widehat{\theta}}+3(\widehat{\theta}+1)s_\theta$ we have,
\begin{eqnarray*}\delta(\lambda^2_{\widehat{\theta}})&=&\delta(\lambda^2_{\tilde{\theta}})-3\delta((\theta+1)s_\theta)
\\&=&-\tilde{\theta}\,\tilde{\theta}'-3\left(-
\lambda^2_{\overline{\theta}}\lambda^2_{\overline{\theta}'}
-s_\theta\theta'-\theta s_{\theta'}
-\theta\lambda^2_{\theta'}-\lambda^2_\theta\theta'\right)
\\&=&-(\widehat{\theta}+3s_{\theta})(\widehat{\theta}'+3s_{\theta'})+3\left(
\lambda^2_{\widehat{\theta}}\lambda^2_{\widehat{\theta}'}
+s_\theta\theta'+\theta s_{\theta'}
+\theta\lambda^2_{\theta'}+\lambda^2_\theta\theta'\right)
\\&=&-\widehat{\theta}\,\widehat{\theta}'+3\left(
\lambda^2_{\widehat{\theta}}\lambda^2_{\widehat{\theta}'}
+\theta\lambda^2_{\theta'}+\lambda^2_\theta\theta'\right)
\end{eqnarray*}
which gives (v).  

For (vi),
$\widehat{s}_\theta$  takes the values $s_\theta\in\{0,1,2\}$, except
taken in $\Z/9\Z$.  Then using the notation in  (i) with
$\widehat{s}(\gamma)=k, \widehat{s}(\gamma)=k'\in\{0,1,2\}
\subset \Z/9\Z$ we still have
$\delta(\widehat{s})(\gamma,\gamma')=[\frac{k}{3}]-[\frac{k+k'}{3}]+
[\frac{k'}{3}]$ which takes the values $0$ or $-1$ and takes the value
$-1$ when one of $k$ or $k'$  has the value $2$ and the other is $2$ or $1$.
In (i) we noted $-(\theta\lambda^2_{\theta'}+\lambda^2_\theta\theta')$  
  (mod $3$)   takes the values $0$ or $-1$ under precisely the same conditions.
However, (mod $9$), $-(\widehat{\theta}\lambda^2_{\widehat{\theta}'}+\lambda^2_
{\widehat{\theta}}\widehat{\theta}')$ takes the value $-1$  when one of $\theta, \theta'$ takes the value $1$ and the other $2$, while it takes the value
$-4$ when both take the value $2$. To compensate we add the value
$3\lambda^2_\theta\lambda^2_{\theta'}$.  We have established
$$\delta(\widehat{s})=-(\widehat{\theta}\lambda^2_{\widehat{\theta}'}+\lambda^2_
{\widehat{\theta}}\widehat{\theta}')+3\lambda^2_\theta\lambda^2_{\theta'}.$$

Part (vii) is straight-forward and this completes the computations. $\Box$

\bs

We remark that if we 
 define $r:\Z/9\Z\rightarrow GL_5(\ZZ)$ 
for $i\in\Z/9\Z$ by 
$$r({i})=\left(\begin{array}{ccccc}1&\overline{i}&\lambda^2_{\overline{i}}&s_{\overline{i}}&\overline{i} s_{\overline{i}}\\0&1&\overline{i}&\lambda^2_{\overline{i}}&s_{\overline{i}}\\
0&0&1&\overline{i}&\lambda^2_{\overline{i}}\\0&0&0&1&\overline{i}\\0&0&0&0&1
\end{array}\right)$$
where ${\overline{i}}=i$ (mod $3)\in\ZZ$,
then one readily checks that $r$ is a homomorphism by noting 
$r(i)=r(1)^i$. Parts  (i),  (iii) and (iv) of the lemma follow by computing
$r(\theta)-r(\theta\theta')+r(\theta')$ using matrix multiplication and 
inspecting the entries.

In the next lemma the cocycle computations 
are elementary and  may be 
 in the literature somewhere but we do not 
have a convenient source.  The result relies on the fact that as 
$\G$ is a Voevodsky group, if 
$[z]\in 
{\rm ker}\left(H^n(\G,\ZZ)\rightarrow H^n({\mathcal  H}_{{\theta}=0},\ZZ)\right)$
then $[z]=[\theta]\smile[w]$ for some $[w]\in H^{n-1}(\G,\ZZ)$.
(This is classical for fields when $n=2$ using the structure of cyclic algebras 
and follows by  
Schley [S] more generally.) We recall for use
below that if  $x\in C^1(\G,\ZZ)$ 
and $y\in C^{n-1}(\G,\ZZ)$ then $\delta(x\cdot y)=x\cdot\delta(y)-x\cdot\delta(y)$.

\bs

\noi {\bf Lemma 3.6.}  {\em  Suppose that $\tilde{\theta}\in H^1(\G,\Z/9\Z)$ is a
  lift
of $\theta\in H^1(\G,\ZZ)$.\\
(i) Suppose $\delta(u_1)=-\theta u_0'$
where $u_1\in C^{n-1}(\G,\ZZ)$ and $u_0\in Z^{n-1}(\G,\ZZ)$. Then
$z_1:=\theta u'_1+\lambda^2_\theta u'_0\in Z^n(\G,\ZZ)$ and there exists
$u_2\in C^{n-1}(\G,\ZZ)$ and $u_1^*\in Z^{n-1}(\G,\ZZ)$ so that
setting $\tilde{u}_1:=u_1-u_1^*$ we have
$\delta(u_2)=-\theta \tilde{u}'_1-\lambda^2_\theta u'_0$ and   $\delta(\tilde{u}_1)=-\theta u_0'$.\\
(ii) Suppose $u_2,\tilde{u}_1\in C^{n-1}(\G,\ZZ)$, 
$\delta(u_2)=-\theta \tilde{u}'_1-\lambda_\theta u'_0$ and   $\delta(\tilde{u}_1)=-\theta u_0'$ as in the conclusion of (i). Then 
$\delta(\theta u_2'+\lambda^2_\theta\tilde{u}'_1)=
(\theta \lambda^2_{\theta'}+\lambda^2_{\theta}\theta')u_0''$, and setting   
$D:=\theta u_2'+\lambda^2_\theta\tilde{u}'_1+s_\theta u_0'$ for $s_\theta$ as in 
Definition 3.3 (v) we have
$D\in Z^n(\G,\ZZ)$.
}

\bs

\noi {\bf Proof.} Since $\delta(u_1)=-\theta u_0'$
 and Lemma 3.5 (iii) gives  $\delta(\lambda^2_{\theta})=-\theta\theta'$, we calculate
$$\delta(z_1)=
\delta(\theta u_1'+\lambda^2_{\theta} u_0')=
-\theta (-\theta' u_0'')+(-\theta\theta') u_0''
= 0
$$ showing $z_1\in Z^n(\G,\ZZ)$.  Next, as
$z_1|_{{\mathcal H}_{{\theta}=0}}=0\in Z^n({\mathcal H}_{{\theta}=0},\ZZ)$
we have $[z_1]\in 
{\rm ker}\left(H^n(\G,\ZZ)\rightarrow H^n({\mathcal  H}_{{\theta}=0},\ZZ)\right)$. 
So we know there exists $u_1^*\in H^{n-1}(\G,\ZZ)$  
with $[z_1]=-[\theta u_1^*{}']\in H^n(\G,\ZZ)$ and $u_2\in C^{n-1}(\G,\ZZ)$ with
$\delta(u_2)=-z_1-\theta u_1^*{}'=-\theta {u}'_1-\lambda^2_\theta u'_0-\theta u_1^*{}'$.
So if $\tilde{u}_1:=u_1+u_1^*$ we have
$\delta(u_2)=-\theta \tilde{u}'_1-\lambda^2_\theta u'_0$ and   $\delta(\tilde{u}_1)=\delta({u}_1)=-\theta u_0'$
giving (i).

 For (ii), as
$\delta(u_2)=-\theta \tilde{u}'_1-\lambda^2_\theta u'_0$ and   $\delta(\tilde{u}_1)=-\theta u_0'$
we find 
\begin{eqnarray*}
\delta(\theta u_2'+\lambda^2_\theta\tilde{u}'_1)&=&
-\theta(-\theta' \tilde{u}''_1-\lambda^2_{\theta'} u''_0)
+(-\theta\theta')\tilde{u}''_1-\lambda^2_{\theta}(-\theta' u_0'')
\\&=&
(\theta \lambda^2_{\theta'}+\lambda^2_{\theta}\theta')u_0''.
\end{eqnarray*}
The statement about $\delta(D)=\delta(\theta u_2'+\lambda^2_\theta\tilde{u}'_1+su_0')=0$
now follows from Lemma 3.5 (i)
since $\delta(su_0')=-(\theta \lambda^2_{\theta'}+\lambda^2_{\theta}\theta')u_0''$.  
%
This gives the result. $\Box$

 \bs
 
 We can generalize Lemma 3.6 to the case of exponent $9$.  This
 gives more information about the class $[D]\in H^n(\G,\ZZ)$
 in part (ii).

\bs

\noi {\bf Theorem 3.7.}   {\em  Suppose that 
$\tilde{\theta}\in H^1(\G,\Z/9\Z)$ is a lift
of $\theta\in H^1(\G,\ZZ)$ and the action of $\G$ on  $\Z[\langle \tau\rangle]$ 
is given by
$\G\rightarrow \G/{\rm ker}(\theta)\cong\langle\tau\rangle$ where
$\theta$'s action is multiplication by $\tau$.
We denote   $t:=\tau-1\in \ZZ[\langle \tau\rangle]$ .\\
(i) Suppose $ u_0\in Z^{n-1}(\G,\ZZ)$, $u_1,u_2\in C^{n-1}(\G,\ZZ)$
with $\delta(u_1)=-\theta u_0'$
and $\delta(u_2)=-\theta {u}'_1-\lambda^2_\theta u'_0$. Then
$$D:=\theta u_2'+\lambda^2_\theta{u}'_1+s_\theta u_0'
\in Z^n(\G,\ZZ)$$ where $s_\theta$ is as in  Definition 3.5.\\
(ii) The conditions on $u_0,u_1,u_2$ in part (i) are equivalent to
$$u_0+u_1t+u_2t^2\in Z^{n-1}(\G,\ZZ[\langle \tau\rangle])$$
(iii) Assuming $\G$ is a Voevodsky group  there exist lifts
of $u_0,u_1,u_2\in C^{n-1}(\G,\Z/3\Z)$ in parts (i) and (ii) to $\tilde{u}_0, \tilde{u}_1,\tilde{u}_2\in C^{n-1}(\G,\Z/9\Z)$ with
$$\tilde{u}_0+\tilde{u}_1t+\tilde{u}_2t^2\in Z^{n-1}(\G,\Z/9\Z[\langle \tau\rangle]).$$
(iv)
The conditions on $\tilde{u}_0, \tilde{u}_1,\tilde{u}_2$ 
in (iii) are equivalent to
\begin{eqnarray*}
\delta(\tilde{u}_0)&=&0\\
\delta(\tilde{u}_1)&=&-\widehat{\theta}{u}_0'+3\left({\theta}{u}_2'
+\lambda^2_{{\theta}}{u}'_1\right)
\\
\delta(\tilde{u}_2)&=&-\widehat{\theta}\tilde{u}_1'-
\lambda^2_{\widehat{\theta}}\tilde{u}'_0+3\left(({\theta}+\lambda^2_{{\theta}}){u}'_2
+\lambda^2_{{\theta}}{u}'_1\right)
\end{eqnarray*}
where we note that $\widehat{\theta}$ is used instead of  $\tilde{\theta}$.
\\
(v) Using  $\tilde{u}_0, \tilde{u}_1,\tilde{u}_2\in C^{n-1}(\G,\Z/9\Z)$ from (iii)
and therefore the conditions in (iv) we have
$$\tilde{D}:=\widehat{\theta}\tilde{u}'_2+\lambda^2_{\widehat{\theta}}\tilde{u}'_1
+\widehat{s}_{{\theta}}\tilde{u}'_0-3\lambda^2_\theta u'_2\in Z^n(\G,\Z/9\Z)$$
where $\widehat{s}_\theta$ is the (mod $9$) lift of $s_\theta$ to elements
in $\{0,1,2\}\in\Z/9\Z$.\\
(vi) For $\tilde{D}$ as in (v) can express $[\tilde{D}]=[\tilde{\theta}w_0']\in
H^n(\G,\Z/9\Z)$ for some $\tilde{w}_0'\in Z^{n-1}(\G,\Z/9\Z)$ and 
 for $D$ from part (i) we have $[D]
=[\theta w'_0]\in H^{n}(\G,\ZZ)$ for some $w_0\in Z^{n-1} (\G,\ZZ)$. 
}

\bs

\noi {\bf Proof.} Part (i) is already contained in Lemma 3.6. For (ii)
by Lemma 3.5 (vii) we have for $\gamma\in\langle\tau\rangle$ 
that $(\gamma-1)=\theta(\gamma)t+\lambda^2_{\theta(\gamma)}t^2
\in\ZZ[\langle\tau\rangle]$
so computing  in $C^n(\G,\ZZ[\langle\tau\rangle]$ and suppressing the
variables $\gamma\in\langle\tau\rangle$, $\gamma'\in\langle\tau\rangle^{n-1}$
we find
\begin{eqnarray*}
\delta(u_0+u_1t+u_2t^2)&=&\delta(u_0)+\theta u_0't+\lambda^2_\theta u_0't^2
+\delta(u_1)t+\theta u_1't^2+\delta(u_2)t^2\\
&=&\delta(u_0)+\left(\delta(u_1)+\theta u_0'\right)t+\left(
\delta(u_2)+\theta u_1't^2+\lambda^2_\theta u_0'\right)t^2.
\end{eqnarray*}
From this we see $\delta(u_0+u_1t+u_2t^2)=0\in
C^n(\G,\ZZ[\langle\tau\rangle]$ if an only if the conditions on 
$u_0$, $u_1$, $u_2$ in part (i) are satisfied.

For 
(iii), the Voevodsky assumption on $\G$ and $\Z[\langle\tau\rangle]$
shows there exist $\tilde{w}_0$, $\tilde{w}_1$, $\tilde{w}_2
\in C^{n-1}(\G,\Z/9\Z)$ with
$\tilde{w}_0+\tilde{w}_1t+\tilde{w}_2t^2\in Z^{n-1}(\G,\Z/9\Z[\langle \tau\rangle])$
and 
$[\overline{\tilde{w}_0}+\overline{\tilde{w}_1}t+\overline{\tilde{w}_2t^2}]
=[u_0+u_1t+u_2t^2]\in H^{n-1}(\G,\Z/3\Z[\langle \tau\rangle])$
where $\overline{\tilde{w}_k}$ is the image of ${\tilde{w}_k}$ induced
under the projection  $\Z/9\Z\rightarrow\ZZ$. So there exists
$x_0+x_1t+x_2t^2 \in C^{n-2}(\G,\Z/3\Z[\langle \tau\rangle])$ with
$\delta(x_0+x_1t+x_2t^2)=(\overline{\tilde{w}_0}-u_0)+
(\overline{\tilde{w}_1}-u_1)t+(\overline{\tilde{w}_2}-u_2)t^2\in Z^{n-1}
(\G,\Z/3\Z[\langle \tau\rangle])$. Now, taking an arbitrary lift of
$x_i$ to $\tilde{x}_i\in C^{n-2}(\G,\ZZ[\langle\tau\rangle])$ we can
then define the $\tilde{u}_i$ by
$$\tilde{u}_0+\tilde{u}_1t+\tilde{u}_2t^2:=
(\tilde{w}_0+\tilde{w}_1t+\tilde{w}_2t^2)-\delta(\tilde{x}_0+\tilde{x}_1t+\tilde{x}_2t^2)
$$
and we now have the desired elements
with $\overline{\tilde{u}_i}=u_i$.

For  (iv) we generalize the computation in (ii), but now 
using $(\gamma-1)=\widehat{\theta}(\gamma)t+\lambda^2_{\widehat{\theta}(\gamma)}t^2\in \Z/9\Z[\langle\tau\rangle]$
where for $\tilde{u}=\tilde{u}_0+\tilde{u}_1t+\tilde{u}_2t^2$ we find in 
$C^n(\G,\Z/9\Z[\langle\tau\rangle])$ (again supressing $\gamma$ and $\gamma'$),
\begin{eqnarray*}
\delta(\tilde{u})&=&
\delta(\tilde{u}_0)+\left(\delta(\tilde{u}_1)+\widehat{\theta}\tilde{u}_0'
-3\left( \theta u_2'+\lambda^2_\theta u'_1\right)\right)t
\\&&+\left(\delta(\tilde{u}_2)+\widehat{\theta}\tilde{u}'_1+\lambda^2_{\widehat{\theta}}\tilde{u}'_0
-3\left( (\theta+\lambda^2_\theta)u_2'+\lambda^2_\theta u'_1\right)
\right)t^2
\end{eqnarray*}
from which (iv) follows.

For (v), by Lemma 3.5 (v),
$\delta(\lambda^2_{\widehat{\theta}})=-\widehat{\theta}\,\widehat{\theta}'+3\left(
\lambda^2_{\widehat{\theta}}\lambda^2_{\widehat{\theta}'}
+\theta\lambda^2_{\theta'}+\lambda^2_\theta\theta'\right)$
and by (vi),
$\delta(\widehat{s})=-(\widehat{\theta}\lambda^2_{\widehat{\theta}'}+\lambda^2_
{\widehat{\theta}}\widehat{\theta}')+3\lambda^2_\theta\lambda^2_{\theta'}$,
so we can   compute
\begin{eqnarray*}
\delta(\widehat{\theta}\tilde{u}_2)&=&3\left(\theta\lambda^2_{\theta'}+\lambda^2_\theta\theta')\right)\tilde{u}''_2\\&&-\widehat{\theta}\left(
-\widehat{\theta}'\tilde{u}_1''-
\lambda^2_{\widehat{\theta}'}\tilde{u}''_0 
+3\left(({\theta}'+\lambda^2_{{\theta}'}){u}''_2
+\lambda^2_{{\theta}'}{u}''_1
\right)\right)  
\\
&=& \widehat{\theta}
\widehat{\theta}'\tilde{u}_1''+\widehat{\theta}
\lambda^2_{\widehat{\theta}'}\tilde{u}''_0 +3\left((-\theta\theta'+\lambda^2_\theta\theta')u_2''-\theta\lambda^2_{{\theta}'}{u}''_1
\right) 
\\
\delta(\lambda^2_{\widehat{\theta}}\tilde{u}_1)&=&  \left(-\widehat{\theta}\widehat{\theta'}
+3(\lambda^2_{\theta}\lambda^2_{\theta'}
+\theta\lambda^2_{\theta'}+\lambda^2_\theta\theta')\right)\tilde{u}''_1\\&& \ -\lambda^2_{\widehat{\theta}}\left(
-\widehat{\theta'}\tilde{u}_0''+3\left({\theta'}{u}_2''
+\lambda^2_{{\theta'}}{u}''_1\right)
\right)\\
&=&-\widehat{\theta}\widehat{\theta'}
\tilde{u}_1''+\lambda^2_{\widehat{\theta}}
\widehat{\theta'}\tilde{u}_0''+3\left(-\lambda^2_{\theta}\theta'u_2''
+(\theta\lambda^2_{\theta'}+\lambda^2_\theta\theta')\tilde{u}''_1\right)
\\
\delta(\widehat{s}\tilde{u}_{0}')&=&-(\widehat{\theta}\lambda^2_{\widehat{\theta}'}+\lambda^2_
{\widehat{\theta}}\widehat{\theta}')\tilde{u}_{0}''+3\lambda^2_\theta\lambda^2_{\theta'}\tilde{u}_{0}''\\
\delta(-3\lambda^2_\theta u_2')&=&3\left(\theta\theta'u_2''-\lambda^2_\theta\theta'u_1''-
\lambda^2_\theta\lambda^2_{\theta'}u_0''\right)
\end{eqnarray*}
The sum of these four is zero, giving (v).

For (vi) the first statement follows as $\tilde{D}=0$ when restricted
to ${\mathcal H}_{\tilde{\theta}=0}$. 
For the second statement we note that by the definitions of the liftings
that $D=\pi^*(\tilde{D})\in C^2(\G,\ZZ)$ where
$\pi:\Z/9\Z\rightarrow\ZZ$ is the projection.
Then expressing
$[\tilde{D}]=[\tilde{\theta}\tilde{w}_0']\in
H^n(\G,\Z/9\Z)$ for some $\tilde{w}_0'\in Z^{n-1}(\G,\Z/9\Z)$, if
${w}_0'\in Z^{1}(\G,\ZZ)$ is defined by
$w_0(\gamma')=\tilde{w}_0(\gamma')$ (mod $3$) for $\gamma'\in\G$
we have $[D]=\pi^*[\tilde{D}]=\pi^*[\tilde{\theta}\tilde{w}_0']
=[{\theta}{w}_0']\in
H^n(\G,\Z/3\Z)$.
This gives the result. $\Box$

\bs

\noi  {\bf Remarks 3.8.} (i) The occurrences of $\widehat{\theta}$ in Theorem
3.7 could be replaced by the Frobenius $\Phi(z):=z^3$
using $\widehat{\theta}=\Phi(\tilde{\theta})+3\lambda^2_\theta$,
$\lambda^2_{\widehat{\theta}}=
\Phi(\lambda^2_{\tilde{\theta}})+3\theta s_\theta$ and $\widehat{s_\theta}=
\Phi(s_{\tilde{\theta}})+3\lambda^2_{s_{\theta}}$.  We have chosen to use
$\widehat{\theta}$ instead of $\Phi(\tilde{\theta})$ 
because it is simpler when applying Lemma 3.5 (vii).

(ii) Regarding the conditions in (iv) one can verify $\delta\circ\delta=0$ on these elements. As noted in Section 1, Shapiro's Lemma
gives  $H^n(\G,\ZZ[\langle\tau\rangle])\cong
H^n(\HH,\ZZ)$ where $\HH={\rm ker}(\theta)$. 
When $n=1$ if we fix $\gamma\in\G$ so that $\gamma\HH=\tau\in\G/\HH$
we then have an explicit isomorphism ${\mathcal T}:
H^1(\HH,\ZZ)\rightarrow H^1(\G,\ZZ[\langle\tau\rangle])$  given by
$$
{\mathcal T}(\sigma)(\gamma^ih):=\sum_{j=0}^2\sigma(h^{\gamma^{-j}})\tau^{i+j}$$
where $h^{\gamma}=\gamma h\gamma^{-1}$.
An inverse ${\mathcal U}$ to ${\mathcal T}$ is given  for $\tilde{\sigma}\in
Z^1(\G,\ZZ[\langle\tau\rangle])$ by decomposing
$\tilde{\sigma}(h)=\sum_{i=0}^2\sigma_{i}(h)\tau^i$
for $\sigma_{i}(h)\in C^1(\G,\ZZ)$ and then setting
${\mathcal U}(\tilde{\sigma})(h)=\sigma_0(h)
$.  With this identification, the action of $\G/\HH\cong\langle\tau\rangle$ on 
$H^1(\G,\ZZ[\langle\tau\rangle])$ can be understood as
an action on $H^1(\HH,\ZZ)$ as follows: For  $\sigma\in Z^1(\HH,\ZZ)$
we can evaluate
$$
\left(\tau\cdot{\mathcal T}(\sigma)\right)(\tau^ih)=\sum_{j=0}^2\sigma(h^{\gamma^{-j}})
\tau^{i+j+1}
=\sum_{j=0}^2\sigma((h^{\gamma^{-1}})^{\gamma^{-j}})\tau^{i+j}.$$
From this, via the identification given by ${\mathcal T}$ we see that $\tau$ acts on 
$H^1(\HH,\ZZ)$ by $(\tau\cdot\sigma)(h)=\sigma(h^{-\gamma})$.
This action is the same as the 
action of $\G/\HH$ on $H^1(\HH,\ZZ)$ in the 
Hochschild-Serre
Spectral Sequence associated to the short exact sequence
$0\rightarrow{\mathcal H}\rightarrow \G\rightarrow\ZZ\rightarrow0$, 
where $E_2^{p,q}=H^p(\G/{\mathcal H},H^q({\mathcal H},\ZZ))\Rightarrow
H^{p+q}(\G,\ZZ)$.  In the special case that
$\tau\mapsto\sigma\in H^1(\HH,\ZZ)$ defines an element
of $H^1(\G/\HH,H^1(\HH,\ZZ))$, as $\G/\HH\cong\langle\tau\rangle$  we would have that
$(1+\tau+\tau^2)\cdot\sigma=0\in H^1(\HH,\ZZ)$, which of course is not
always the case.  In the special case where $\sigma$ is the restriction of
$\tilde{\sigma}\in Z^1(\G,\ZZ)$ this is true where it defines the element
$\tilde{\sigma}t_2\in H^1(G,\ZZ[\langle\tau\rangle])$. Then the class of
$\tau\mapsto\tilde{\sigma}t^2\in H^1(\G/\HH,H^1(\HH,\ZZ))$ can be seen
to correspond to the element $\theta\smile\tilde{\sigma}\in H^2(\G,\ZZ)$ via
the spectral sequence.  Unfortunately, outside of this case, in the  general situation the Spectral Sequence does
not seem to yield information to clarify parts (i) and (ii) of Theorem 3.7 where
an interpretation of the correspondence $H^1(\G,\ZZ[\langle\tau\rangle])
\rightarrow  H^2(\G,\ZZ)$ given by
$[ u_0+u_1t+u_2t^2]\mapsto[\theta u_2'+\lambda^2
_\theta u_1'+s_\theta u_0']$
is still lacking.

\bs

We continue to study $\eta:H^{n-1}(\G,\M_4\otimes_\Z\ZZ)\rightarrow
H^{n}(\G,\ZZ)$, whose image in the field context when $n=2$ is the relative
Brauer group ${\rm Br}_3(E/F)$. We start with $[c]\in H^{n-1}(\G,\M_4\otimes_\Z\ZZ)$ 
where, by abuse of notation, we view $c\in C^{n-1}(\G,\M_3\otimes_\Z\ZZ)$,
that is, $c$ is a pullback of an element representing the equivalence class $[c]$.
The condition that $[c]\in H^{n-1}(\G,\M_4\otimes_\Z\ZZ)$ means that for some
$\omega\in C^{n}(\G,\M_2\otimes_\Z\ZZ)$ we have 
$\delta(c)=d_2(\omega)$.  In the computation we will split
$\omega=\omega_1+\omega_2$ into two pieces; only $\omega_2$
will contribute a nonzero term to  $\eta([c])$.
We recall by Corollary 2.4 we have an exact sequence
$$0\longrightarrow (\ZZ[G]\cdot\alpha_1\oplus\ZZ[G]\alpha_2)\longrightarrow
{\mathcal M}_4\otimes_\Z\ZZ\stackrel{\pi_4}{\longrightarrow}
\ZZ[G]\cdot T_G\longrightarrow0$$
with associated map $\pi_4$  given by
$\pi_4([\eta_1 t_1,\eta_2 t_2,\eta_3 T_3,\eta_4 T_4)]=
(\eta_3-\eta_4) T_G$ and the induced map
$\pi^*_4:H^{n-1}(\G,\M_4\otimes_\Z\ZZ)\rightarrow H^{n-1}(\G,\ZZ[G]\cdot T_G)
\cong H^{n-1}(\G,\ZZ)$.
 We suppose
that $\pi_4^*([c])=\chi\cdot T_G\in H^{n-1}(\G,\ZZ\cdot T_G)$ 
for a cohomology class $\chi\in H^{n-1}(\G, \ZZ)$. We have the following lemma.

\bs

\noi {\bf Lemma 3.9.} Assume $[c]\in H^{n-1}(\G,\M_4\otimes_\Z\ZZ)$ with
$\pi_4^*([c])=\chi  T_G$. Without loss of generality we can represent $[c]$ using  
$c=(u,v,0,-\chi T_4)\in C^{n-1}(\G,M_3\otimes_\Z\ZZ)$
for  some $u,v\in C^{n-1}(\G,\ZZ[G])$.

\bs

\noi {\bf Proof.}  If  $c=(w_1,w_2,w_3T_3,w_4T_4)\in 
C^{n-1}(\G,\M_3\otimes_\Z\ZZ)$
represents $[c]$ then we set $c_2:=c-d_2(w_3,0)=
(w_1+w_3t_1,w_2+w_3t_2,0,(w_4-w_3)T_3)\in
C^{n-1}(\G,\M_3\otimes_\Z\ZZ)$.  It follows 
$[c]=[c_2]$ as they agree as functions valued
in $\M_4\otimes_\Z\ZZ$. Relabeling, $c=(w_5,w_6,0,w_7T_4)$. We can express $w_7=r_0+r_1\tau_3+r_2\tau_3^2=
(r_0+r_1+r_2)-r_1(1-\tau_3)-r_2(1-\tau_3^2)$ for $r_i\in
C^{n-1}(\G,\ZZ)$ and then
if $r:=-r_1(1-\tau_3)-r_2(1-\tau_3^2)$ we set
\begin{eqnarray*}
c_3&:=&c-d_2(r,0)\\&=&
(w_5-rt_1,w_6-rt_2,0,(r_0+r_1+r_2)T_4)\in C^{n-1}(\G,\M_3/\otimes_\Z\ZZ)
\end{eqnarray*}
where again
 $[c]=[c_2]=[c_3]$.
Relabeling yet again, we can express $c=(u,v,0,w_8T_4)$ where $w_8\in C^{n-1}(\G,\ZZ)$.
Since $[c]\in H^{n-1}(\G,\M_4\otimes_\Z\ZZ)$ we know that, as functions valued
in $\M_3$ we have $$\delta(c)=(\delta(u),\delta(v),0,
\left((\gamma-1)s'_8+\delta(w_8)\right)T_4)\in
d_2(C^{n-1}(\G,\M_2\otimes_\Z\ZZ)).$$
However, inspecting the basis ${\mathcal D}=\cup_{j=1}^6{\mathcal D}_j$ for $d_2(\M_2)$ given in Corollary 2.4 (i) we see that
only ${\mathcal D}_4, {\mathcal D}_5$ and ${\mathcal D}_6$ contribute nonzero
third or fourth components, and of those only ${\mathcal D}_6$ can
contribute an entry outside of  multiples of $t_4T_3$ in the third component
or of multiples of $t_3T_4$ in the fourth component.  
   In particular, as the third component of
 $\delta(c)$ is zero, the fourth component of $\delta(c)$ 
 is congruent to $\delta(w_8)T_4$ (mod $t_3T_4$) and as 
 ${\mathcal D}_6$ contributes the same multiple of $T_3$ and $T_4$ in these
 components,  this forces $\delta(w_8)=0$.  Applying $\pi^*_4$ we find $w_8=-\chi\in H^{n-1}(\G,\ZZ)$, which gives
 the lemma. $\Box$
 
\bs

We proceed with $c=(u,v,0,-\chi\cdot T_4)\in C^{n-1}(\G,\M_3\otimes_\Z\ZZ)$.
We calculate  in $C^n(\G,\M_3\otimes_\Z\ZZ)$, 
\begin{eqnarray*}\delta\left(c\right)(\gamma,\gamma')&=&
(\delta(u)(\gamma,\gamma'),\delta(v)(\gamma,\gamma'),
0,-(\gamma-1)\chi(\gamma')T_4)
\\&=&(\delta(u)(\gamma,\gamma'),\delta(v)(\gamma,\gamma'),0,-(\tau_3^{\psi(\gamma)}-1)\chi(\gamma')T_4)
\\&=&(\delta(u)(\gamma,\gamma'),\delta(v)(\gamma,\gamma'),0,-(\psi(\gamma)t_3
+\lambda^2_{\psi(\gamma)}T_3)\chi(\gamma')T_4)
\end{eqnarray*}
where $t_3=\tau_1\tau_2-1$, $T_3=t^2_3$,
and $\psi(\tau_1^i\tau_2^j)=\frac{1}{2}(i+j)=-(i+j)\in\ZZ$.  

We next set  
$$\omega_1:=(-(\psi(\gamma)t_3\
+\lambda^2_{\psi(\gamma)}T_3)\chi(\gamma'),0)\in 
C^{n-1}(\G,\M_2\otimes_\Z\ZZ)$$ 
and then 
$$d_2(\omega_1)(\gamma,\gamma')=-\left(
\psi(\gamma)\chi(\gamma')t_3+
\lambda^2_{\psi(\gamma)}\chi(\gamma')T_3\right)(-t_1,-t_2,T_3,T_4).
$$
In $\ZZ[G]$ we have
$t_3=\tau_1\tau_2-1=t_1+t_2+t_1t_2$ and $T_3=t^2_3=(t_1+t_2+t_1t_2)^2$,
$t^3_3=t_3T_3=0$. So
  suppressing the variables $\gamma,\gamma'$  we have
   \begin{eqnarray*}d_2(\omega_1)&=&-\left(
\psi\chi't_3+
\lambda^2_{\psi}\chi'T_3\right)(-t_1,-t_2,T_3,T_4)\\
&=&(\psi\cdot\chi')((t_1t_2+T_1+T_1t_2),(
t_1t_2+T_2+t_1T_2),0,-t_3T_4)\\&& +\,
(\lambda^2_{\psi}\cdot\chi')(-t_1T_2+T_1t_2-T_G,t_1T_2-T_1t_2-T_G,0,-T_G).
\end{eqnarray*}

We set $\M_{3;00}:=(\ZZ[G]t_1T_2,\ZZ[G]T_1t_2,0,0)\subset\M_3$.
Summarizing the above we have the following.

\bs

\noi {\bf Formula 3.10.} {\em 
Suppose $[c]=[(u,v,0,-\chi T_4)]\in H^{n-1}(\G,\M_4\otimes_\Z\ZZ)$ where $u,v\in
C^{n-1}(\G,\ZZ[G])$ and $\chi\in Z^{n-1}(\G,\ZZ)$. 
Suppressing the variables $\gamma$ and $\gamma'$ we set
 $\omega_1:=(-(\psi t_3+
\lambda_{\psi}T_3)\cdot \chi',0)\in C^n(\G,\M_2\otimes_\Z\ZZ)$. Then
\begin{eqnarray*}\delta\left(c\right)&=&(\delta(u),\delta(v),0,-(\psi t_3
+\lambda^2_{\psi}T_3)\chi'T_4),\\
d_2(\omega_1)&=&\left(
\psi\chi't_1t_3+
\lambda^2_{\psi}\chi't_1T_3\right)(t_1,t_2,0,-T_4) \ .
\end{eqnarray*}
Therefore, setting 
\begin{eqnarray*}T_\psi&:=&(t_1t_2+T_1+T_1t_2,t_1t_2+T_2+t_1T_2,0,0) \ \ {\rm and}\\
T_{\lambda^2_\psi}&:=&(t_1T_2-T_1t_2+T_G,-t_1T_2,+T_1t_2+T_G,0,0) \in \M_{3;00}
\end{eqnarray*}
we find
$$\delta(c)-d_2(\omega_1)=(\delta(u),\delta(v),0,0)-
\psi\chi'T_\psi-\lambda^2_{\psi}
\chi'T_{\lambda^2_\psi}\in C^n(\G, \M_{3;00}\otimes_\Z\ZZ) \ .
$$}

We now turn to finding 
$u,v\in C^{n-1}(\G,\ZZ[G])$  such  that
if $[c]=[(u,v,0,-\chi T_4)]\in Z^{n-1}(\G,\M_4\otimes_\Z\ZZ)$ we can find  $\omega\in C^n(\G,\M_2\otimes_\Z\ZZ)$
with $\delta(c)=d_2(\omega)$.
To find $\omega_2$ so that $\omega:=\omega_1+\omega_2$ works
we need the following result describing $d_2(\M_2)\cap\M_{3;00}$.

\bs

\noi {\bf Theorem 3.11.} {\em Suppose 
$[c]\in H^{n-1}(\G,\M_4\otimes_\Z\ZZ)$ and $\pi_4^*([c])=\chi T_G$ so
 we may assume $c=(u,v,0,-\chi T_4)\in C^{n-1}(\G,\M_3\otimes_\Z\ZZ)$. Suppose also
 $\omega_1$ is as in Formula 3.10 so
 $\delta(c)-d_2(\omega_1)\in C^n(\G,\M_{3;00}\otimes_\Z\ZZ)$.
  Then\\
 (i) We have (working over $\Z$)
\begin{eqnarray*}d_2(\M_2)\cap\M_{3;00}&=&
(\Z[G]t_1T_2,\Z[G]T_1t_2,0,0)
\\&=&\{d_2(\eta_1T_1+\eta_2T_2,-(\eta_1+\eta_2)T_G)|
\eta_1,\eta_2\in\Z[G]\}
\end{eqnarray*}
with the same result $\otimes_\Z\Z/9\Z$ and $\otimes_\Z\ZZ$.\\
(ii) There exist $\omega_2\in C^{n-1}(\M_2\otimes_\Z\ZZ)$
and $A_0,A_1,B_0,B_1\in C^{n}(\G,\ZZ)$  with
\begin{eqnarray*}
d_2(\omega_2)&=&\delta(c)-d_2(\omega_1)\\
&=&(-A_0t_1T_2-A_1 T_G,
-B_0T_1t_2-B_1T_G,0,0).
\end{eqnarray*}
(iii) For $A_0$ and $B_0$ as in (ii) we have
$A_0+B_0\in Z^n(\G,\ZZ)$ and 
$\eta([c])=[A_0+B_0]\in H^n(\G,\ZZ)$.
}

\bs

\noi {\bf Proof.}
For (i), by definition  $d_2(T_2,-T_G)=(-t_1T_2,0,0,0)$ and $d_2(T_1,-T_G)=(0,-T_1t_2,0,0)$
we have the first $\supseteq$ inclusion.  Next, if $d_2(\eta,nT_G)$ is in the
stated intersection, then $\eta T_3=-nT_G$ and $\eta T_4=-nT_G$ so that
$\eta(T_3-T_4)=0$. But $T_3-T_4=\tau_1\tau_2(\tau_1-1)(\tau_2-1)$
so $\eta(\tau_1-1)(\tau_2-1)=\eta t_1t_2=0$.  
By Lemma 2.1 (iv) we have an exact sequence
$$0\rightarrow\Z[G]\cdot\{T_1,T_2\}\rightarrow\Z[G]\rightarrow\Z[G]\cdot t_1t_2\rightarrow0$$
so it follows that $\eta\in\Z[G]\cdot\{T_1,T_2\}$ and from this we have
$d_2(\eta,nT_G)=$ $(-t_1\eta,-t_2\eta,0,0)\in
(\Z[G]t_1T_2,\Z[G]T_1t_2,0,0)$ giving the first equality.
 The second equality
follows from the proof of the first. Passing to
$\otimes_\Z\Z/9\Z$ and $\otimes_\Z\ZZ$ is clear.

For (ii), as $[c]\in H^{n-1}(\G,\M_4\otimes_\Z\ZZ)$, 
$\delta(c)\in d_2(C^{n}(\G,\M_2\otimes_\Z\ZZ))$.  
By Formula 3.10 we know
that $\delta(c)-d_2(\omega_1)$ is valued in $d_2(\M_2)\cap\M_{3;00}$.
We note that $\ZZ[G]t_1T_2=\ZZ\cdot t_1T_2\oplus \ZZ\cdot T_G$ and 
$\ZZ[G]T_1t_2=\ZZ\cdot T_1t_2\oplus \ZZ \cdot T_G$ 
where each are $2$-dimensional
as $\ZZ$-modules.  
So by part (i) the desired $A_0,A_1,B_0,B_1\in C^n(\G,\ZZ)$ exist.

For (iii),  we have 
$d_2(\omega_1+\omega_2)=
\delta(c)=\delta(u,v,0,-\chi T_4)$ where $\pi_4^*(c)=\chi T_G$.  
By the Positselski [P] construction we have that
$\eta([c])=[h_1(\omega_1+\omega_2)]$ where 
$h_1(\xi,n\cdot T_G)=(\xi-n)\cdot T_G$.
We have $\omega_1\in (C^{n-1}(\G,\ZZ[G](\tau_1\tau_2-1)),0)$ so
that  $h_1(\omega_1)=0$
and as $\omega_2=((A_0+A_1t_1)T_2+(B_0+t_1B_1)T_1,-(A_0+B_0)T_G)
$  with $T_1T_4=T_2T_4=0$ we find $h_1(\omega_1+\omega_2)=A_0+B_0$. Therefore $A_0+B_0$  must lie in 
$Z^n(\G,\ZZ)$ with $\eta([c])= [A_0+B_0]\in H^n(\G,\ZZ)$.
The result follows.   $\Box$

\bs

The next task is to find 
$\omega_2\in C^{n-1}(\G,\M_2\otimes_\Z\ZZ)$ with $$d_2(\omega_2)=\delta(c)-d_2(\omega_1)\in C^n(\G,\M_{3;00}\otimes_\Z\ZZ)$$
and $u,v\in C^1(\G,\ZZ[G]\otimes_\Z\ZZ)$
for which we have a 
cocycle $[c]=[(u,v,0,-\chi T_G)]$ in $H^{n-1}(\G,\M_4)$ 
with $\eta([c])= [A_0+B_0]\in H^n(\G,\ZZ)$.

 For this we consider  
$(u,v,0,0)\in C^{n-1}(\G,\M_3\otimes_\Z\ZZ)$ and compute $\delta(u,v,0,0)$. We  
 set $$u=\sum_{0\leq \ell,m\leq2}
u_{\ell m}t_1^\ell t_2^m\in C^{n-1}(\G,\ZZ[G])$$ where  the
$u_{\ell m}\in C^{n-1}(\G,\ZZ)$.  For $k=1,2$ we 
recall $T_k:=t_k^2$ and  $t_k^3=0\in\ZZ[G]$, and  
 denoting $\Lambda_k(\gamma):=
\chi_k(\gamma)t_k+\lambda^2_{\chi_k(\gamma)}T_k$,  by  
direct calculation we generalize Lemma 3.5 (vii) and for $\gamma\in G$ we have
$$(\gamma-1)=\Lambda_1+\Lambda_2+\Lambda_1\Lambda_2\in \ZZ[G].$$
From this we find that
\begin{eqnarray*}
\delta(u)&=&\sum_{0\leq \ell,m\leq2}
\delta(u_{\ell m})t_1^\ell t_2^m+(\gamma-1)
\sum_{0\leq \ell,m\leq2}
u_{\ell m}'t_1^\ell t_2^m\\
&=&\sum_{0\leq \ell,m\leq2}
\delta(u_{\ell m})t_1^\ell t_2^m+(\Lambda_1+\Lambda_2+\Lambda_1\Lambda_2)
\sum_{0\leq \ell,m\leq2}
u'_{\ell m}t_1^\ell t_2^m.
\end{eqnarray*}
To simplify these expressions we abbreviate
$\chi_1(\gamma):=i$ and $\chi_2(\gamma):=j$, and we 
denote by $\lambda^1_i:=i$ and  $\lambda^0_i:=1$.
Collecting the  $t_1^\ell t_2^m$ summands we find that
$$\delta(u)=\sum_{0\leq \ell,m\leq2}
\left(\delta(u_{\ell m})+\sum_{
\begin{array}{c}0\leq r\leq \ell\\
0\leq s\leq m\end{array}
}\lambda^{\ell-r}_i\lambda^{m-s}_ju'_{r\,s}
\right)t_1^\ell t_2^m \ .
$$

To find    $u$, $v$ with
$\delta(u,v,0,0)- d_2(\omega_1)\in d_2(\M_2)\cap C^2(\G,M_{3,00}\otimes_\Z\ZZ$))
we start with finding possible $u$, where the case of $v$
can be handled by symmetry. Formula 3.10 and Theorem 3.11 (i) show that
$$\delta(u)-\psi\chi'(t_1t_2+t_1^2+t_1t_2^2)-\lambda^2_\psi\chi'
(t_1t_2^2-t_1^2t_2+t_1^2t_2^2)\in C^n(\G,\ZZ[G]t_1t_2^2)$$
so the computation of $\delta(u)$ just given  gives the system
  $(*_u)$:
\begin{eqnarray*}
\delta(u_{00})&=&0\\
\delta(u_{10})&=&-iu_{00}'\\
\delta(u_{01})&=&-ju_{00}'\\
\delta(u_{20})&=&-(iu_{10}'+\lambda_iu'_{00})+\psi\chi'\\
\delta(u_{02})&=&-(ju_{01}'+\lambda_ju_{00}')\\
\delta(u_{11})&=&-(ju_{10}'+iu_{01}+iju_{00})+\psi\chi'\\
\delta(u_{21})&=&-(ju_{20}'+iu'_{11}+iju'_{10}+\lambda_iu_{01}'
+\lambda_iju_{00}')+(\psi-\lambda_\psi)\chi'\\
\delta(u_{12})&:=&-(iu_{02}'+ju'_{11}+iju_{01}'+\lambda_ju_{10}'
+i\lambda_ju_{00}')+\lambda_{\psi}\chi'+A_0
\\
\delta(u_{22})&:=&-(\lambda_iu_{02}'+\lambda_ju_{20}'+iju'_{11}+iu_{12}+ju_{21}\\
&& \ +\lambda_iju_{01}'+i\lambda_ju_{10}'
+\lambda_i\lambda_ju_{00}' )+\lambda_{\psi}\chi'+A_1.
\end{eqnarray*}
Here, to find possible values of $u$ 
we need to solve the first seven equations for the functions
$u_{00}$, $u_{10}$, $u_{01}$, $u_{20}$, $u_{11}$, $u_{02}$, $u_{21}$, 
and then using these,  
the values of $A_0$ and $A_1$ 
are determined 
(mod $\delta(C^{n-1}(\G,\ZZ)$)).
Since all we need is the value of $A_0$ (mod $\delta(C^1(\G,\ZZ)))$
we see that the goal is to describe the solutions to the first
seven equations of $(*_u)$ and use these to obtain 
an expression for $A_0$.
We have the
analogous equations for $v_{ij}$, which we will denote as $(*_v)$,
where since $t_1$ and $t_2$ are reversed, the $i$ and $j$ are also
reversed, and we will use $B_0$ and $B_1$ in place of $A_0$
and $A_1$. 

\bs



To analyze $A_0$ 
we 
 %
apply Lemma 3.6 and Theorem 3.7 three times
using the equations in $(*_u)$ to obtain  three 
cohomology classes,  $D_1$, $D_2$
and $D_3$ as follows.
\begin{eqnarray*}
\delta(u_{10})&=&-iu'_{00}\\
\delta(u_{20})&=&-(iu'_{10}+\lambda_iu_{00})+\psi\chi' 
\ \ \Rightarrow\\
\exists y_1,\chi_3 \ \ \delta(y_1)&=&\psi\chi'-i\chi_3'
\ \ \Rightarrow\\
\bullet_1 \ \ \delta(u_{20}-y_1)&=&-i(u_{10}-\chi_3)'-\lambda_iu_{00}
\\ \bullet_2 \ \ 
\delta(u_{10}-\chi_3)&=&-iu'_{00}\\
 \bullet_1 \ , \ \bullet_2 \ 
{\rm give:} \ \ D_1&=&
i(u_{20}-y_1)'+\lambda_i(u_{10}-\chi_3)'+s_iu_{00}'
\end{eqnarray*}
\begin{eqnarray*}
\bullet_3 \ \ \delta(u_{01})&=&-ju'_{00}\\
\bullet_4 \ \  \delta(u_{02})&=&-(ju'_{01}+\lambda_ju_{00})
\\ \bullet_3 \ , \ \bullet_4 \ 
{\rm give:} \ \ D_2&=&
ju_{02}'+\lambda_ju_{01}'+s_ju_{00}'
\end{eqnarray*}
Using $\lambda_i-\lambda_{i+j}+\lambda_j=-ij$ 
and $\psi=-(i+j)$ we find
\begin{eqnarray*}
\bullet_5 \ \ \delta(u_{20}+u_{02}+u_{11})&=&-(i+j)(u_{10}+u_{01})'\\&& \ \ -(\lambda_i+\lambda_j+ij)u_{00}'+2\psi\chi'
\\
&=&-(i+j)(u_{10}+u_{01}-\chi)'-\lambda_{i+j}u_{00}'\\
\bullet_6 \ \ \delta(u_{10}+u_{01}-\chi)&=&-(i+j)u_{00}\\
\bullet_5 \ , \ \bullet_6 \ 
{\rm give:} \ \ \ \ D_3&=&
(i+j)(u_{20}+u_{02}+u_{11})'\\&& \ \ +\lambda_{i+j}(u_{10}+u_{01}
-\chi)'+s_{i+j}u_{00}'.
\end{eqnarray*}

\bs

We remark a fourth cocycle $D_4$ can be obtained as follows:
\begin{eqnarray*}
\bullet_7 \ \ \delta(u_{01}+u_{20}+u_{02}-u_{11})&=&-(i-j)(u_{10}-u_{01})' -(\lambda_{i-j})u_{00}'
\\
\bullet_8 \ \ \delta(u_{10}-u_{01})&=&-(i-j)u_{00}\\
\bullet_7 \ , \ \bullet_8 \ 
{\rm give:} \ \ \ \ D_4&=&
(i-j)(u_{01}+u_{20}+u_{02}-u_{11})'\\&& \ \ +\lambda_{i-j}(u_{10}-u_{01}
-\chi)'+s_{i-j}u_{00}'.
\end{eqnarray*}
This cocycle is of interest because of the calculations of the
Bockstein in Section 4.

\bs

We recall that the $3$-$3$ subgroup
of $\G$ is given as ${\mathcal H}:={\rm ker}(\chi_1)\cap{\rm ker}(\chi_2)\subset\G$ for the two characters $\chi_1$, $\chi_2\in H^1(\G,\ZZ)$
and in analogy with Brauer group theory we denote
$${\rm  Dec}^n({\mathcal G},\chi_1,\chi_2)
:=\chi_1\smile H^{n-1}(\G,\ZZ)+\chi_2\smile H^{n-1}(\G,\ZZ)\subseteq H^n(\G,\ZZ).$$
Combining the $D_i$ gives the following.

\bs

\noi {\bf Lemma 3.12} {\em We have
\begin{eqnarray*}
D_1+D_2- D_3&=&
-(ju_{20}'+iu_{02}'+(i+j)u_{11}') -(\lambda^2_j+ij)u_{10}'-(\lambda^2_i+ij)u_{01}'\\&& \ \ 
-(i\lambda^2_j+\lambda^2_ij)u_{00}'-
(iy_1'+\lambda^2_i\chi_3')+\lambda^2_{i+j}\chi'
\end{eqnarray*}
where $[D_1+D_2- D_3]
\in \ {\rm  Dec}^n({\mathcal G},\chi_1,\chi_2)$.}

\bs

\noi {\bf Proof.}  Combining the $D_i$ together  gives
\begin{eqnarray*}
D_1+D_2- D_3&=&
-(ju_{20}+iu_{02}+(i+j)u_{11})'+(\lambda^2_i-\lambda^2_{i+j})u_{10}'+(\lambda^2_j-\lambda^2_{i+j})u_{01}'\\&& \ \ 
+(s_i+s_j-s_{i+j})u_{00}' -
(iy_1'+\lambda^2_i\chi_3')+\lambda^2_{i+j}\chi'.
\end{eqnarray*}
Using Lemma 3.5 (iii) we have
$\lambda^2_i-\lambda^2_{i+j}+\lambda^2_j=-ij$ which gives 
$\lambda^2_i-\lambda^2_{i+j}=-(\lambda^2_j+ij)$, and using
Lemma 3.5 (i)
we have
 $s_i+s_j-s_{i+j}=-(i
\lambda^2_j+\lambda^2_ij)$. 
With these relations 
we see
\begin{eqnarray*}
D_1+D_2- D_3&=&
-(ju_{20}'+iu_{02}'+(i+j)u_{11}') -(\lambda^2_j+ij)u_{10}'-(\lambda^2_i+ij)u_{01}'\\&& \ \ 
-(i\lambda^2_j+\lambda^2_ij)u_{00}'-
(iy_1'+\lambda^2_i\chi_3')+\lambda^2_{i+j}\chi'
\end{eqnarray*}
as required. 
According to Theorem 3.7 (vi) each of
 the classes $[D_1]$, $[D_2]$, and $[D_3]$ lie in
${\rm  Dec}^n({\mathcal G},\chi_1,\chi_2)$ (here $\theta$ is $i$, $j$ and $i+j$
respectively), giving the result. $\Box$

\bs


With this previous lemma we obtain a calculation of $A_0$ and $B_0$.

\bs

\noi {\bf Lemma 3.13.} {\em For $y_1$ and $\chi_3$ given by the
construction of $D_1$ and  $y_2$ and $\chi_4$ given by the symmetry and
$(*_v)$, we have
$A_0+(iy_1'+\lambda^2_i\chi_3')-\lambda^2_\psi\chi,
B_0+(jy_2'+\lambda^2_j\chi_4')-\lambda^2_\psi\chi\in Z^n(\G,\ZZ)$ with
$$
[A_0+(iy_1'+\lambda^2_i\chi_3')-\lambda^2_\psi\chi'] \in \  {\rm  Dec}^n({\mathcal G},\chi_1,\chi_2)\subseteq H^n(\G,\ZZ)$$
and
$$[B_0+(jy_2'+\lambda^2_j\chi_4')-\lambda^2_\psi\chi'] \in \
  {\rm  Dec}^n({\mathcal G},\chi_1,\chi_2)\subseteq H^n(\G,\ZZ).$$
 }

\noi {\bf Proof.}   
We consider the  two $(*_u)$ equations
\begin{eqnarray*}
\delta(u_{21})&=&-(ju_{20}'+iu'_{11}+iju'_{10}+\lambda^2_iu_{01}'
+\lambda^2_iju_{00}')+(\psi-\lambda^2_\psi)\chi'\\
\delta(u_{12})&=&A_0-(iu_{02}'+ju'_{11}+iju_{01}'+\lambda^2_ju_{10}'
+i\lambda^2_ju_{00}')+\lambda^2_{\psi}\chi.
\end{eqnarray*}
Adding these and subtracting using the Lemma 3.12
 gives
\begin{eqnarray*}
\delta(u_{12}+u_{21})
-(D_1+D_2-D_3)&=&A_0+(iy_1'+\lambda^2_i\chi_3')+(\psi-\lambda^2_{i+j})\chi'\\
&=&A_0+(iy_1'+\lambda^2_i\chi_3')-\lambda^2_\psi\chi'
\end{eqnarray*}
where the latter equation uses
$\psi-\lambda_{i+j}=\psi-\lambda^2_{-\psi}=-\lambda^2_\psi$.
As $[D_1+D_2- D_3]
\in \ {\rm  Dec}^n({\mathcal G},\chi_1,\chi_2)$
this give the first assertion.
The analogous equation for $B_0$ follows by symmetry,
reversing $i$ and $j$.
With this the result follows. $\Box$

\bs



We note that using $\psi=\frac{1}{2}(i+j)=-(i+j)$ these calculations give that
\begin{eqnarray*}
\delta(A_0)&=&-[-i\delta(y_1')-ii'\chi_3''-\psi\psi''\chi'']\\
&=&i(\psi'\chi''-i'\chi_3'')+ii'\chi_3''+(i+j)\psi'\chi''
\\&=&(2i+j)\psi'\chi''
\end{eqnarray*}
 and similarly, interchanging $i$ and $j$ gives
$$\delta(B_0)=(i+2j)\psi'\chi''=-\delta(A_0).$$
That $\delta(A_0+B_0)=0$  must happen follows from the calculation of $\eta([c])$ 
in Theorem 3.11 (iii), but
it is good to have  computational confirmation.
Given these results 
we obtain
\begin{eqnarray*}
[A_0+B_0]&\equiv&-[(iy_1'+\lambda^2_i\chi_3')-\lambda^2_\psi\chi'+(jy_2'+\lambda^2_j\chi_4')-\lambda^2_\psi\chi']
\\&\equiv&-
[(iy_1'+\lambda^2_i\chi_3')+(jy_2'+\lambda^2_j\chi_4')+\lambda^2_\psi\chi'] \ \ ({\rm mod \  {\rm  Dec}^n({\mathcal G},\chi_1,\chi_2)}).
\end{eqnarray*}

Combining  the above gives the main result of this section.
 
 \bs

\noi {\bf Theorem 3.14.} {\em Suppose $\G$ is a Vovodsky group, $\chi_1,\chi_2\in
H^1(\G,\ZZ)$ are independent characters,
and  ${\HH}={\rm ker}(\chi_1)\cap{\rm ker}(\chi_2)$.  For $\M_i$ as in Section 1
we assume  $c=(u,v,0,-\chi T_4)$ is such that 
$[c]\in H^{n-1}(\G,\M_4\otimes_\Z\ZZ)$ so $\pi^*_4([c])=\chi\in H^{n-1}(\G,\ZZ)$. 
We recall that $\eta:H^{n-1}(\G,\M_4\otimes_\Z\ZZ)\rightarrow H^n(\G,\ZZ)$ 
denotes the Positselksi map with image ${\rm ker}\left(H^n(\G,\ZZ)\rightarrow
H^n(\HH,\ZZ)\right)$.  
As above we abbreviate
$\chi_1=i$, $\chi_2=j$ in cocycle formulae. Then,\\
(i) The image $\eta([c])$ $({\rm mod}  \ 
 {\rm Dec}^n({\mathcal G},\chi_1,\chi_2))$ can be
expressed
as $[z]$ where
$$z=iw_{1\,2}'+\lambda^2_iw_{0\,2}'+jw_{2\,1}'+
\lambda^2_jw_{2\,0}'+
ijw_{1\,1}'$$
for $w_{2\,0},w_{1\,1},w_{0\,2}\in H^{n-1}(\G,\ZZ)$, $w_{1\,2}, w_{2\,1}\in C^{n-1}(\G,\ZZ)$
with $w_{1\,1}=\chi$,
$\delta(w_{1\,2})=-iw_{0\,2}'-jw_{1\,1}$ and
$\delta(w_{2\,1})=-iw_{1\,1}'-jw_{2\,0}$.
Moreover, the conditions  on the $w_{r\,s}$ independent of the
previous hypotheses guarantee
that $z\in Z^n(\G,\ZZ)$. 
\\
(ii) The image  $\eta\left(H^{n-1}(\G,\M_4\otimes_\Z\ZZ)
\right)\subseteq H^n(\G,\ZZ)$ consists of
all  $[z]$ where
$$z=iw_{1\,2}'+\lambda^2_iw_{0\,2}'+jw_{2\,1}'+
\lambda^2_jw_{2\,0}'+
ijw_{1\,1}'\in H^n(\G,\ZZ)$$
with $w_{2\,0},w_{1\,1},w_{0\,2}\in Z^{n-1}(\G,\ZZ)$, 
$w_{2\,1},w_{1\,2}\in C^{n-1}(\G,\ZZ)$,
$\delta(w_{1\,2})=-iw_{0\,2}'-jw_{1\,1}$ and
$\delta(w_{2\,1})=-iw_{1\,1}'-jw_{2\,0}$.\\
(iii)  If $[\psi\chi']=[-(i+j)\chi']=0$ 
or $[(i-j)\chi']=0\in H^n(\G,\ZZ)$,
then 
$[z]\in
 {\rm Dec}^n({\mathcal G},\chi_1,\chi_2)$.
In particular if  $\chi=0$ then $[z]\in
 {\rm Dec}^n({\mathcal G},\chi_1,\chi_2)$.}

\bs

\noi {\bf Proof.}  For part (i) the preceding calculations  give the elements
$y_1, y_2\in C^{n-1}(\G,\ZZ)$
with
$\delta(y_1)=-i\chi_3'+\psi\chi'$ and 
$\delta(y_2)=-j\chi_4'+\psi\chi'$
where the image $\eta([c])$ $({\rm mod}  \ 
 {\rm Dec}^n({\mathcal G},\chi_1,\chi_2))$ can be
expressed
as $[z]$ where
$$z=(iy_1'+\lambda_i\chi_3')+(jy_2'+\lambda_j\chi_4')
+\lambda^2_\psi\chi'.$$
with $\chi_3, \chi_4\in Z^{n-1}(\G,\ZZ)$,
  $y_1, y_2\in C^{n-1}(\G,\ZZ)$.
 Recalling that $\psi=-(i+j)$,
we note that as $\lambda^2_a+\lambda^2_b-\lambda^2_{a+b}=-ab$ and
as $\lambda^2_{-a}=\lambda^2_a+a$ we have 
$$\lambda^2_\psi=\lambda^2_{-i-j}=\lambda^2_{-i}+\lambda^2_{-j}+ij=
i+j+ij+\lambda^2_i+\lambda^2_j
$$
so we can express
$$z=i(y_1+\chi)'+\lambda^2_i(\chi_3+\chi)'+j(y_2
+\chi)'+\lambda^2_j(\chi_4+\chi)'+
ij\chi'.$$
This give the desired format for $z$ and the conditions on the $w_{r\,s}$ follow 
readily from the conditions on $\chi_3, \chi_4\in Z^{n-1}(\G,\ZZ)$ and
  $y_1, y_2\in C^{n-1}(\G,\ZZ)$.  The final statement in (i) is a direct computation.

For part (ii) we apply part (i) noting that any element of 
${\rm Dec}^n({\mathcal G},\chi_1,\chi_2)$ is of the form 
$[i\mu_1'+j\mu'_2]\in H^n(\G,\ZZ)$ where
$\mu_1,\mu_2\in Z^{n-1}(\G,\ZZ)$.  So we can replace
$u_{1\,2}$ with $u_{1\,2}+\mu_1$ and replace
$u_{2\,1}$ with $u_{2\,1}+\mu_2$ in the formulation of the second version
of part (i) and  the ${\rm Dec}^n({\mathcal G},\chi_1,\chi_2)$ term
goes away.

For (iii), returning to expression for $z$ using
 $y_1,y_2,\chi_3,\chi_4$ as above we let $\mu$ be such that
$\delta(\mu)=\psi\chi'\in Z^n(\G,\ZZ)$. Then $\delta(i\mu')=-i\delta(\mu')=
 -i\psi'\chi''$ and $\delta(j\mu')=-j\delta(\mu')=
 -j\psi'\chi''$. But then $\delta(i(y_1-\mu)'+\lambda^2_i\chi'_3)=0$
 and $\delta(j(y_2-\mu)'+\lambda^2_j\chi'_4)=0$.  From this
 \begin{eqnarray*}
 [z]&=&[
iy_1'+\lambda^2_i\chi_3'+jy_2'
+\lambda^2_j\chi_4'+
\lambda^2_{\psi}\chi']\\
&=&[i(y_1-\mu)'+\lambda^2_i\chi'_3]+[(j(y_2-\mu)'+\lambda^2_j\chi'_4]\\&& \ \ 
+[(i+j)\mu'+\lambda^2_{\psi}\chi']
\in H^2(\G,\ZZ).
\end{eqnarray*}
The hypothesis on $\G$ that all $H^n(\G,\ZZ)$ classes vanishing in a cyclic $3$-extension
are in ${\rm Dec}^n({\mathcal G},\chi_1,\chi_2)$ for that extension gives
 the first case, for  the latter three summands
lie in three (different) $3$-extensions and therefore   $[z]\in 
{\rm Dec}^n({\mathcal G},\chi_1,\chi_2).$

The second case proceeds similarly. We first note that
$\lambda^2_{-(i+j)}+\lambda^2_{-(i-j)}=\lambda^2_i-(i+j)(i-j)
=\lambda^2_i-i^2+j^2=\lambda^2_i-(i-\lambda^2_i)+(j-\lambda^2_j)=
-i+j-\lambda^2_i-\lambda^2_j=-i-\lambda^2_i-j-\lambda^2_{-j}$
using $\lambda^2_j=\lambda^2_{-j}-j$. From this we find
\begin{eqnarray*}
 [z]&=&[
iy_1'+\lambda^2_i\chi_3'+jy_2'
+\lambda^2_j\chi_4'+
\lambda^2_{\psi}\chi']\\
&=&[
i(y_1-\chi)'+\lambda^2_i(\chi_3-\chi)'\\&& \ -j(-y_2'+\chi_4+\chi')
+\lambda^2_{-j}(\chi_4'-\chi)-
\lambda^2_{-(i-j)}\chi']
\end{eqnarray*}
If $\delta(\nu)=(i-j)\chi'$ then
$\delta((i-j)\nu')=-(i-j)(i-j)'\chi''$ and therefore
$(i-j)\nu'-\lambda^2_{-(i-j)}\chi'\in Z^2(\G,\ZZ)$.
Also, $\delta(y-\nu-\chi)=-(i+j)\chi'-i\chi'_3-(i-j)\chi'=-i\chi'_3+i\chi'$
so we find $\delta(i(y_1-\nu-\chi)+\lambda^2_i(\chi_3-\chi)')=
-i(-i'\chi''_3+i'\chi'')+-ii'(\chi_3-\chi)''=0$.  From this we find that
\begin{eqnarray*}
 [z]&=&[ i(y_1-\nu-\chi)+\lambda^2_i(\chi_3-\chi)']\\&&+
[-j(-y_2-\nu+\chi_4+\chi)'+\lambda^2_{-j}(\chi_4-\chi)']+
[(i-j)\nu'-\lambda^2_{-(i-j)}\chi']
\end{eqnarray*}
which as discussed in the previous case lies in 
${\rm Dec}^n({\mathcal G},\chi_1,\chi_2)$  giving (iii).  $\Box$

\bs

\noi {\bf Remarks 3.15.} (i) The conditions  
that $\delta(w_{1\,2})=-iw_{0\,2}'-jw_{1\,1}$ and
$\delta(w_{2\,1})=-iw_{1\,1}'-jw_{2\,0}$ for 
 $w_{2\,0},w_{1\,1},w_{0\,2}\in Z^{n-1}(\G,\ZZ)$
in the expression
$$z=iw_{1\,2}'+\lambda_iw_{0\,2}'+jw_{2\,1}'+
\lambda_jw_{2\,0}'+
ijw_{1\,1}'\in Z^n(\G,\ZZ)$$ are
readily checked to be equivalent to the assertion
that $\delta(u)=zT_G\in C^{n}(\G\ZZ[G])$ where
$u=w_{1\,2}t_1t_2^2+w_{0\,2}t_2^2+w_{2\,1}t_1^2t_2+
w_{2\,0}t_1^2+w_{1\,1}t_1t_2\in C^{n-1}(\G\ZZ[G])$.
This is to be expected, for the latter expression is, in view of
Shapiro's Lemma as discussed in Lemma 2.1, equivalent to 
$[z]\mapsto0\in H^n(\G,{\rm CoInd}^G_{\{1\}}(\ZZ))\cong H^n({\mathcal H},\ZZ)$.  It is for this
reason that the notation using $w_{r\,s}$ is used in the theorem.\\
(ii) Here we quickly check the formulae for $\delta$ on $\Z/9\Z[\langle t\rangle]$
 in Theorem 3.7 (iv) above by applying $\delta$ twice.
By that result  we have for $\tilde{u}=\tilde{u}_0+\tilde{u}_1t+\tilde{u}_2t^2$,
\begin{eqnarray*}
\delta(\tilde{u})&=&
\delta(\tilde{u}_0)+\left(\delta(\tilde{u}_1)+\widehat{\theta}\tilde{u}_0
-3\left( \theta u_2+\lambda^2_\theta u_1\right)\right)t
\\&&+\left(\delta(\tilde{u}_2)+\widehat{\theta}\tilde{u}_1+\lambda^2_{\widehat{\theta}}\tilde{u}_0
-3\left( (\theta+\lambda^2_\theta)u_2+\lambda^2_\theta u_1\right)
\right)t^2
\end{eqnarray*}
To check  $\delta\circ\delta(\tilde{u})=0$ we compute
the coefficients  $\Delta_1$ and $\Delta_2$ where $\delta\circ\delta(\tilde{u})=\delta(\delta(\tilde{u}_0))+
\Delta_1t+\Delta_2t^2$. Using $3\cdot 3=0$ 
and $\delta(\widehat{\theta}\,)\tilde{u}_0=3\left(\theta\lambda^2_{\theta'} 
+\lambda^2_\theta\theta'\right) u_0''$ we find
\begin{eqnarray*}
\Delta_1&=&\delta\left(\delta(\tilde{u}_1)+\widehat{\theta}\tilde{u}_0
-3\left( \theta u_2+\lambda^2_\theta u_1\right)\right)+\widehat{\theta}\delta(\tilde{u}_0)\\
&& \ -3\left(\theta(\delta(u_2)+\theta u_1'+\lambda^2_\theta u_0')
+\lambda^2_\theta(\delta(u_1)+\theta u_0')
\right)\\
&=&\delta(\widehat{\theta})\tilde{u}_0-\widehat{\theta}\delta(\tilde{u}_0)
+3\left(\theta\delta(u_2)+\theta\theta'u_1+\lambda^2_\theta\delta(u_1)\right)
+\widehat{\theta}\delta(\tilde{u}_0)\\
&& -3\left(\theta\delta(u_2)+\theta\theta' u_1''+\theta\lambda^2_{\theta'} u_0'
+\lambda^2_\theta\delta(u_1)+\lambda^2_\theta\theta' u_0'')
\right)\\
&=&\delta(\widehat{\theta})\tilde{u}_0-3\left(\theta\lambda^2_{\theta'} u_0'
+\lambda^2_\theta\theta' u_0''\right)\\&=&0 
\end{eqnarray*}
and
\begin{eqnarray*}
\Delta_2&=&\delta\left(\delta(\tilde{u}_2)+\widehat{\theta}\tilde{u}_1+
\lambda^2_{\widehat{\theta}}\tilde{u}_0
-3\left( (\theta+\lambda^2_\theta) u_2+\lambda^2_\theta u_1\right)\right)\\
&& \ +\widehat{\theta}\left(\delta(\tilde{u}_1)+\widehat{\theta}\tilde{u}_0
-3\left( \theta u_2+\lambda^2_\theta u_1\right)\right)+\lambda^2_{\widehat{\theta}}\delta(\tilde{u}_0)\\
&& \ -3\left((\theta+\lambda^2_\theta)(\delta(u_2)+\theta u_1'+\lambda^2_\theta u_0')
+\lambda^2_\theta(\delta(u_1)+\theta u_0')
\right)\\
&=&\delta(\widehat{\theta})\tilde{u}_1-\widehat{\theta}\delta(\tilde{u}_1)+
\widehat{\theta}\delta(\tilde{u}_1)
+\delta(\lambda^2_{\widehat{\theta}})\tilde{u}_0-
\lambda^2_{\widehat{\theta}}\delta(\tilde{u}_0)+\widehat{\theta}\widehat{\theta}'
\tilde{u}_0+\lambda^2_{\widehat{\theta}}\delta(\tilde{u}_0)
\\&&-3\left(-\theta\theta' u_2''-(\theta+\lambda^2_\theta)\delta(u_2)
-\theta\theta' u_1''-\lambda^2_\theta \delta(u_1')\right)\\
&& \ -3\left( \theta\theta' u_2''+\theta\lambda^2_{\theta'} u_1''\right)\\
&& -3\left((\theta+\lambda^2_\theta)(\delta(u_2)+\theta u_1'+\lambda^2_\theta u_0')
+\lambda^2_\theta(\delta(u_1)+\theta u_0')
\right)\\
&=&3(\theta\lambda^2_{\theta'}+\lambda^2_\theta\theta')u_1''+
3(\lambda^2_\theta\lambda^2_{\theta'}+\theta\lambda^2_{\theta'}+\lambda^2_\theta\theta')u_0''\\
&& \ -3\left((\theta\lambda^2_{\theta'}+\lambda^2_\theta\theta')u_1''+
(\theta\lambda^2_{\theta'}+\lambda^2_\theta\lambda^2_{\theta'}+\lambda^2_\theta\theta')u_0''
\right)\\
&=&0
\end{eqnarray*}
as expected.  Note that this also double checks the formulae for
$\delta(\widehat{\theta}\,)$ and $\delta(\lambda^2_{\widehat{\theta}})$.

\bs

\bs

\newpage

\noi {\bf \S 4.  Calculating the Bochstein.}

\bs

In this section  the Bochstein 
${\mathcal B}_4:H^n(\G,\M_4\otimes_\Z\ZZ)\rightarrow
H^{n+1}(\G,\M_4\otimes_\Z\ZZ)$ is computed. For $i=1,2,3$ the value ${\mathcal B}_i$ is
zero because $\G$ assumed to be is a  Voevodsky group (in the applications they
are pro-$p$ Galois groups of fields.) The computation here shows
that when $i=4$, where $\M_4$ is the cokernel of
$d_2:\M_2\rightarrow\M_3$,
the value of ${\mathcal B}_4$ is nonzero. However, what we actually compute is
$h_3\circ{\mathcal B}_4$, which is zero, and this (by careful inspection
of his proof) suffices to apply the
Positselski results we need.  For the remainder of this section
we will abbreviate ${\mathcal B}_4$ by ${\mathcal B}$ since it is the only
case that concerns us.

We define $t_1:=\tau_1-1$, $t_2:=\tau_2-1$,
$t_3:=\tau_1\tau_2-1$, $T_1=1+\tau_1+\tau_1^2$,
$T_2=1+\tau_2+\tau_1^2$,
$T_3=1+\tau_1\tau_2+\tau_1^2\tau_2^2$,
and $T_4=1+\tau_1\tau_2^2+\tau_1^2\tau_2$ in $\Z/9\Z[G]$ exactly as they 
were used in $\ZZ[G]$. They will not be distinguished by using 
a different notation.  With this, we emphasize that now 
$t_1^2\neq T_1\in \Z/9\Z[G]$. This next result collects some basic facts we need
in $\Z/9\Z[G]$ for computing $\delta:
C^1(\G,\Z/9\Z[G])\rightarrow C^2(\G,\Z/9\Z[G])$.

\bs


\noi {\bf Facts 4.1.} {\em In $\Z/9\Z[G]$ we have\\
(i) $t_1T_1=t_2T_2=t_3T_3=0$,\\
(ii) $t_1^3=-3(t_1+t_1^2)$, $t_2^3=-3(t_2+t_2^2)$,\\
(iii) $\tau_2T_3=\tau_1^2T_3$ so $t_2T_3=(2t_1+t_1^2)T_3$
and $t_2^2T_3=(t_1^2-3(t_1+t_1^2))T_3$ \\
(iv)
$\tau_2T_4=\tau_1T_4$ so $t_1T_4=t_2T_4$\\
Suppose  $\gamma=\tau_1^i\tau_2^j$. Then\\
(v)
$(\gamma-1)T_1=(\widehat{j}t_2+\lambda^2_{\widehat{j}}t_2^2)T_1$\\
(vi) $(\gamma-1)T_2=(\widehat{i}t_1+\lambda^2_{\widehat{i}}t_1^2)T_2$,\\
(vii) $(\gamma-1)T_3=(\widehat{(i-j)}t_1+\lambda^2_{\widehat{i-j}}t_1^2)T_3$, \ {\rm and} \\
(viii) $(\gamma-1)T_4=(\widehat{(i+j)}t_1+\lambda^2_{\widehat{i+j}}t_1^2)T_4$.
}

\bs

\noi {\bf Proof.} (i) to (iv) are straight-forward and for (v) through (viii) we note that $\gamma=\tau_1^i\tau_2^j=\tau_1^{i-j}(\tau_1\tau_2)^j=
\tau_1^{i+j}(\tau_1\tau_2^2)^{-j}$ and the results then follow from 
(i) through (iv). $\Box$

\bs

We remark that in  the previous section we used the  equation
$(\gamma-1)T_4=(\psi t_3+\lambda^2_{\psi}t_3^2)T_4
\in C^1(\G,\ZZ[G])$ where $\psi=\frac{1}{2}(i+j)=-(i+j)$. But it is different here
because we are working with multiples of $t_1T_4$ instead of $t_3T_4$
as well as working (mod $9$) which explains why we need the ``hat'' functions.
As an application  we have the following.

\bs

\noi 
{\bf Theorem 4.2.} {\em Suppose $\tilde{\chi} \in Z^{n-1}(\G,\Z/9\Z)$.
Then\\
(i)
$\delta(\tilde{\chi}T_4)=
\left(\widehat{(i+j)}t_1+\lambda^2_{\widehat{i+j}}t_1^2\right)\tilde{\chi}'T_4.
$\\
(ii)
$\delta(\tilde{\chi}T_4)t_1=\left(
\widehat{(i+j)}t^2_1-3(\lambda^2_{{i+j}}(t_1+t_1^2)\right){\chi}'T_4.
$\\
(iii) $\delta(st_1^2T_4)\in3\delta(C^{n-1}(\G,\ZZ[G]))$
whenever $s\in Z^{n-1}(G,\Z/9\Z)$. The same calculation gives
$\delta(st_1^2T_2),\delta(st_2^2T_1),\delta(st_1^2T_3)\in3\delta(C^{n-1}(\G,\ZZ[G]))$.
}

\bs

\noi {\bf Proof.}
Part (i) is clear by the Fact 4.1 (viii), for as $\tilde{\chi}\in Z^{n-1}(\G,\Z/9\Z)$
we have $\delta(\tilde{\chi})=0$ and therefore
we have
$$\delta(\tilde{\chi}T_4)=(\gamma-1)\tilde{\chi}'T_4+\delta(\tilde{\chi})'T_4=
\left(\widehat{(i+j)}t_1+\lambda^2_{\widehat{i+j}}t_1^2\right)\tilde{\chi}'T_4.
$$
Part (ii) follows from (i) multiplying by $t_1$ and using Fact 4.1 (ii).
For (iii) we note that 
\begin{eqnarray*}
\delta(t_1^2T_4)&=&(\gamma-1)t_1^2T_4=
\left(\widehat{(i+j)}t^3_1+\lambda^2_{\widehat{i+j}}t_1^4\right)T_4\\
&=&-3\left((i+j)t_1+(i+j+\lambda^2_{i+j})t_1^2\right)T_4
\\
&=&-3\delta\left((1+t_1)T_4\right).
\end{eqnarray*}
From this, whenever $s\in Z^{n-1}(G,\Z/9\Z)$ we see   $\delta(st_1^2T_4)
=-3\delta\left(s(1+t_1)T_4\right)$ lies in $3\delta(C^{n-1}(\G,\ZZ[G]))$. Similarly,
$\delta(st_1^2T_2),\delta(st_2^2T_1),\delta(st_1^2T_3)
$ $\in3\delta(C^{n-1}(\G,\ZZ[G]))$.
$\Box$

\bs

We now turn to explicitly formalizing the lifts given by the hypothesis that $\G$ is a Voevodsky group.
We start with 
the equations that are necessary for the construction
of the algebras  $D_1$, $D_2$, $D_3$, and $D_4$ which we will refer to
as the nine $(**_u)$ equations since they are derived from the $(*_u)$
equations:
\begin{eqnarray*}
\delta(u_{00})&=&0\\
\delta(u_{10}-\chi_3)+iu'_{00}&=&0\\
\delta(u_{20}-y_1)+i(u'_{10}-\chi_3)+\lambda^2_iu_{00}&=&0 \\
\delta(u_{01})+ju'_{00}&=&0\\
  \delta(u_{02})+ju'_{01}+\lambda^2_ju_{00}&=&0\\
  \delta(u_{10}+u_{01}-\chi)+(i+j)u_{00}&=&0\\
  \delta(u_{20}+u_{02}+u_{11})+(i+j)(u_{10}+u_{01}-\chi)'+\lambda^2_{i+j}u_{00}
  &=&0\\
  \delta(u_{10}-u_{01})+(i-j)u_{00}&=&0\\
 \delta(u_{01}+u_{20}+u_{02}-u_{11})+(i-j)(u_{10}-u_{01})' +(\lambda^2_{i-j})u_{00}'&=&0
\end{eqnarray*}
where $y_1,\chi_3$ satisfy $ \delta(y_1)=\psi\chi'-i\chi_3'$.  Here
the first seven $(**_u)$ equations are rearrangements of the corresponding
first seven $(*_u)$ equations,  while
the latter two $(**_u)$ equations are not.  Instead, they are they are readily checked to be 
linear combinations of the first seven $(*_u)$ equations
using
$\lambda^2_{-j}=j+\lambda^2_j$ which in turn
gives $\lambda^2_{i-j}=-\lambda^2_i+\lambda^2_j+j-ij$.

We use the lifting criteria described in Theorem 3.7 (iv) three times;
once for each of $\cdot T_1$, $\cdot T_3$ and $\cdot T_4$.  These are
given next.

\bs

\noi {\bf 4.3.1. The case of $\cdot T_1$}. \ 
The system given by the first, fourth and fifth $(**_u)$ equations is
\begin{eqnarray*}
\delta(u_{00})&=&0\\
\delta(u_{01})+ju'_{00}&=&0\\
  \delta(u_{02})+ju'_{01}+\lambda^2_ju_{00}'&=&0\, ,
\end{eqnarray*}
which by Theorem 3.7 (ii) is equivalent to
$$(u_{00}+u_{01}t_2+u_{02}t_2^2)T_1\in Z^{n-1}(\G,\ZZ[G]\cdot T_1).$$
By Theorem 3.7 (iii) 
 we can lift to
$$(\tilde{u}_{00}+\tilde{u}_{01}t_2+\tilde{u}_{02}t_2^2)T_1
\in Z^{n-1}(\G,\Z/9\Z[G]\cdot T_1)$$
where $\tilde{u}_{00}\equiv u_{00}$, $\tilde{u}_{01}\equiv u_{01}$,
and $\tilde{u}_{02}\equiv u_{02}$ (mod $3$).
By   Theorem 3.7 (iv) this lifted condition is equivalent to
\begin{eqnarray*}
\delta(\tilde{u}_{00})&=&0\\
\delta(\tilde{u}_{01})+\widehat{j}\tilde{u}'_{00}-3(j{u}_{02}'+
\lambda^2_{{j}}{u}'_{01})&=&0\\
  \delta(\tilde{u}_{02})+\widehat{j}\tilde{u}'_{01}+\lambda^2_
  {\widehat{j}}\tilde{u}_{00}'-3(({j}+\lambda^2_{j}){u}_{02}'+
\lambda^2_{{j}}{u}'_{01})&=&0.
\end{eqnarray*}

\bs


\noi {\bf 4.3.3. The case of $\cdot T_3$}. \ 
Next, using the first, eighth and ninth $(**_u)$ equations, 
where if we denote $\mu_1=u_{10}+u_{01}$ and
$\mu_2=u_{01}+u_{20}+u_{02}-u_{11}$ we find
\begin{eqnarray*}
\delta(u_{00})&=&0\\
  \delta(\mu_1)+(i-j)u_{00}&=&0\\
 \delta(\mu_2)+(i-j)\mu_1' +\lambda^2_{i-j}u_{00}'&=&0 \, .
\end{eqnarray*}
Noting that  $t_1T_4=t_2^2T_4$ gives $(\gamma-1)T_4=(i-j)t_1T_4+\lambda^2_{i-j}t_1^2T_4$, this system is equivalent to
$$\left(u_{00}+\mu_1t_1+\mu_2t_1^2\right)T_3\in Z^{n-1}(\G,\ZZ[G]\cdot T_3).$$
By Theorem 3.7 (iii) we can lift to
$$\left(\tilde{u}_{003}+\tilde{\mu}_{1}t_1+\tilde{\mu}_{2}t_1^2
\right)T_3\in Z^{n-1}(\G,\Z/9\Z[G]\cdot T_3)$$
where
${\tilde{u}_{003}}\equiv u_{00}$, ${\tilde{\mu}_{1}}\equiv {\mu}_1$, 
and ${\tilde{\mu}_{2}}\equiv {\mu}_{2}$ (mod $3$). 
The lifted equation  is equivalent to
\begin{eqnarray*}
\delta(\tilde{u}_{003})&=&0\\
\delta(\tilde{\mu}_1)+(\widehat{i-j})\tilde{u}_{003}'-3\left(({i-j})\tilde{\mu}_2'+
\lambda^2_{{i-j}}\tilde{\mu}_1'
\right)&=&0\\
  \delta(\tilde{\mu}_{2})+(\widehat{i-j})\tilde{\mu}_1'+\lambda^2_
  {\widehat{i-j}}\tilde{u}_{003}'
  -3\left((({i-j})+\lambda^2_{{i-j}})\tilde{\mu}_2'+
\lambda^2_{{i+j}}\tilde{\mu}_1'
\right)&=&0 \ .
\end{eqnarray*}
We can express $\tilde{u}_{003}=\tilde{u}_{00}+3u_{003}^*$ with 
$u_{003}^*\in Z^{n-1}(\G,\ZZ)$ and then the latter system can be written
as 
\begin{eqnarray*}
\delta(\tilde{u}_{00})&=&0\\
\delta(\tilde{\mu}_1-3u_{003}^*)+(\widehat{i-j})\tilde{u}_{00}' 
-3\left(({i-j})(\tilde{\mu}_2'-u_{003}^*)'+
\lambda^2_{{i-j}}(\tilde{\mu}_1-3u_{003}^*)'
\right)&=&0\\
  \delta(\tilde{\mu}_{2}-u_{003}^*)+(\widehat{i-j})
  (\tilde{\mu}_1-3u_{003}^*)'+\lambda^2_
  {\widehat{i-j}}\tilde{u}_{00}'
  \hspace{1.75in}\\-3\left((({i-j})+\lambda^2_{{i-j}})(\tilde{\mu}_2-u_{003}^*)'+
\lambda^2_{{i-j}}(\tilde{\mu}_1-3u_{004}^*)'
\right)&=&0 \ .
\end{eqnarray*}
where these equations are equivalent to
$$\left(\tilde{u}_{00}+(\tilde{\mu}_1-3u_{003}^*)t_1+
(\tilde{\mu}_2-u_{003}^*)'t_1^2\right)T_4\in Z^{n-1}(\G,\Z/9\Z[G]\cdot T_3).$$

\bs


\noi {\bf 4.3.4. The case of $\cdot T_4$}. \  
Third,  using the first, sixth and seventh $(**_u)$ equations, 
where if we denote $\sigma_1=u_{10}+u_{01}$ and
$\sigma_2=u_{20}+u_{02}+u_{11}$ we find
\begin{eqnarray*}
\delta(u_{00})&=&0\\
  \delta(\sigma_1-\chi)+(i+j)u_{00}&=&0\\
  \delta(\sigma_2)+(i+j)(\sigma_1-\chi)'+\lambda^2_{i+j}u_{00}
  &=&0 \ .
\end{eqnarray*}
Noting that  $t_1T_4=t_2T_4$ gives $(\gamma-1)T_4=(i+j)t_1T_4+\lambda^2_{i+j}t_1^2T_4$, this system is equivalent to
$$\left(u_{00}+(\sigma_1-\chi)t_1+\sigma_2t_1^2\right)T_4\in Z^{n-1}(\G,\ZZ[G]
\cdot T_4).$$
By Theorem 3.7 (iii) we can lift to
$$\left(\tilde{u}_{004}+\tilde{\sigma}_{1}t_1+\tilde{\sigma}_{2}t_1^2
\right)T_4\in Z^{n-1}(\G,\Z/9\Z[G]\cdot T_4)$$
where
${\tilde{u}_{004}}\equiv u_{00}$, ${\tilde{\sigma}_{1}}\equiv {\sigma}_{1}-\chi$, 
and ${\tilde{\sigma}_{2}}\equiv {\sigma}_{2}$ (mod $3$). 
The lifted equation  is equivalent to
\begin{eqnarray*}
\delta(\tilde{u}_{004})&=&0\\
\delta(\tilde{\sigma}_1)+(\widehat{i+j})\tilde{u}_{004}'-3\left(({i+j})\tilde{\sigma}_2'+
\lambda^2_{{i+j}}\tilde{\sigma}_1'
\right)&=&0\\
  \delta(\tilde{\sigma}_{2})+(\widehat{i+j})\tilde{\sigma}_1'+\lambda^2_
  {\widehat{i+j}}\tilde{u}_{004}'
  -3\left((({i+j})+\lambda^2_{{i+j}})\tilde{\sigma}_2'+
\lambda^2_{{i+j}}\tilde{\sigma}_1'
\right)&=&0 \ .
\end{eqnarray*}
We can express $\tilde{u}_{004}=\tilde{u}_{00}+3u_{004}^*$ with 
$u_{004}^*\in Z^{n-1}(\G,\ZZ)$ and then the latter system can be written
as 
\begin{eqnarray*}
\delta(\tilde{u}_{00})&=&0\\
\delta(\tilde{\sigma}_1-3u_{004}^*)+(\widehat{i+j})\tilde{u}_{00}' 
-3\left(({i+j})(\tilde{\sigma}_2'-u_{004}^*)'+
\lambda^2_{{i+j}}(\tilde{\sigma}_1-3u_{004}^*)'
\right)&=&0\\
  \delta(\tilde{\sigma}_{2}-u_{004}^*)+(\widehat{i+j})
  (\tilde{\sigma}_1-3u_{004}^*)'+\lambda^2_
  {\widehat{i+j}}\tilde{u}_{00}'
  \hspace{1.75in}\\-3\left((({i+j})+\lambda^2_{{i+j}})(\tilde{\sigma}_2-u_{004}^*)'+
\lambda^2_{{i+j}}(\tilde{\sigma}_1-3u_{004}^*)'
\right)&=&0 \ .
\end{eqnarray*}
where these equations are equivalent to
$$\left(\tilde{u}_{00}+(\tilde{\sigma}_1-3u_{004}^*)t_1+
(\tilde{\sigma}_2-u_{004}^*)'t_1^2\right)T_4\in Z^{n-1}(\G,\Z/9\Z[G]\cdot T_4).$$


\bs

Given these three lifts, we may combine them to give the lift of $u$.

\bs

\noi {\bf Definition 4.4.} For $u\in C^{n-1}(\G,\ZZ[G])$ satisfying the equations $(**_u)$ 
above we define
\begin{eqnarray*}
\tilde{u}&:=&\tilde{u}_{00}+(\tilde{\mu}_1+\tilde{u}_{01})t_1
+\tilde{u}_{01}t_2+\frac{1}{2}(\tilde{\sigma}_2+\tilde{\mu}_2-\tilde{u}_{01}-2\tilde{u}_{02})t_1^2\\
&& \ +\frac{1}{2}(\tilde{\sigma}_2-\tilde{\mu}_2+\tilde{u}_{01})t_1t_2
+\tilde{u}_{02}t_2^2+\tilde{u}_{21}t_1^2t_2+\tilde{u}_{12}t_1t_2^2+\tilde{u}_{22}t_1^2t_2^2.
\end{eqnarray*}
where $\tilde{u}_{01}$, $\tilde{u}_{02}$, $\tilde{\mu}_{1}$, $\tilde{\mu}_{2}$, $\tilde{\sigma}_{1}$
and $\tilde{\sigma}_{2}$ are given by the  above liftings 
in 4.3.1, 4.3.3 and 4.3.4, and where $\tilde{u}_{21}$,
$\tilde{u}_{12}$ and $\tilde{u}_{22}$ are arbitrary lifts
of ${u}_{21}$,
${u}_{12}$ and ${u}_{22}$ respectively.  

For $v\in C^{n-1}(\G,\ZZ[G])$ satisfying the  analogous equations $(**_v)$ we define $\tilde{v}$ according to exactly the same plan; here the
roles of $t_1$ and $t_2$ are interchanged as are $i$ and $j$. 


\bs

The next result tabulates the $\delta(\tilde{u}\cdot T_i)$ and
$\delta(\tilde{v}\cdot T_i)$ for $i=1,2,3,4$. We remark that $t_1$ and
$t_2$ are interchanged when we go from $\tilde{u}$ to $\tilde{v}$ so the
same is true for $T_1$ and $T_2$ when we apply the results of
4.3.  But we also note that
$T_3$ and $T_4$ remain stable under this interchange.

\bs

\noi {\bf Theorem 4.5.}
{\em For $\tilde{u}, \tilde{v}\in C^{n-1}(\G,\Z/9\Z[G])$ as in Definition 4.4 we have 
\begin{eqnarray*}
\delta(\tilde{u}T_1) \ = \ \delta(\tilde{v}T_2)&=&0\\
\delta(\tilde{u}T_3) \ \equiv \ \delta(\tilde{v}T_3)&\equiv &0
\\
\delta(\tilde{u}T_4) \ \equiv \ \delta(\tilde{v}T_4)&\equiv&
\widehat{(i+j)}\tilde{\chi}'t^2_1T_4-3\lambda^2_
{i+j}\chi'(t_1+t_1^2)T_4\\
\delta(\tilde{u}(t_1+t_1^2)T_4) &\equiv&0\\
\delta(\tilde{v}(t_2+t_2^2)T_4) &\equiv&0\\
\delta(\tilde{\chi}T_4)t_1 \ \equiv \ \delta(\tilde{\chi}T_4)t_2&\equiv&
\widehat{(i+j)}\tilde{\chi}'t^2_1T_4-3\lambda^2_{{i+j}}{\chi}'(t_1+t_1^2)T_4
\ \ 
\end{eqnarray*}
where all congruences mean $\left({\rm mod} \ 3\delta(C^{n-1}(\G,\Z/9\Z[G]))\right)$.}

\bs

\noi {\bf Proof.} For
 $\delta(\tilde{u}T_1)$, using the definition of $\tilde{u}$ and the equations in
 4.3.1  we find
$$
\delta(\tilde{u}T_1)=\delta\left((\tilde{u}_{00}+\tilde{u}_{01}t_2
+\tilde{u}_{02}t_2^2)T_1\right)=0 \ . $$
For $\delta(\tilde{u}T_3)$, by the equations in
 4.3.3 we find
$$\left(\tilde{u}_{00}+(\tilde{\mu}_1-3u_{003}^*)t_1+
(\tilde{\mu}_2-u_{003}^*)'t_1^2\right)T_3\in Z^{n-1}(\G,\Z/9\Z[G]\cdot T_3),$$
and using the definition of  $\tilde{u}$ along with  $t_2T_3=(2t_1+t_1^2)T_3$,  $t_2^2T_3=(t_1^2-3(t_1+t_1^2))T_3$ and $t^3\in3\Z/9\Z[G]$ from Facts 4.1, we have
for some $E_1\in C^{n-1}(\G,\Z/9\Z[G])$,
\begin{eqnarray*}
\delta(\tilde{u}T_3)&=&\delta\left((\tilde{u}_{00}+
(\tilde{\mu}_{1}+3\tilde{u}_{01})t_1
+\tilde{\mu}_{2}t_1^2+3E_1)T_3\right)\\
&\equiv&\delta\left(\tilde{u}_{00}+(\tilde{\mu}_1-3u_{003}^*)t_1+
(\tilde{\mu}_2-u_{003}^*)'t_1^2\right)T_3\\
&& \ + \delta\left((
3\tilde{u}_{003}^*t_1
+\tilde{u}_{003}^*t_1^2)T_3\right)\\
&\equiv&\delta\left(\tilde{u}_{003}^*t_1^2T_3\right) \ \ \left({\rm mod} \ 
3\delta(C^{n-1}(\G,\Z/9\Z[G]))\right)
\end{eqnarray*}
However, by the Theorem 4.2 (iii),
$\delta\left(\tilde{u}_{003}^*t_1^2T_3\right) \in 
3\delta(C^{n-1}(\G,\Z/9\Z[G]))$
so we find that
$$\delta(\tilde{u}T_3)\in 
3\delta(C^{n-1}(\G,\Z/9\Z[G])).
$$

Next for $\delta(\tilde{u}T_4)$, by the equations in
 4.3.4 we find 
$$\left(\tilde{u}_{00}+(\tilde{\sigma}_1-3u_{004}^*)t_1+
(\tilde{\sigma}_2-u_{004}^*)'t_1^2\right)T_4\in Z^{n-1}(\G,\Z/9\Z[G]\cdot T_4),$$
so using the definition of  $\tilde{u}$ along with $t_1T_4=t_2T_4$ 
and $t^3\in3\Z/9\Z[G]$,
\begin{eqnarray*}
\delta(\tilde{u}T_4)&=&\delta\left((\tilde{u}_{00}+
(\tilde{\mu}_{1}+2\tilde{u}_{01})t_1
+\tilde{\sigma}_{2}t_1^2+(\tilde{u}_{21}+\tilde{u}_{12})t_1^3+\tilde{u}_{22}t_1^4)T_4\right)\\
&\equiv&\delta\left(\tilde{u}_{00}+(\tilde{\sigma}_1-3u_{004}^*)t_1+
(\tilde{\sigma}_2-u_{004}^*)'t_1^2\right)T_4\\
&& \ + \delta\left((
(\tilde{\mu}_{1}+2\tilde{u}_{01}-\tilde{\sigma}_{1}+3\tilde{u}_{004}^*)t_1
+\tilde{u}_{004}^*t_1^2\right)
\\&\equiv&\delta\left((
(\tilde{\mu}_{1}+2\tilde{u}_{01}-\tilde{\sigma}_{1})t_1+\tilde{u}_{004}^*t_1^2)T_4\right) 
\ \ \left({\rm mod} \ 
3\delta(C^{n-1}(\G,\Z/9\Z[G]))\right)
\end{eqnarray*}
However, 
as $\tilde{\mu}_{1}\equiv u_{10}-u_{01}$ and
$\tilde{\sigma}_1\equiv u_{10}+u_{01}-\chi$ (mod $3C^1(\G,\Z/9\Z)$)
we can express
$\tilde{\mu}_{1}+2\tilde{u}_{01}-\tilde{\sigma}_{1}=\tilde{\chi}+3E_2$
for some $E_2\in C^{n-1}(\G,\Z/9\Z[G])$.
As $\delta(\tilde{\chi})=0$ and $\delta(\tilde{u}_{004}^*)=0$  we can  continue
$\left({\rm mod} \ 
3\delta(C^{n-1}(\G,\Z/9\Z[G])\right)$,
\begin{eqnarray*}
\delta(\tilde{u}T_4)
&\equiv& \delta\left((\tilde{\chi}t_1+\tilde{u}_{004}^*t_1^2)T_4\right)\\
&=&((\widehat{i+j})t_1+\lambda^2_
{\widehat{i+j}}t_1^2)(\tilde{\chi}t_1+\tilde{u}_{004}^*t_1^2)'T_4
\\
&=&(\widehat{i+j})\tilde{\chi}'t_1^2T_4-3\lambda^2_{i+j}\chi'(t_1+t_1^2)T_4
\\&&  \ -3\left(({i+j})(t_1+t_1^2)+\lambda^2_
{{i+j}}t_1^2\right)\tilde{u}_{004}^*{}'T_4
 \ .
\end{eqnarray*}
 But $\delta((1+t_1)\tilde{u}_{004}^* T_4)=
((i+j) t_1+\lambda^2_{i+j}t_1^2+(i+j)t_1^2)\tilde{u}_{004}^*)T_4$ so we find that
$$
\delta(\tilde{u}T_4)\equiv
\widehat{(i+j)}\tilde{\chi}'t^2_1T_4-3\lambda^2_
{i+j}\chi'(t_1+t_1^2)T_4 \ \ \left({\rm mod} \ 3\delta(C^{n-1}(\G,\Z/9\Z))\right) \ .
$$

Next,  by the above $\delta(\tilde{u}T_4)
\equiv \delta((\tilde{\chi}t_1+\tilde{u}_{004}^*t_1^2)T_4)$, so
for some $E_3\in C^{n-1}(\G,\Z/9\Z[G])$,
$$
\delta(\tilde{u}(t_1+t_1^2)T_4)
\equiv \delta((\tilde{\chi}t^2_1+3E_3)T_4)\ \ 
\equiv \  \delta(\tilde{\chi}t^2_1T_4)
$$
where by Theorem 4.2 (iii) we have $\delta(\tilde{\chi}t^2_1T_4)\in 3\delta
(C^{n-1}(\G,\Z/9\Z[G]))$ and from this we see that
$\delta(\tilde{u}(t_1+t_1^2)T_4)
\equiv 0 \ \left({\rm mod} \  3\delta(C^{n-1}(\G,\Z/9\Z))\right)$.

Finally,
 we have by Theorem 4.2 (i), 
$$\delta(\tilde{\chi}T_4)=(\gamma-1)\tilde{\chi}'T_4=
(\widehat{(i+j)}t_1+\lambda^2_{\widehat{i+j}}t_1^2)\tilde{\chi}'T_4.
$$
The results for $\tilde{v}$ follow by symmetry
where $t_1$ and $t_2$ are interchanged
and this  concludes the proof. $\Box$

\bs

We now can give the main result of this section.

\bs

\noi {\bf Theorem 4.6.} {\em Suppose $[c]\in H^{n-1}(\G,\M_4\otimes_\Z\ZZ)$.
Then  the composite  
$h_3({\mathcal B}([c]))=0\in
H^n(\G,\M_3\otimes_\Z\ZZ)$.
}

\bs

\noi {\bf Proof.}
Following Theorem 3.11, labeling $c=(u,v,0,-\chi T_4)$, 
in order to to compute the composite $h_3({\mathcal B}([(u,v,0,-\chi T_4))])\in H^2(\G,\M_3\otimes_\Z\ZZ)$
 we use   $\tilde{u}, \tilde{v}\in  C^{n-1}(\G,\Z/9\Z[G])$
 and $\tilde{\chi}\in C^{n-1}(\G,\Z/9\Z)$ described above
to find the Bochstein ${\mathcal B}([(u,v,0,-\chi T_4)])$.
With this, $h_3({\mathcal B}([c]))=[h_3(\delta((\tilde{u},\tilde{v},0,-\tilde{\chi} T_4)))]$.
Then as
$h_2(\tilde{c}):=(\delta(\tilde{u})\kappa_1+\delta(\tilde{v})\kappa_2+\delta(\tilde{\chi}) T_4,0)
\in Z^n(\G,3\Z/9\Z[G])\cong Z^n(\G,\ZZ[G])$ we have
\begin{eqnarray*}h_3(\delta((\tilde{u},\tilde{v},0,-\tilde{\chi} T_4)))
&=&3\delta(\tilde{u},\tilde{v},0,-\tilde{\chi} T_4)-d_2(h_2(\delta(\tilde{u},\tilde{v},0,-\tilde{\chi} T_4))) \\
&\equiv&-d_2((\delta(\tilde{u})\kappa_1+\delta(\tilde{v})\kappa_2+\delta(\tilde{\chi}T_4),0)
 \\
&\equiv& \left(A(\tilde{u},\tilde{v},\tilde{\chi}),B(\tilde{u},\tilde{v},\tilde{\chi}),
0,0\right)  \ \ 
\end{eqnarray*}
$\left({\rm mod} \ 3\delta(C^{n-1}(\G,\M_3\otimes_\Z\Z/9\Z))\right)$
where  $A(\tilde{u},\tilde{v},\tilde{\chi})=
\delta(\tilde{u})\kappa_1t_1+\delta(\tilde{v})\kappa_2t_1+\delta(\tilde{\chi}T_4)t_1$
and $ B(\tilde{u},\tilde{v},\tilde{\chi})=\delta(\tilde{u})\kappa_1t_2+
\delta(\tilde{v})\kappa_2t_2+\delta(\tilde{\chi}T_4)t_2$.
By symmetry it will suffice to compute 
$\delta(\tilde{u})\kappa_1t_1$, $\delta(\tilde{u})\kappa_1t_2$
and  $\delta(\tilde{\chi}T_4)t_i$.  

At this point we recall that $\kappa_1\equiv(t_1+t_1^2)$ and $T_1\equiv t_1^2$
(mod $3(\Z/9\Z)[G]$) and therefore
$\kappa_1t_1\equiv t_1^2+t_1^3\equiv t_1^2\equiv T_1$
(mod $3(\Z/9\Z)[G]$). Also
\begin{eqnarray*}T_3&\equiv& t_1^2+2t_1t_2+t_2^2
+2t_1^2t_2+2t_1t_2^2+t_1^2t_2^2 \ \ \ {\rm (mod} \  3(\Z/9\Z)[G]) \ ,
\\ T_4&\equiv& t_1^2+t_1t_2+t_2^2+t_1^2t_2+t_1t_2^2
 \hspace{.8in}  {\rm (mod} \  3(\Z/9\Z)[G]) \ ,
\end{eqnarray*}
so using $t_2^2\equiv T_2\equiv-(T_1+T_3+T_4-T_G) \ {\rm (mod} \  3(\Z/9\Z)[G])$ 
we have
\begin{eqnarray*}
\kappa_1t_2&\equiv& t_1t_2+t_1^2t_2\\
&\equiv&T_3-T_4-t_1t_2^2-t_1^2t_2^2\\
&\equiv&T_3-T_4+(t_1+t_1^2)(T_3+T_4)\\ \ 
&\equiv&(1+t_1+t_1^2)T_3-(1-t_1-t_1^2)T_4 \ \ {\rm (mod} \  3(\Z/9\Z)[G]) \ .
\end{eqnarray*}
By symmetry (interchanging $t_1$ and $t_2$),
$$
\kappa_2t_1\equiv
(1+t_2+t_2^2)T_3-(1-t_2-t_2^2)T_4 \ \ {\rm (mod} \  3(\Z/9\Z)[G]) \ .
$$
From this, applying Theorem 4.5 we find
\begin{eqnarray*}
\delta(\tilde{u})\kappa_1t_1&=&\delta(\tilde{u})T_1 \ = \ 0 \ ,\\
\delta(\tilde{v})\kappa_2t_1&=&\delta(\tilde{v})\left(
(1+t_2+t_2^2)T_3-(1-t_2-t_2^2)T_4\right) \\
&=&-\left(
\widehat{(i+j)}\tilde{\chi}'t^2_1T_4-3\lambda^2_{{i+j}}{\chi}'(t_1+t_1^2)T_4
\right) \ ,\\
\delta(\tilde{\chi}T_4)t_1&\equiv&
\widehat{(i+j)}\tilde{\chi}'t^2_1T_4-3\lambda^2_{{i+j}}{\chi}'(t_1+t_1^2)T_4 \ .
\end{eqnarray*}
We find that
 $A(\tilde{u},\tilde{v},\tilde{\chi})=
\delta(\tilde{u})\kappa_1t_1+\delta(\tilde{v})\kappa_2t_1+\delta(\tilde{\chi}T_4)t_1=0$
and similarly $B(\tilde{u},\tilde{v},\tilde{\chi})=0$, giving the result. $\Box$

\bs

We conclude with the desired applications of Theorem 4.6.

\bs

\noi {\bf Theorem 4.7.} {\em Suppose that $G\cong \ZZ\oplus\ZZ$
and $\G$ is a pro-$3$ group
acting on $G$ via $\G/{\mathcal H}\cong G$ where
${\mathcal H}={\rm ker}(\chi_1)\cap {\rm ker}(\chi_2)$
and $\chi_1,\chi_2\in H^1(\G,\ZZ)$.
We assume the Voevodsky condition is satisfied
 for $\G$ and $\Z[G]$ $({\rm mod} \  9)$.
 Suppose that $\M_1$, $\M_2$, $\M_3$ and $\M_4$
 are defined as in Section 1 where $p=3$.
 Abbreviating $\M_i\otimes_\Z\ZZ$ as $\M_i/3$
  there is a six-term exact sequence
  \small
  \[
\begin{tikzpicture}[descr/.style={fill=white,inner sep=1.5pt}]
        \matrix (m) [
            matrix of math nodes,
            row sep=3em,
            column sep=2.3em,
            text height=1.5ex, text depth=0.25ex
        ]
        { H^{n-1}(\G,\M_2/3)\oplus H^{n-1}(\G,\M_4/3) & 
        H^{n-1}(\G,\M_3/3) & H^{n-1}(\G,\M_4/3) &  \\ 
        H^{n}(\G,\M_1/3) & H^{n}(\G,\M_2/3)& 
        H^{n}(\G,\M_1/3)\oplus H^{n}(\G,\M_3/3) \\ 
        };

        \path[overlay, font=\scriptsize,>=latex]
        (m-1-1) edge [->]  node[descr,yshift=1.5ex] {$d_2^*\! +\! h_3^*$} (m-1-2) 
        (m-1-2) edge [->]  node[descr,yshift=1.5ex] {$d_3^*$} (m-1-3) 
        (m-2-1) edge [->]  node[descr,yshift=1.5ex] {$d_1^*$} (m-2-2) 
        (m-2-2) edge [->]  node[descr,yshift=1.5ex] {$h_1^*\! \oplus\! d_2^*$} (m-2-3) 
        (m-1-3) edge[out=355,in=175,looseness = .8,->] node[descr,yshift=0.3ex] {$\eta$} (m-2-1)
        
;
\end{tikzpicture} 
\] 
\normalsize
where the maps are induced by the tools of Positselski in [P].
}

\bs

\noi {\bf Proof.} This follows from  [P, Th. 6] provided
the proof is read carefully to see that the hypothesis that the
vanishing of the Bockstein ${\mathcal B}_4$ for $\M_4$ can be replaced
by $h_3\circ\,{\mathcal B}_4=0$, which is true by Theorem 4.6. $\Box$

\bs

Applying Theorem 4.7 to the field case gives the following
result. 

\bs

\noi {\bf Theorem 4.8.} {\em Suppose that
$F$ is a field containing a primitive $9$th root of unity and that
$E=F(\sqrt[3]{b_1},\sqrt[3]{b_2})$ is a $\ZZ\oplus\ZZ$
extension of $F$.
Suppose that $\M_1$, $\M_2$, $\M_3$ and $\M_4$
 are defined as in Section 1 where $p=3$.
 Abbreviating $\M_i\otimes_\Z\ZZ$ as $\M_i/3$
  there is a six-term exact sequence
  \small
  \[
\begin{tikzpicture}[descr/.style={fill=white,inner sep=1.5pt}]
        \matrix (m) [
            matrix of math nodes,
            row sep=3em,
            column sep=2.3em,
            text height=1.5ex, text depth=0.25ex
        ]
        { H^{n-1}(F,\M_2/3)\oplus H^{n-1}(F,\M_4/3) & 
        H^{n-1}(F,\M_3/3) & H^{n-1}(F,\M_4/3) &  \\ 
        H^{n}(F,\mu_3) & H^{n}(F,\M_2/3)& 
        H^{n}(F,\mu_3)\oplus H^{n}(F,\M_3/3) \\ 
        };

        \path[overlay, font=\scriptsize,>=latex]
        (m-1-1) edge [->]  node[descr,yshift=1.5ex] {$d_2^*\! +\! h_3^*$} (m-1-2) 
        (m-1-2) edge [->]  node[descr,yshift=1.5ex] {$d_3^*$} (m-1-3) 
        (m-2-1) edge [->]  node[descr,yshift=1.5ex] {$d_1^*$} (m-2-2) 
        (m-2-2) edge [->]  node[descr,yshift=1.5ex] {$h_1^* \oplus d_2^*$} (m-2-3) 
        (m-1-3) edge[out=355,in=175,->] node[descr,yshift=0.3ex] {$\eta$} (m-2-1)
        
;
\end{tikzpicture} 
\]
\normalsize
where the maps are induced by the tools of Positselski in [P].
Here $H^{n+1}(F,\M_2/3)\cong H^{n+1}(E,\mu_3) \oplus H^{n+1}(F,\mu_3)$
and $d_1^*=i_{E/F}\oplus0$ so it can be understood as the map given by 
field extension 
$H^{n+1}(F,\mu_3)\rightarrow H^{n+1}(E,\mu_3)$.
}

\bs

\bs

\noi \S 5. {\bf The  quotient ${\rm Br}_3(E/F)/{\rm Dec}(E/F)$.}

\bs

As noted in the introduction the work of Tignol [T] shows 
that  ${\rm Br}_3(E/F)/{\rm Dec}(E/F)$ can be nontrivial.
In this section we use the tools developed above to give
an explicit description of the structure of this quotient. 
The  results given here are more general and
apply to the quotient $H^n(E/F)/{\rm Dec}^n(E/F)$ for $n\geq2$.
As in Section 3 for readability
of the computations we continue to use $\chi$ to   denote both 
 a function $\chi\in Z^{n-1}(\G,\ZZ)$ and the cohomology class
$\chi\in H^{n-1}(\G,\ZZ)$ it represents.
To begin  we give the following definitions. 

\bs

\noi {\bf Definitions 5.1.}  Given the notation in Theorem 3.14, we define $X^{n-1}\subseteq
H^{n-1}(\G,\ZZ)$ by 
\begin{eqnarray*}
X^{n-1}&:=&\left\{\chi\in H^{n-1}(\G,\ZZ)\right| \exists w_{0\,2},w_{2\,0}
\in Z^{n-1}(\G,\ZZ), \ w_{2\,1},w_{1\,2} \\&&  \ \ \ \in C^{n-1}(\G,\ZZ)  \left.{\rm with} \ \ \delta(w_{2\,1})=-jw'_{2\,0}-i\chi'  \ \ {\rm and} \ \delta(w_{1\,2})=
-iw'_{0\,2}-j\chi'
\right\} .
\end{eqnarray*} 
For $k=3,4$  we set $N_k^{n-1}:=\left\{\chi\in X^{n-1}\right|\left.[\psi_k]\smile\chi=0
\in H^n(\G,\ZZ)\right\}$ 
where $\psi_3=i+j$ and $\psi_4=i-j$. Finally we set ${\mathcal O}^{n-1}:=X
^{n-1}/\left(N^{n-1}_3+N^{n-1}_4\right)$.

\bs

The subscript notation in $w_{r\, s}$ was chosen since the Definition 5.1 conditions
are equivalent to 
$w_{2\,0}t_1^2+\chi t_1t_2+ w_{0\,2}t_2^2+w_{1\,2}t_1t_2^2+w_{2\,1}t_1^2t_2+\in
H^{n-1}(\G,\ZZ[G])$, where it could have made sense to label
$\chi$ as $w_{1\,1}$ except that we want to emphasize the role of
the cocycle $\chi$ as it relates to Sections 3 and 4. This means that
in terms of field theory, the subquotient ${\mathcal O}^{n-1}$ of $H^{n-1}(F,\ZZ)$ could
be reformulated in terms of $H^{n-1}(E,\ZZ)$ where $E/F$ is the
$3$-$3$ extension under consideration.  But we prefer this formulation
because it focuses on $\chi$ as the most relevant invariant.

As a first result we see how the conditions on $\chi\in X^{n-1}$ 
can be used to give elements
of $H^n(\G,\ZZ)$.

\bs

\noi {\bf Lemma 5.2.} {\em Suppose $\chi,w_{2\,0},w_{0\,2}\in H^{n-1}(\G,\ZZ)$, 
$\delta(w_{2\,1})=-i\chi'-jw_{2\,0}'$ and $\delta(w_{1\,2})=- j\chi'-iw_{0\,2}'.$ 
We define $\alpha\in C^n(\G,\ZZ)$ by
\begin{eqnarray*}
\alpha&=&iw_{1\,2}'+jw_{2\,1}'+\lambda^2_iw_{0\,2}'+\lambda^2_jw_{2\,0}'+ij\chi'\\
&=&i(w_{1\,2}-\chi)'+j(w_{2\,1}-\chi)'+\lambda^2_i(w_{0\,2}-\chi)'+\lambda^2_j(w_{2\,0}-\chi)'
+\lambda^2_{-(i+j)}\chi'.
\end{eqnarray*}
Then $\alpha\in Z^n(\G,\ZZ)$
and moreover, if $[c]\in H^{n-1}(\G,\M_4\otimes_\Z\ZZ)$, and if $\pi_4^*([c])=\chi$, then
we can express $\eta([c])=[\alpha]\in H^n(\G,\ZZ)$  in the above form.
}

\bs

\noi {\bf Proof.}  We have $\delta(ij\chi')=-(ij'+ji')\chi''$, 
$\delta(iw_{2\,1}')=-i(-j'\chi''-i'w_{0\,2}'')$,
$\delta(\lambda^2_iw_{0\,2}')=-ii'w_{0\,2}''$,
$\delta(jw_{2\,1}')=-j(-i'\chi''-j'w_{2\,0}'')$,
$\delta(\lambda^2_jw_{2\,0}')=-jj'w_{2\,0}''$. Adding shows
that the first expression representing $\alpha$ is a cocycle.
The second expression follows from the first using
$\lambda^2_{-(i+j)}=ij+i+j+\lambda^2_i+\lambda^2_j$. 

For the second statement, Theorem 3.12 shows if 
 $\pi_4^*([c])=\chi$ then
$$\eta([c])=[(iy_1'+\lambda^2_i\chi_3')+(jy_2'+\lambda^2_j\chi_4')
+\lambda^2_\psi\chi'].$$
where $\chi_3, \chi_4\in H^{n-1}(\G,\ZZ)$,
  $y_1, y_2\in C^{n-1}(\G,\ZZ)$
and the following equations are true
\begin{eqnarray*}
\delta(y_1)&=&-i\chi_3'+\psi\chi'\\
\delta(y_2)&=&-j\chi_4'+\psi\chi'
\end{eqnarray*}
This gives the result. $\Box$

\bs

We next give  a lemma that shows that the $N_i$ are ``norms''.

\bs

\noi {\bf Lemma 5.3.} {\em Suppose that $\chi\in N^{n-1}_3$. Then there exist $n_{3\,0}$,  $n_{3\,1}$, $n_{3\,2}\in C^{n-1}(\G,\ZZ)$ with 
$n=(n_{3\,0}+n_{3\,1}t_1+n_{3\,2}t^2_1)T_3\in H^{n-1}(\G,\ZZ[G]\cdot T_3)$ and
with $\chi=n_{3\,0}$. The analogous result holds for $\chi\in N^{n-1}_4$.}

\bs

\noi {\bf Proof.} Setting $\chi=n_{3\,0}$, as $\chi\in N^{n-1}_3$, there exists
$n_{3\,1}\in C^{n-1}(\G,\ZZ)$ with $\delta(n_{3\,1})=-\psi_3n'_{3\,0}$.
By
Lemma 3.6, modifying $n_{3\,1}$ by an element of $H^{n-1}(\G,\ZZ)$  as needed, there exists $n_{3\,2}$ with $\delta(n_{3\,2})=-\psi_3n'_{3\,1}
-\lambda^2_{\psi_3}n'_{3\,0}$.  It follows that
$nT_3=(n_{3\,0}+n_{3\,1}t_1+n_{3\,2}t^2_1)T_3\in H^{n-1}(\G,\ZZ[G]T_3)$
as required. The case of $\chi\in N_4^{n-1}$ is entirely analogous. $\Box$

\bs

\noi {\bf Remark 5.4.} The previous lemma gives a computational
version of Shapiro's lemma to produce an element of $H^{n-1}(\G,\ZZ[G]\cdot T_3)$
with ``norm'' $\chi$ when interpreted in the case of field theory.

\bs

Now a lemma that allows us to utilize  the map $\pi_4$.

\bs

\noi {\bf Lemma 5.5.} {Suppose $\alpha\in {\rm ker}(H^n(\G,\ZZ)\rightarrow
H^n(\G,\M_2))$. The map given by $\alpha\mapsto\overline{\pi}_4(\alpha):=
[{\pi_4^*([c])}]\in {\mathcal O}^{n-1}=X^{n-1}/(N^{n-1}_3+N^{n-1}_4)$ where $\alpha = \eta([c])$ for
$c\in Z^{n-1}(\G,\M_4)$ is well-defined. 
}

\bs

\noi {\bf Proof.} Since $\M_4={\rm cok}(\M_2\rightarrow \M_3)$ we 
  represent $[c]$ by $c\in C^{n-1}(\G,\M_3)$. 
By the Positselski sequence in Theorem 4.7  
${\rm ker}(\eta)={\rm im}(d_3^*)$,
so that if $\eta([c'])=\eta([c])$ then $c-c'\in
Z^{n-1}(\G,\M_3)$. Writing $c-c'=(c_1,c_2,c_3T_3,c_4T_4)$ so  
$c_3T_3\in H^{n-1}(\G,\ZZ[G]T_3)$, $c_4T_4\in H^{n-1}(\G,\ZZ[G]\cdot T_4)$  we find that 
$\pi^*_4([c-c'])=[(c_3-c_4)T_G]\in H^{n-1}(\G,\ZZ[G]\cdot T_G)$.
By Lemma 5.3 we can express $c_3T_3=(c_{3\,0}+c_{3\,1}t_1+c_{3\,2}t_1^2)\cdot T_3$
and $c_4T_4=(c_{4\,0}+c_{4\,1}t_1+c_{4\,2}t_1^2)T_4$ with 
$\delta(c_3T_3)=0$ and $\delta(c_4T_4)=0$. But then we have
$\delta(c_{3\,0})=0$, $\delta(c_{3\,1})+\psi_3c_{3\,0}'=0$, 
$\delta(c_{4\,0})=0$ and $\delta(c_{4\,1})+\psi_4c_{4\,0}'=0$. From this
$ c_{3\,0},c_{4\,0}\in H^{n-1}(\G,\ZZ)$ and 
$[-\psi_3 c_{3\,0}']=[-\psi_4 c_{4\,0}']=0\in H^n(\G,\ZZ)$. These latter 
conditions show that $c_{3\,0}\in N_3$ and $c_{4\,0}\in N_4$.
The well-definition of $\overline{\pi}_4([c])\in {\mathcal O}^{n-1}$ follows.
$\Box$

\bs

We can now give the main result  about ${\mathcal O}^{n-1}$ and its relation
to Dec$^{n}$ expressed in the case of $E/F$ a $3$-$3$ extension of fields
where $F$contain a primitive $9$th root of unity.

\bs

\noi {\bf Theorem 5.6.} {\em Suppose $\alpha\in {\rm ker}(H^n(\G,\ZZ)\rightarrow
H^n(\G,\M_2))$.   Then $\eta(\alpha)\in {\rm Dec}^{n}(E/F)$ if an only if 
$\overline{\pi}_4(\alpha)=0\in {\mathcal O}^{n-1}$, that is,
$\left({\rm ker}(H^n(\G,\ZZ)\rightarrow
H^n(\G,\M_2))\right)/{\rm Dec}^n(\G,\chi_1,\chi_2)\cong{\mathcal O}^{n-1}$.
  In particular, for the $3$-$3$ extension
$E/F$,
$${\rm Br}_3(E/F)/{\rm Dec}(E/F)\cong{\mathcal O}^1 \ .$$}

\noi {\bf Proof.}  First we assume $\overline{\pi}_4^*(\alpha)=0\in {\mathcal O}^{n-1}$, that
is $\pi_4^*(\alpha)=n_{3\,0}-n_{4\,0}$ where for $i=3,4$, $n_{i\,0}\in N_i\subseteq H^{n-1}(\G,\ZZ)$. By Lemma 5.5 we have
 $n_{i\,1},n_{i\,2}\in C^{n-1}(\G,\ZZ)$ with
 $n_iT_i=(n_{i\,0}+n_{i\,1}t_1+n_{i\,2}t_1^2)T_i\in H^{n-1}(\G,\ZZ[G]T_i)$.  
It follows that $n:=(0,0,n_3T_3,n_4T_4)\in
H^{n-1}(\G,\M_3)$ and that also
$\pi_4^*([n])=n_{3\,0}-n_{4\,0}$. 
But then $\pi^*_4(\alpha-[n])=0\in H^{n-1}(\G,\ZZ)$ and
moreover, as $\delta(n_{i\,1})=-\psi_i n_{i,0}'$ we 
 see by 
Theorem 3.14 (iii) that $\eta([n])\in {\rm Dec}^{n}(E/F)$.
Next  by 
Theorem 3.14 (iii)  again, as $\pi_4^*(\alpha-[n])
=0$ we find $\eta(\alpha-[n])\in {\rm Dec}^{n}(E/F)$ and therefore
$\eta(\alpha)=\eta(\alpha-[n])+\eta([n])\in {\rm Dec}^{n}(E/F)$.

Converelsly, we assume  $\eta(\alpha)\in {\rm Dec}^{n}(E/F)$.  This means we have
$\eta(\alpha)=i\chi'_3+j\chi'_4=
ij\cdot0'+i\chi'_3+j\chi'_4+\lambda^2_i\cdot0+\lambda^2_j\cdot0\in H^2(\G,\ZZ)$ for $\chi_3,\chi_4\in H^{n-1}(\G,\ZZ)$. 
But then by computation of $\eta$ in 3.13 we have that
$[(\chi_3 t_2^2,\chi_4 t_1^2,0,0)]\in H^{n-1}(\G,\M_4\otimes_\Z\ZZ)$ with
$\eta([(\chi_3 t_2^2,\chi_4 t_1^2,0,0)])=[i\chi_3'+j\chi_4']$. Then
by the well-definition of $\overline{\pi}_4^*$,
$\overline{\pi}_4^*([\alpha])=\overline{\pi}_4^*([(\chi_3 t_2^2,\chi_4 t_1^2,0,0)])=0$.  So we see that $\overline{\pi}_4^*([\alpha])=0$ if and only if
$\eta(\alpha)\in {\rm Dec}^{n}(E/F)$.
$\Box$





  \bibliographystyle{amsalpha}

\section*{Acknolwedgement} 
The second author was supported in part by a James B. Ax Postdoctoral Fellowship.

\bs

\end{document}